\input amstex\input amsppt.sty
\magnification=\magstep1
\advance\vsize0cm\voffset=-1cm\advance\hsize1.5cm\hoffset-0.5cm
\NoBlackBoxes
\define\R{{\Bbb R}}\define\Z{{\Bbb Z}}
\def\pr{\mathop{\fam0 pr}}
\def\rk{\mathop{\fam0 rk}}
\def\forg{\mathop{\fam0 forg}}
\def\SO{\mathop{\fam0 SO}}

\def\adiag{\mathop{\fam0 adiag}}

\def\con{\mathop{\fam0 con}}
\def\id{\mathop{\fam0 id}}

\def\st{\mathop{\fam0 st}}
\def\lk{\mathop{\fam0 lk}}
\def\delet{\mathaccent"7017 }
\def\rel#1{\allowbreak\mkern8mu{\fam0rel}\,\,#1}
\def\Int{\mathop{\fam0 Int}}
\def\Cl{\mathop{\fam0 Cl}}
\def\Emb{\mathop{\fam0 Emb}}

\def\link{\mathop{\fam0 link}}
\def\im{\mathop{\fam0 im}}

\def\t{\widetilde}

\topmatter
\title Classification of embeddings below the metastable dimension \endtitle
\author A. Skopenkov \endauthor
\address Department of Differential Geometry, Faculty of Mechanics and
Mathematics, Moscow State University, Moscow, Russia 119992, and
Independent University of Moscow, B. Vlasy\-ev\-skiy, 11, 119002, Moscow, Russia.
e-mail: skopenko\@mccme.ru \endaddress
\thanks \endthanks

\abstract 
It is developed a new approach to the classical problem on 
isotopy classification of embeddings of manifolds into Euclidean spaces.
This approach involves studying of a {\it new embedding invariant}, of 
{\it almost-embeddings} and of {\it smoothing}, as well as {\it explicit 
constructions} of embeddings.
Using this approach we obtain complete concrete classification results 
{\it below the metastable dimension range}, i.e. where the 
configuration spaces invariant of Haefliger-Wu is incomplete. 
Note that all known complete concrete classification results, except for the 
Haefliger classification of links and smooth knots, can be obtained using the 
Haefliger-Wu invariant. 

More precisely, we classify embeddings $S^p\times S^{2l-1}\to\R^{3l+p}$ for 
$p<l$ in terms of homotopy groups of Stiefel manifolds 
(up to minor indeterminancies for $p>1$ and for the smooth category).
A particular case states that 
{\it the set of piecewise-linear isotopy classes of piecewise-linear embeddings
(or, equivalently, the set of almost smooth isotopy classes of smooth 
embeddings) $S^1\times S^3\to\R^7$ has a geometrically defined group structure, 
and with this group structure is isomorphic to $\Z\oplus\Z\oplus\Z_2$.} 

We exhibit {\it an example disproving the conjecture proposed by Viro and 
others} on the completeness of the multiple Haefliger-Wu invariant for 
classification of PL embeddings of connected manifolds in codimension at least 
3.
\endabstract
\endtopmatter

\document
\head 1. Main results \endhead

\smallskip
{\bf Motivation.} 

This paper is on the classical Knotting Problem: {\it for an $n$-manifold $N$ 
and a number $m$ describe isotopy classes of embeddings $N\to\R^m$}.
For recent surveys [RS99, Sk07]; whenever possible we refer to these surveys 
not to original papers. 

If $N$ is closed $d$-connected, then the Knotting Problem is easier for 
$2m\ge3n+4$ in the smooth category and for 
$$2m\ge3n+3-d$$ 
in the piecewise-linear (PL) category.  
This is so because the {\it Haefliger-Wu invariant} is bijective (by van 
Kampen, Shapiro, Wu, Haefliger, Weber and the author [RS99, \S4, Sk07, \S5]). 
This simplest configuration space invariant is defined below in \S1. 
If $N$ is closed $d$-connected, then the Knotting Problem is much harder for 
$$2m\le3n+2-d$$
both in PL and smooth categories. 
Indeed, for $N$ distinct from disjoint union of (homology) spheres no concrete 
complete description of isotopy classes was known (except recent results [KS05, 
Sk06, CRS]), in spite of the existence of interesting approaches [Br68, Wa70, 
GW99, We]. 
In this paper the case $2m=3n+2-d$ is studied, i.e. the first case when 
$2m\le3n+2-d$. 

Many interesting counterexamples in the theory of embeddings [Al24, Ko62, Hu63, 
Wa65, Ti69, BH70, Bo71, MR71, Sk02, CRS04] (see also the Whitehead Torus 
Example 1.7 below) are 
embeddings $S^p\times S^q\to\R^m$, i.e. {\it knotted tori}. 
Classification of knotted tori is a natural next step (after the link theory 
[Ha66']) 
towards classification of embeddings of {\it arbitrary} manifolds. 
Thus classification of knotted tori gives some insight or even precise 
information concerning arbitrary manifolds (cf. [Hu63, Hu69, \S12, Vr77, Sk05]) 
and reveals new interesting relations to algebraic topology. 
Since the general Knotting Problem is recognized to be unsolvable, it is 
very interesting to solve it for the important particular case of knotted tori.

\bigskip
{\bf Informal description of main results.} 

Our main results are the Main Theorems 1.3 and 1.4 below, which are 
{\it classification of embeddings $S^p\times S^{2l-1}\to\R^{3l+p}$} for $0<p<l$ 
in terms of homotopy groups of Stiefel manifolds (for $p>1$ or for the smooth 
case up to minor indeterminacy). 
An interesting particular case is Theorem 1.3.aPL for $k=1$: this is the 
first complete concrete classification of embeddings of a {\it
non-simply-connected} 4-manifold into $\R^7$, cf. [KS05, Sk05]. 
Our ideas can be extended to give rational classification for $2m<3n+2-d$ [CRS] 
or estimations for arbitrary manifolds [Sk05].

A nice feature of our classification results is an explicit and simple 
construction of representatives for all isotopy classes, see \S2.
The simplicity of these constructions suggests that such embeddings 
could well appear in other branches of mathematics. 
These constructions are used in the proof and imply interesting results such as 
Symmetry Remark in the appendix.

The main difficulty in obtaining such classification is that the Haefliger-Wu 
invariant is not complete in our situation. 
This is shown by the Whitehead Torus Example 1.7 {\it disproving the conjecture proposed by Viro} 
on the completeness of the multiple Haefliger-Wu invariant for classification 
of PL embeddings of connected manifolds in codimension at least 3 
(precise formulation is given below; the same or analogous conjectures were
proposed by Dranishnikov, Sz\"ucs, Schepin and perhaps others). 
The construction (but not the proof) of this example is very simple and 
explicit, see the end of \S1.

The main results are obtained by studying a {\it new embedding
invariant}, {\it almost-embeddings}, {\it smoothing} and {\it explicit 
construction} of all the isotopy classes of embeddings. 
The construction of the new invariant is based on known ideas (for references 
see \S2).
Thus the main point of this paper is not the invention of a new invariant but
{\it proof of new non-trivial properties of the invariant, allowing to obtain 
concrete complete classification results}.
These ideas presented in \S2 are hopefully interesting in themselves.

A difficult point is that at some points we use smoothing theory to deduce a 
result in the PL category from a result in {\it almost smooth} category, and at 
other points vice versa.  
The results in the smooth category are obtained using those in almost smooth 
category and then applying and developing smoothing theory, see \S8.

\bigskip
{\bf Plan of the paper and notation.} 

In the rest of \S1 we present statements of main results. 
Although they are independent on the previous work, we review most closely 
related known results in \S1. 
The two remaining subsections of \S1 could be read independently on each other, 
although in fact they are interrelated.
Remarks of \S1 (and of other sections) can be omitted during the first reading. 

In \S2 the main ideas are exposed (except for smoothing) and used
to deduce the main results in the PL and almost smooth categories. 
The theorems used for this deduction are proved in \S\S3--7. 
In \S8 smoothing theory is developed and applied to prove the main 
results in the smooth category.   

Sections \S3--\S8 are independent on each other in the sense that if in one 
section we use the result proved in the other, then we do not use the proof 
but only the statement which is presented in \S2 (except that in \S6 we use 
part of the Restriction Lemma 5.2). 
So a reader can pick up one of the main results of \S1 (or a specific case 
in a specific category), read its proof in \S2 and then read the proof of only 
those parts of \S\S3--8 which are required.

Let us fix some notation.  
Fix a codimension 0 ball $B$ in a smooth manifold.  
A piecewise smooth embedding of this manifold is called {\it almost smooth} if
it is a smooth embedding outside this ball. 
A piecewise smooth concordance of this manifold is called {\it almost smooth} 
if it is a smooth embedding outside $B\times I$.

Denote CAT = DIFF (smooth) or PL (piecewise linear) or AD (almost smooth). 
(The study of AD embeddings is not only interesting in itself, but is required 
to prove the main results in PL and DIFF categories.) 
If CAT is omitted, then a statement holds in all the three categories, unless 
specified otherwise.
Let $\Emb_{CAT}^m(N)$ be the
set of CAT embeddings $N\to\R^m$ up to CAT isotopy. 
For terminology concerning PL topology we refer to [RS72].

Denote by $V_{k,l}$ the Stiefel manifold of $l$-frames in $\R^k$.
Denote by $\Z_{(s)}$ the group $\Z$ for $s$ even and $\Z_2$ for $s$ odd.
This notation should not be confused with the notation for localization.
Denote $\pi^S_i=0$ for $i<0$.

\bigskip
{\bf Knotted tori.} 

Denote
$$C_q^{m-q}:=\Emb\phantom{}^m_{DIFF}(S^q)\qquad\text{and}
\qquad KT^m_{p,q,CAT}:=\Emb\phantom{}^m_{CAT}(S^p\times S^q).$$
The 'connected sum' commutative group structure on $C_q^{m-q}$ and on 
$KT^m_{0,q}$ was defined for $m\ge q+3$ in [Ha66, Ha66'].
In \S2 we define an '$S^p$-parametric connected sum' commutative group 
structure on $KT^m_{p,q,CAT}$ for $m\ge2p+q+3$ or $(m,p,q)=(7,1,3)$, cf. 
[Sk05]. 

Concrete complete classification of knotted tori was known only 
$$\text{either if }p=0\text{ and }3m\ge4q+6,
\quad\text{or if }2m\ge3q+2p+4\text{ and }p\le q.$$
We state it as the Haefliger Theorem 1.1 and Theorem 1.2. 
Our main results are the Main Theorems 1.3 and 1.4 that give {\it classification of 
knotted tori for}
$$2m=3q+2p+3\quad\text{and}\quad 1\le p<q/2.$$

{\bf The Haefliger Theorem 1.1.}

{\it (PL) The group of PL embeddings $S^q\sqcup S^q\to\R^m$ up to PL isotopy is 
$$KT^m_{0,q,PL}\cong\pi_q(S^{m-q-1})\oplus\pi_{2q+2-m}(V_{M+m-q-1,M})
\quad \text{for}\quad 3m\ge4q+6,$$ 
where $M$ is large enough [Sk07, Theorem 3.6.b, Sk07'].

(DIFF) [Sk07, Theorem 3.6.a] 
The group of smooth embeddings $S^q\sqcup S^q\to\R^m$ up to smooth 
isotopy is 
$$KT^m_{0,q,DIFF}\cong KT^m_{0,q,PL}\oplus C^{m-q}_q\oplus C^{m-q}_q\quad
\text{for}\quad m-q\ge3.$$}  

{\bf Remarks.} (a) There is a description of 
$\Emb^m(S^{p_1}\sqcup\dots\sqcup S^{p_s})$ for $m-3\ge p_1,\dots,p_s$ in terms 
of exact sequences [Ha66', Ha86]. 
For $3m<4q+6$ this description is not as explicit as for the case $3m\ge4q+6$ 
above, but it seems for the author to be 'the best possible': it 
reduces the classification of links to certain calculations in homotopy groups 
of spheres, and when such calculations cannot be completed, then no other 
approach would give an explicit answer. 
However, an alternative classification of links (e.g. using $\beta$-invarant of 
this paper [Sk07'] or isovariant configuration spaces [Me]) could give clearer 
construction of the invariants and could be applicable to related problems 
[Sk07', Me]. 

(b) Particular cases of the Haefliger Theorem 1.1PL state that 
$KT^m_{0,q,PL}\cong\pi_{2q+1-m}^S$ for $2m\ge3q+4$ [Ze62, Ha62'] 
and that $KT^{3l}_{0,2l-1,PL}\cong\pi_{2l-1}(S^l)\oplus\Z_{(l)}$ for 
$l\ge2$ [Ha62', Ha66', Theorem 10.7]. 

\smallskip
{\bf Theorem 1.2.}
{\it (PL) [Sk07, Theorem 3.9] 
The set of PL embeddings $S^p\times S^q\to\R^m$ up to PL isotopy is 
$$KT^m_{p,q,PL}=\pi_q(V_{m-q,p+1})\oplus\pi_p(V_{m-p,q+1})\quad\text{for}\quad 
p\le q\quad\text{and}\quad m\ge3q/2+p+2.$$} 
{\it (DIFF) The group of smooth embeddings $S^p\times S^q\to\R^m$ up to smooth 
isotopy is 
$$KT^m_{p,q,DIFF}
\cong\pi_q(V_{m-q,p+1})\oplus C^{m-p-q}_{p+q}\quad\text{for}\quad 
p\le q\quad\text{and}\quad m\ge\max\{2p+q+3,3q/2+p+2\}.$$} 

{\bf Remarks.} (a) Theorem 1.2DIFF was only known for $2m\ge3q+3p+4$, when 
$C^{m-p-q}_{p+q}=0$ and so $KT^m_{p,q}\cong\pi_q(V_{m-q,p+1})$ [Sk02, Corollary 
1.5] (or for $p=0$). 
The extension to $2m<3q+3p+4$ is a new result of this paper (proved in \S8) 
but not the main one. 
See [KS05, Sk06] for motivation (and explanation why Theorem 1.2DIFF is not a 
trivial corollary of Theorem 1.2PL and smoothing theory). 
Cf. [Sk06, Higher-dimensional Classsification Theorem (a)].  

(b) Theorem 1.2 is a generalization to $p>0$ of the case $2m\ge3q+4$ of 
the Haefliger Theorem 1.1.

(c) The restriction 
$m\ge2p+q+3$ is sharp in Theporem 1.2DIFF 
(because the analogue of Theorem 1.2DIFF 
for $m=2p+q+2$ is false by [Sk06, Classification 
Theorem and Higher-dimensional Classification Theorem (b)]). 

(d) For $m-3\ge2q\ge2p$ and $2m\ge3q+3p+4$ Theorem 1.2 follows by [BG71, 
Corollary 1.3]. 

(e) Cf. [LS02] for knotted tori in codimension 1.

\smallskip
{\bf Main Theorem 1.3.} 
{\it (aPL) The group of PL embeddings
$S^1\times S^{4k-1}\to\R^{6k+1}$ up to PL isotopy is 
$$KT^{6k+1}_{1,4k-1,PL}\cong\pi_{2k-2}^S\oplus\pi_{2k-1}^S\oplus\Z.$$
(aDIFF) For $k>1$ the group of smooth embeddings 
$S^1\times S^{4k-1}\to\R^{6k+1}$ up to smooth isotopy is 
$$KT^{6k+1}_{1,4k-1,DIFF}\cong
\pi_{2k-2}^S\oplus\pi_{2k-1}^S\oplus\Z\oplus G_k,$$
where $G_k$ is an abelian group of order 1, 2 or 4. 
 
(b) The group of PL or smooth embeddings $S^1\times S^{4k+1}\to\R^{6k+4}$ up to 
PL or smooth isotopy is 
$$KT^{6k+4}_{1,4k+1,PL}\cong KT^{6k+4}_{1,4k+1,DIFF}\cong 
\Z_2^a\oplus\Z_4^b\quad
\text{for some integers $a=a(k)$, $b=b(k)$ such that}$$
$$a+2b-\rk(\pi_{2k}^S\otimes\Z_2)-\rk(\pi_{2k-1}^S\otimes\Z_2)=
\cases 0         & k\in\{1,3\}\\
1                & k+1\text{ is not a power of 2} \\
1\text{ or }0    & k+1\ge8\text{ is a power of 2}\endcases.$$}

{\bf Remarks.} 
(a) More precisely, for $k>1$ we have $G_k\in\{0,\Z_2,X_k\}$, where
$X_k=\Z_2\oplus\Z_2$ for $k$ even and $X_k=\Z_4$ for $k$ odd. 

(b) For $l\ge2$ the forgetful map 
$KT^{3l+1}_{1,2l-1,AD}\to KT^{3l+1}_{1,2l-1,PL}$ is an isomorphism by 
the Almost Smoothing Theorem 2.3. 

(c) 
We conjecture that $KT^7_{1,3,DIFF}$ has a geometrically defined group 
structure and is isomorphic to 
$\Z\oplus\Z\oplus\Z_2\oplus\Z_{12}$.  

(d) In the Main Theorem 1.3.b $b=b(k)$ is not greater than each of the ranks
of $\pi_{2k}^S\otimes\Z_2$ and of $\rk(\pi_{2k-1}^S\otimes\Z_2$. 

(e) The Main Theorem 1.3.b for $k+1$ not a power of 2 can be reformulated as 
follows (corresponding reformulation exists for $k+1$ a power of 2):
{\it $KT^{6k+4}_{1,4k+1}\cong
\dfrac{\pi_{4k+1}(V_{2k+3,2})}{w_{2k+1,1}}\oplus\Z_2$, where
$w_{2k+1,1}$ is defined below.} 
Hence $KT^{6k+4}_{1,4k+1}\cong\pi_{4k+1}(V_{2k+3,2})$
if $k\not\equiv3\mod4$ (by [MS04] for $k>1$ and by \S2 for $k=1$).
We conjecture that the latter isomorphism holds for each $k$. 

In general there is an exact sequence
$0\to\pi_{2k}^S\otimes\Z_2\to \dfrac{\pi_{4k+1}(V_{2k+3,2})}{w_{2k+1,1}}
\to \pi_{2k-1}^S\otimes\Z_2\to0$. 

(f) By the Main Theorem 1.3, the above remark and [To62, Pa56, DP97] we have the 
following table, where $u^v$ means $(\Z_u)^v$ (see the details for $l=5$ in 
\S2).
$$\minCDarrowwidth{2pt}\CD l @= 2       @=3 @=4                @=5 @=6 @=7 @=8 @=9 @=10\\
KT^{3l+1}_{1,2l-1,PL}@= \Z\oplus\Z\oplus2@=4 @=\Z\oplus24\oplus2
@=2\oplus2 @=\Z @=2 @=\Z\oplus240\oplus2 @=2\t\times2^2
@=\Z\oplus2^5\endCD$$ 

In order to formulate the Main Theorem 1.4 let
$S^l\overset{\mu''}\to\to V_{l+p+1,p+1}\overset{\nu''}\to\to
V_{l+p+1,p}$ be 'forgetting the last vector' bundle. 
Denote 
$$w_{l,p}=\mu''_*[\iota_l,\iota_l]\in\pi_{2l-1}(V_{l+p+1,p+1}).$$ 
Denote by $\widehat{\Z_2}$ an unknown group which is either $0$ or
$\Z_2$, and by $\widehat\Z$ an unknown cyclic group. 
(These groups can vary with the values of $l$ and $p$.) 
By $\widehat2$ we denote the image of 2 under the quotient map 
$\Z\to\widehat\Z$. 

\smallskip
{\bf Main Theorem 1.4.} {\it Suppose that $1\le p\le l-2$.
$$(AD)\qquad KT^{3l+p}_{p,2l-1,AD}\cong\cases
\pi_{2l-1}(V_{l+p+1,p+1}) &l\in\{3,7\}\\
\dfrac{\pi_{2l-1}(V_{l+p+1,p+1})}{w_{l,p}}\oplus\Z_2&l\text{ is odd }\ne2^s-1\\
\dfrac{\pi_{2l-1}(V_{l+p+1,p+1})}{w_{l,p}}\oplus\widehat{\Z_2}& l+1=2^s\ge16\\
\dfrac{\pi_{2l-1}(V_{l+p+1,p+1})\oplus\Z}{w_{l,p}\oplus2}& l\text{ is even}
\endcases$$
$$(PL)\qquad KT^{3l+p}_{p,2l-1,PL}\cong\cases
\pi_{2l-1}(V_{l+p+1,p+1}) &l\in\{3,7\}\\
\dfrac{\pi_{2l-1}(V_{l+p+1,p+1})}{w_{l,p}}\oplus\widehat{\Z_2}&l\text{ is odd}\\
\dfrac{\pi_{2l-1}(V_{l+p+1,p+1})\oplus\widehat{\Z}}{w_{l,p}\oplus\widehat2} 
& l\text{ is even}\endcases$$
(DIFF) \quad There is an exact sequence which splits for $l$ odd
$$C^{l+1}_{2l+p-1}\overset{\zeta}\to\to KT^{3l+p}_{p,2l-1,DIFF}\overset
{\forg}\to\to KT^{3l+p}_{p,2l-1,AD}\to0.$$ 
Here the map $\zeta:C^{m-p-q}_{p+q}\to KT^m_{p,q,DIFF}$ is defined by setting
$\zeta(\varphi)$ to be the embedded connected sum of the standard torus with
the embedding $\varphi:S^{p+q}\to\R^m$.
The map $\forg$ is the obvious forgetful map.}

\smallskip
{\bf Remarks.} 
(a) Since $w_{l,p}$ has infinite order for $l$ even, by the Main Theorem 1.4AD 
and [DP97] we have 
$$\rk KT^{3l+p}_{p,2l-1,AD}=\rk\pi_{2l-1}(V_{l+p+1,p+1})=
\cases 2 & l=p+1\text{ is even}\\
1 & \text{either }l\ge p+2\text{ is even or }l=p+1\text{ is odd}\\
0  & l\text{ is odd and }l\ge p+2\endcases.$$
For $l$ odd the same relation holds also in the PL category; for $2l+p$ not
divisible by 4 also in the smooth category. Cf. [CRS]. 

(b) $KT^{3l+p}_{p,2l-1,PL}\cong KT^{3l+p}_{p,2l-1,AD}/K_{l,p}$, where 
$K_{l,p}$ is contained in $0$, in $0\oplus\Z_2$, in $0\oplus\widehat{\Z_2}$ 
and in the class of $0\oplus\Z$ for the four cases of the formula for 
$KT^{3l+p}_{p,2l-1,AD}$, respectively.

(c) In the Main Theorem 1.4 the AD category can be replaced to the category of
{\it smooth} embeddings up to {\it almost smooth} isotopy (because all 
representatives of AD embeddings are in fact smooth by the Realization Theorem 
2.4 below).

(d) The conditions that $l+1$ and $k+1$ are not powers of 2 in 
the Main Theorems 1.4, 1.3.b can be weakened to $w_{l,p}\ne0$ and 
$w_{2k+1,1}\ne0$, respectively. 
Note that $w_{l,p}=0$ for $l\in\{1,3,7,15\}$ (the equality $w_{15,p}=0$ is
due to J. Mukai).

(e) We conjecture that
$KT^{3l+p}_{p,2l-1,AD}\cong\pi_{2l-1}(V_{l+p+1,p+1})$ for $l$ distinct from
2, 4 and 8 (with possible exception of those $l$ for which $w_{l,p}=0$).
This is so either if $p=1$ and $l\equiv1,3,5\mod8$ [MS04], or if $l$ even and
the projection of $w_{l,p}$ onto the $\Z$-summand of
$\pi_{2l-1}(V_{l+p+1,p+1})$ is odd (the latter is so for $p=1$, $l$ even,
$l\ge6$, $l\ne8$; note that if the projection of $w_{l,p}$ is even, then
$KT^{3l+p}_{p,2l-1,AD}$ is a $\Z_2$-extension of $\pi_{2l-1}(V_{l+p+1,p+1})$).



(f) The Main Theorem 1.4 is a generalization to $p>0$ of the case $m=3l$, 
$q=2l-1$ of the Haefliger Theorem 1.1.  
The formulas of the Main Theorem 1.4 are false for $p=0$ by the Haefliger 
Theorem 1.1, and are true but covered by the Main Theorem 1.3 for $p=1$.   
Our proof of the Main Theorem 1.4 work with minor modifications to obtain 
a {\it new} proof of the Haefliger Theorem 1.1. 
See details in Remark (13) of \S9 and a generalization in [Sk07']. 

(g) By [Ha66, 8.15, Mi72, Theorem F and Corollary G] we have
$C^{l+1}_{2l+p-1}=0$ and so
$KT^{3l+p}_{p,2l-1,DIFF}\cong KT^{3l+p}_{p,2l-1,AD}$ if
$$\text{either}\quad p=1,\ l\equiv1\ mod2\quad \text{or}
\quad p=3,\ l\equiv\pm2\ mod12\quad \text{or}\quad p=5,\ l\equiv13,21\ mod24.$$


(h) For $1\le p=l-1$ the Main Theorem 1.4PL is true and the sequence of Main 
Theorem 1.4.DIFF is an exact sequence of sets with an action $\zeta$, and the 
map $\forg$ has a right inverse.

\bigskip
{\bf A counterexample to the Multiple Haefliger-Wu Invariant Conjecture.} 

For classification of embeddings, as well
as in other branches of mathematics, the approach using configuration spaces,
or 'complements to diagonals' proved to be very fruitful [Va92].
This approach gives the {\it Haefliger-Wu invariant} defined as follows.
Let
$$\t N=\{(x,y)\in N\times N\ |\ x\ne y\}$$
be the {\it deleted product} of $N$, i.e. the configuration space
of ordered pairs of distinct points of $N$.
For an embedding $f:N\to\R^m$ one can define a map $\t f:\t N\to S^{m-1}$ 
by the Gauss formula 
$$\t f(x,y)=\frac{fx-fy}{|fx-fy|}.$$
This map is equivariant with respect to the 'exchanging factors' involution 
$t(x,y)=(y,x)$ on~$\t N$ and antipodal involution on~$S^{m-1}$.  
The {\it Haefliger-Wu invariant} $\alpha(f)$ of the embedding $f$ is the 
equivariant homotopy class of the map $\t f$, cf. [Wu65, Gr86, 2.1.E]. 
This is clearly an isotopy invariant. 
Let $\pi^{m-1}_{eq}(\t N)$ be the set of equivariant maps $\t N\to S^{m-1}$
up to equivariant homotopy.
Thus the Haefliger-Wu invariant is a map 
$$\alpha:\Emb\phantom{}^m(N)\to\pi^{m-1}_{eq}(\t N).$$
It is important that using algebraic topology methods the set 
$\pi^{m-1}_{eq}(\t N)$ can be explicitly calculated in many cases 
[Ad93, 7.1, BG71, Ba75, RS99, Sk07].
So it is very interesting to know under which conditions the Haefliger-Wu 
invariant is bijective, and the bijectivity results we are going to state 
have many specific corollaries [BG71, Sk07]. 


\smallskip 
{\bf Classical Isotopy Theorem 1.5.}
{\it (a) [Sk07, the Haefliger-Weber Theorem 5.4] 
The Haefliger-Wu invariant is bijective for embeddings 
$N\to\R^m$ of a smooth $n$-manifold or an $n$-polyhedron $N$, if   
$$2m\ge3n+4.$$ 

(b) [Sk07, Theorem 5.5]
The Haefliger-Wu invariant is bijective for embeddings $N\to\R^m$ of a 
$d$-connected PL $n$-manifold $N$, if 
$$2m\ge3n+3-d\quad\text{and}\quad m\ge n+3.$$} 

Classical examples of Haefliger and Zeeman [Ha62', \S3, RS99, \S4, Sk02, 
Examples 1.2, Sk07, \S5] show that the Haefliger-Wu invariant is
not bijective without the above metastable dimension assumption
$2m\ge3n+4$, both in the smooth case and for non-connected $N$.

(For classical and recent examples of the {\it non-surjectivity} of
the Haefliger-Wu invariant without the restriction $2m\ge3n+3$ see [RS99, \S4, 
Sk02, Examples 1.2, MRS03, Theorem 1.3, GS06, Proposition].) 

\smallskip 
{\bf The Whitehead Link Example 1.6.} [Sk07, 3.3]
{\it For each $l$ the Haefliger-Wu invariant is not injective for 
embeddings $S^0\times S^{2l-1}\to\R^{3l}$, i.e. there exists an embedding 
$\omega_0:S^0\times S^{2l-1}\to\R^{3l}$ which is not isotopic to the standard 
embedding $f_0$ but for which $\alpha(\omega_0)=\alpha(f_0)$.}  


\smallskip 
(The construction of $\omega_0$ is presented at the end of \S1.)

On the other hand, the PL Unknotting (and Embedding) Theorems were proved first
for $2m\ge3n+4$ (or similar restriction) and then for $m\ge n+3$ [PWZ61, Ze62, 
Ir65]. 
Therefore the following conjecture was 'in the air' since 1960's (the author 
learned it from A. N. Dranishnikov and E. V. Schepin): 
{\it the Haefliger-Wu invariant is bijective for PL embeddings of connected 
PL $n$-manifolds into $\R^m$ and $m\ge n+3$}.
Even better known is the following refined version of this conjecture
(the author learned it from O. Ya. Viro and A. Sz\"ucs).

\smallskip 
{\bf Multiple Haefliger-Wu Invariant Conjecture.} {\it The multiple Haefliger-Wu invariant 
is bijective for PL embeddings of connected PL $n$-manifolds into $\R^m$ and $m\ge n+3$.}

\smallskip
Recall the definition of the multiple Haefliger-Wu invariant [Wu65] 
(for generalizations see [Wu65, VII, \S5, Va92, FM94, RS99, \S5, Sk02, \S1, 
Sk07, \S5, Me]). 
Consider the configuration space of $i$-tuples of pairwise distinct points in 
$N$:
$$\t N^i=\{(x_1,\dots,x_i)\in N^i\ |\ x_j\not=x_k\text{ if }\ j\ne k\}.$$
The space $\t N^i$ is called the {\it deleted $i$-fold product} of $N$.
The group $S_i$ of permutations of $i$ symbols obviously acts on the space 
$\t N^i$.
For an embedding $f:N\to\R^m$ define the map
$$\t f^i:\t N^i\to\t{(\R^m)}^i\quad\text{by}\quad
\t f^i(x_1,\dots,x_i)=(fx_1,\dots,fx_i).$$
Clearly, the map $\t f^i$ is $S_i$-equivariant.
Then we can define the multiple Haefliger-Wu invariant
$$\alpha_\infty=\alpha\phantom{}^m_{\infty,CAT}(N):\Emb\phantom{}^m_{CAT}(N)
\to\oplus_{i=2}^\infty[\t N^i,\t{(\R^m)}^i]_{S_i} \quad\text{by}\quad 
\alpha_\infty(f)=\oplus_{i=2}^\infty[\t f^i].$$

\smallskip
{\bf Remarks.} 
(a) {\it The connectedness assumption is essential in the Conjecture} because 
of the Whitehead Link Example 1.6 and because the example works for multiple 
Haefliger-Wu invariant (because of the following essentially known results: 
the Invariance Lemma of \S4 and the existence of an almost isotopy $\Omega_0$ from 
the Whitehead link $\omega_0$ to the standard embedding; the construction of 
$\Omega_0$ is sketched in the subsection 'almost embeddings and their 
Haefliger-Wu invariant' of \S2 and presented in \S4). 

(b) {\it The PL category is also essential in the above conjecture} because in 
the smooth case the multiple Haefliger-Wu invariant is not injective for 
$2m<3n+4$ and as highly-connected manifold as the sphere $S^n$ [Sk07, \S3].

(c) Analogously to the construction of the multiple Haefliger-Wu invariant 
but using {\it isovariant} rather than {\it equivariant} maps one can 
obtain another invariants. 
They are possibly stronger but apparently hard to calculate (at least for 
more complicated manifolds than disjoint unions of spheres; for disjoint union 
of spheres see Remark (a) to the Haefliger Theorem 1.1). 
This seems to be one of the reasons why these invariants were not mentioned in 
[Wu65] (where in the last section even more complicated generalizations were 
considered). 

(d) The Multiple Haefliger-Wu Invariant Conjecture was supported by the fact that 
analogous invariants were successfully used to classify links or link maps, and to 
construct new obstructions to embeddability [Ha62', Ma90, Ko91, HK98, Kr00].

(e) Classical Isotopy Theorem 1.5.b shows that the Multiple Haefliger-Wu 
Invariant conjecture is true 
under additional $(3n-2m+3)$-connectedness assumption, even for the ordinary (i.e. 
not multiple) Haefliger-Wu invariant. 
The inequality $2m\ge3n+2-d$ appeared in the Hudson PL version of the
Browder-Haefliger-Casson-Sullivan-Wall Embedding Theorem for closed manifolds
(proved by engulfing), but was soon proved to be superfluous (by surgery).
So it was natural to expect that the $(3n-2m+3)$-connectivity assumption in
Classical Isotopy Theorem 1.5.b is superfluous (Classical Isotopy Theorem 1.5.b 
was proved using generalization of the engulfing approach).
This information again supported the above conjecture.

(f) The restrictions like $2m\ge3n+3-d$ did appear in results for manifolds 
with boundary both for PL and smooth cases [Sk02, Theorem 1.1.$\alpha\partial$] 
(they are sharp by the Filled-Tori Theorem 7.3 below). 
But such restriction did {\it not} previously appear for closed manifolds.

(g) By the above remarks it is surprising that the assumption of 
$d$-connectedness, not only of connectedness, comes into the dimension range 
where the Haefliger-Wu invariant is complete: the connectedness assumption
is {\it essential} in Classical Isotopy Theorem 1.5.b [Sk02, Examples 1.4]. 
This shows that Classical Isotopy Theorem 1.5.b is not quite the result expected in the 1960's, 
and that its proof cannot be obtained by direct generalization of the 
Haefliger-Weber proof without invention of new ideas.

\smallskip 
Our next main result disproves the Multiple Haefliger-Wu Invariant Conjecture and 
shows that in Classical Isotopy Theorem 1.5.b the dimension restriction, or the connectivity 
assumption, is {\it sharp} (but not only {\it essential}).

\smallskip 
{\bf The Whitehead Torus Example 1.7.}  {\it If $l+1$ is not a power of 2, then the (multiple) 
Haefliger-Wu invariant is not injective for embeddings 
$S^1\times S^{2l-1}\to\R^{3l+1}$, i.e. there exists an embedding 
$\omega_1:S^1\times S^{2l-1}\to\R^{3l+1}$ which is not isotopic to the standard 
embedding $f_0$ but for which $\alpha(\omega_1)=\alpha(f_0)$ (and even 
$\alpha_\infty(\omega_1)=\alpha_\infty(f_0)$).} 

\smallskip
{\bf Remarks.} 
(a) The smooth case of the Whitehead Torus Example 1.7 is not so interesting as the PL case. 
Indeed, by smoothing theory it is known that $\alpha^{3l+p}_{DIFF}(S^{2l-1+p})$ 
is not injective for some cases, e.g. $p=1$ and $l$ even [Ha66', Mi72].

(b) The passage from the Whitehead Link Example 1.6 to the Whitehead Torus 
Example 1.7 is non-trivial for the following reasons. 
First, in the Whitehead Link Example 1.6 $l$ is arbitrary, while {\it the 
analogue of the Whitehead Torus Example 1.7 is false for $l\in\{3,7\}$} by the 
Torus Theorem 2.8 and the triviality of the $\beta$-invariant Theorem 2.9 below, 
cf. [Sk05, Isotopy Theorem]. 
Second, the Whitehead Link Example 1.6 for most $l$ requires only as simple 
an invariant as the linking coefficient, but a direct analogue of this 
invariant for tori gives only a weaker example [Sk02, Examples 1.4.i]. 


\smallskip
The construction (but not the proof!) of the Whitehead Torus Example 1.7 is as 
simple and explicit as to be given right away. 
We present here the PL version.

Proof of the Whitehead Torus Example 1.7 is based on this construction and a 
new embedding invariant (\S2,\S6). 
Formally, the Whitehead Torus Example 1.7 follows from the Whitehead 
Torus Theorem 2.5 and the Almost Smoothing Theorem 2.3 (both below).

\smallskip
{\it Construction of the PL Whitehead link 
$\omega_{0,l,PL}:S^0\times S^{2l-1}\to\R^{3l}$.} Denote coordinates in
$\R^{3l}$ by 
$(x,y,z)=(x_1,\dots,x_l,y_1,\dots,y_l,z_1,\dots,z_l)$. The
Borromean rings is the linking $S_x^{2l-1}\sqcup S_y^{2l-1}\sqcup
S_z^{2l-1}\subset\R^{3l}$ of the three spheres given by the
equations
$$\cases x=0\\ y^2+2z^2=1\endcases, \qquad
\cases y=0\\ z^2+2x^2=1\endcases \qquad\text{and}\qquad \cases
z=0\\ x^2+2y^2=1 \endcases,$$ 
respectively. 
Nowthe PL whitehead link is obtained from the Borromean rings by joining
the components $S_x^{2l-1}$ and $S_y^{2l-1}$ with a tube.

\smallskip
{\it Construction of the PL Whitehead torus
$\omega=\omega_{1,l,PL}:S^1\times S^{2l-1}\to\R^{3l+1}$.} Add a strip to the
Whitehead link $\omega_{0,l,PL}$, i.e. extend it to an embedding
$$\omega':S^0\times S^{2l-1}\bigcup\limits_{S^0\times D^{2l-1}_+=\partial
D^1_+\times D^{2l-1}_+} D^1_+\times D^{2l-1}_+\to\R^{3l}.$$ This
embedding contains connected sum of the components of the Whitehead link. 
The union of $\omega'$ and the cone over the connected sum forms an embedding 
$D^1_+\times S^{2l-1}\to\R^{3l+1}_+$. 
This latter embedding can clearly be shifted to a proper embedding. 
The PL Whitehead torus is the union of this proper embedding and its mirror image with respect to
$\R^{3l}\subset\R^{3l+1}$.

\bigskip
{\bf Acknowledgements.} 
This paper was presented by parts at the International Conference 'Knots, Links 
and Manifolds' (Siegen, Germany, 2001), the International Congress of 
Mathematicians (Beijing, China, 2002), the Kolmogorov Centennial 
Conference (Moscow, Russia, 2003) and at other events.
I would like to acknowledge P.  Akhmetiev, M. Cencelj, F. Cohen, T. Eckholm, 
U. Koschorke, M. Kreck, S. Melikhov, J. Mukai, V. Nezhinskiy, M. Skopenkov and 
A.  Sz\"ucs for useful discussions and H. Zieschang for his invitation to 
Ruhr-Universit\"at Bochum in 2001, where the first version of this paper 
was written.


\head 2. Main ideas \endhead

{\bf Notation and conventions.}

An embedding $F:N\times I\to\R^m\times I$ is a {\it concordance} if
$N\times0=F^{-1}(\R^m\times0)$ and $N\times1=F^{-1}(\R^m\times1)$.
We tacitly use the facts that in codimension at least 3

{\it concordance implies isotopy} [Hu70, Li65], and

{\it every concordance or isotopy is ambient} [Hu66, Ak69, Ed75].


Clearly, every smooth (i.e.  differentiable) map is piecewise differentiable.
The forgetful map from the set of piecewise linear embeddings (immersions) 
up to piecewise linear isotopy (regular homotopy) to the set of piecewise 
differentiable embeddings (immersions) up to piecewise differentiable isotopy 
(regular homotopy) is a 1--1 correspondence [Ha67].
Therefore we can consider any smooth map as PL one, although this is
incorrect literally.

By $[x]$ we denote the equivalence class of $x$. 
However, we often do not distinguish an element and its equivalence class.

Let $\sigma^i_k$ be the mirror-symmetry of $S^k$ or $\R^k$ with
respect to the hyperplane $x_i=0$. 
Denote by $\sigma_i$ the mirror-symmetry with respect to a hyperplane in $\R^i$ 
or with respect to an equator in $S^i$, when the concrete choice of the
hyperplane or an equator does not play a role. 
Denote by the same symbol $\sigma_i$ the map induced by this symmetry on the 
set of embeddings. 
An embedding $f:T^{p,q}\to\R^m$ is called {\it mirror-symmetric} if 
$\sigma_pf=\sigma_mf$.

Consider the standard decomposition 
$S^p=D^p_+\bigcup\limits_{\partial D^p_+=S^{p-1}=\partial D^p_-}D^p_-$. 
Analogously define $\R^m_\pm$ and $\R^{m-1}$ by equations $x_1\ge0$, $x_1\le0$ 
and $x_1=0$, respectively. 
Denote
$$T^{p,q}:=S^p\times S^q\quad\text{and}\quad T^{p,q}_\pm:=D^p_{\pm}\times S^q.$$

{\bf Construction of a commutative group structure on $KT^m_{p,q}$.}

We denote by 0 the {\it standard embedding} 
$T^{p,q}\cong S^q\times S^p\to\R^{q+1}\times\R^{p+1}\subset\R^m$.

A map $f:T^{p,q}\to S^m$ is called {\it standardized} if

$f(S^p\times D^q_-)\subset D^m_-$ is the restriction of the above standard 
embedding and 

$f(S^p\times\Int D^q_+)\subset\Int D^m_+$.

Roughly speaking, a map $T^{p,q}\to\R^m$ is standardized if its image is put 
on hyperplane $\R^{m-1}$ so that the image intersects the hyperplane by 
standardly embedded 
$S^p\times D^q_-$ (for such a map the image of $S^p\times\Int D^q_-$ can 
be pulled below the hyperplane to obtain a standardized embedding in the above 
sense). 

A concordance $F:T^{p,q}\times I\to S^m\times I$ between standardized maps is 
called {\it standardized} if 

$F(S^p\times D^q_-\times I)\subset D^m_-\times I$ is the identical concordance 
and

$F(S^p\times\Int D^q_+\times I)\subset\Int D^m_+\times I$.

\smallskip
{\bf Standardization Lemma 2.1.}
{\it Assume that 

$m\ge\max\{2p+q+2,q+3\}$ in the PL category; 

either $m\ge2p+q+3$ or $(m,p,q)=(7,1,3)$ in the AD category; 

$m\ge2p+q+3$ in the DIFF category. 

Then 

any embedding $T^{p,q}\to S^m$ is concordant to a standardized embedding, and 

any concordance between standardized embeddings $T^{p,q}\to S^m$ 
is concordant to a standardized 
embedding, and any concordance between standardized embeddings $T^{p,q}\to S^m$ 
is concordant relative to the ends to a standardized concordance.} 

\smallskip
Let $R_k:\R^k\to\R^k$ be the symmetry with respect to the $(k-2)$-hyperplane 
given by equations $x_1=x_2=0$. 

\smallskip
{\bf Group Structure Theorem 2.2.}
{\it Under the assumptions of the Standardization Lemma 2.1 for $q\ge2$ a 
commutative group 
structure on $KT^m_{p,q}$ is well-defined by the following construction.

Let $f_0,f_1:T^{p,q}\to S^m$ be two embeddings.
Take standardized embeddings $f_0',f_1'$ isotopic to them.
Let
$$(f_0+f_1)(x,y):=
\cases f_0'(x,y)& y\in D^q_+\\ R_m(f_1'(x,R_qy))& y\in
D^q_-\endcases.$$ 
Let $(-f_0)(x,y):=\sigma^2_mf_0'(x,\sigma^2_qy)$. 
Let 0 be the above standard embedding. } 

\smallskip
(The image of embedding $f_0+f_1$ is the union of $f_0'(S^p\times D^q_+)$ and
the set obtained from $f_1'(S^p\times D^q_+)$ by axial
symmetry with respect to the plane $x_1=x_2=0$.)

Under the assumptions of the Standardization Lemma 2.1 for $q=1$ we have 
$\# KT^m_{p,q}=1$. 

The Standardization Lemma 2.1 (except the case $(m,p,q)=(7,1,3)$ in the AD 
category) and the Group Structure Theorem 2.2 are proved in \S3. 
The proof is a non-trivial generalization of [Ha62, Lemma 1.3 and 1.4].

The following result (proved in \S7) reduces the case $(m,p,q)=(7,1,3)$ of the 
Standardization Lemma 2.1 in the AD category to that in the PL category. 

\smallskip
{\bf Almost Smoothing Theorem 2.3.}
{\it The forgetful map $KT^m_{1,q,AD}\to KT^m_{1,q,PL}$ is a 1--1 
correspondence for $m\ge q+4$.}

\bigskip
{\bf Explicit construction of generators of $KT^{3l+p}_{p,2l-1}$.} 

\smallskip
{\bf Realization Theorem 2.4.a for $k=1$.}
{\it For smooth embeddings $\tau^1,\tau^2,\omega:T^{1,3}\to\R^7$ constructed 
below and at the end of \S1 
$$KT^7_{1,3,PL}=KT^7_{1,3,AD}=
\left<\tau^1,\tau^2,\omega_1\ |\ 2\tau^2=2\omega_1\right>.$$}

{\it Construction of $\tau^1$ and $\tau^2$.} The maps $\tau^i$ are
defined as compositions $S^1\times S^3\overset{\pr_2\times
t^i}\to\rightarrow S^3\times S^3\subset\R^7$, where $\pr_2$ is the
projection onto the second factor, $\subset$ is the standard
inclusion and maps $t^i:S^1\times S^3\to S^3$ are defined below.
We shall see that $t^i|_{S^1\times y}$ are embeddings for
each $y\in S^3$, hence $\tau^i$ are embeddings.

Define $t^1(s,y):=sy$, where $S^3$ is identified with the set of
unit length quaternions and $S^1\subset S^3$ with the set of unit
length complex numbers. 

Define $t^2(e^{i\theta},y):=H(y)\cos\theta+\sin\theta$, where $H:S^3\to
S^2$ is the Hopf map and $S^2$ is identified with the 2-sphere formed by unit 
length quaternions of the form $ai+bj+ck$. 

\smallskip
The Realization Theorem 2.4.a for $k=1$ is a particular case of the Realization Theorem 
2.3.a below because the above construction of $\tau^1$ and $\tau^2$ is 
equivalent to the geometric construction below. 
 

\smallskip
{\bf Realization Theorem 2.4.a.} {\it Consider the map  
$$\tau^1\oplus\tau^2\oplus\omega:\pi_{2k-2}^S\oplus\pi_{4k-1}(S^{2k})\oplus\Z
\to KT^{6k+1}_{1,4k-1,DIFF},$$ 
where $\tau^1$, $\tau^2$ and $\omega$ are defined below and by 
$\omega(s):=s\omega_{1,2k,PL}$, see the end of \S1. 
The compositions of this map with the forgetful maps to $KT^{6k+1}_{1,4k-1,AD}$ 
and to $KT^{6k+1}_{1,4k-1,PL}$ are epimorphisms with the 
kernel generated by $\tau^2[\iota_l,\iota_l]-2\omega(1)$.} 

\smallskip
{\it Construction of $\tau^1$ and $\tau^2$.} 
The maps $\tau^i$ are
defined as compositions $S^1\times S^{4k-1}\overset{\pr_2\times
t^i}\to \rightarrow S^{4k-1}\times S^{2k+1}\subset\R^{6k+1}$,
where $\pr_2$ is the projection onto the second factor, $\subset$
is the standard inclusion and maps $t^i:S^1\times S^{4k-1}\to S^{4k-1}$ are
defined below. 
We shall see that $t^i|_{S^1\times y}$ are embeddings for each $y\in S^{4k-1}$, 
hence $\tau^i$ are embeddings.

Take any $\varphi\in\pi_{2k}^S\cong\pi_{4k-1}(S^{2k+1})$ realized by a smooth 
map $\varphi:S^{4k-1}\to S^{2k+1}$. Define $t^1$ on
$S^0\times S^{4k-1}$ by $t^1(\pm1,y)=\pm\varphi(y)$. Take a
non-zero vector field on $S^{2k+1}$ whose vectors at each pair of antipodal
points are opposite. Using this vector field for each $y\in
S^{4k-1}$ we can extend the map $t^1:S^0\times y\to S^{2k+1}$ to a
map $t^1:S^1\times y\to S^{2k+1}$.

Take any $\varphi\in\pi_{4k-1}(S^{2k})$ realized by a smooth map
$\varphi:S^{4k-1}\to S^{2k}$. Define $t^2$ on $S^0\times S^{4k-1}$
by $t^2(\pm1,y)=\pm\varphi(y)$. Then for each $y\in S^{4k-1}$
extend the map $t^2:S^0\times y\to S^{2k}$ as a suspension to a
map $t^2:S^1\times y\to S^{2k+1}$.

\smallskip
{\it Proof of the Main Theorem 1.3.aPL modulo Realization Theorem 2.4.a.} 
By the Realization Theorem 2.4.a it suffices to prove that 
$\dfrac{\pi_{2l-1}(V_{l+2,2})\oplus\Z}{w_{l,1}\oplus2}\cong
\pi^S_{l-2}\oplus\pi^S_{l-1}\oplus\Z$ for $l$ even.
The isomorphism is probably known, but let us prove it for completeness.
For $l$ even the bundle $\nu''$ defined before the Main Theorem 1.4 has a 
cross-section.
Hence $\pi_{2l-1}(V_{l+2,2})\cong\pi^S_{l-2}\oplus\pi_{2l-1}(S^l)$.

Clearly, $w_{l,1}$ goes under this isomorphism to $0\oplus[\iota_l,\iota_l]$. 
For $l=2,4,8$ by [To62] and the algorithm of finding a quotient of $\Z^t$ over
linear relations we have
$\dfrac{\pi_{2l-1}(S^l)\oplus\Z}{[\iota_l,\iota_l]\oplus2}\cong
\pi_{l-1}^S\oplus\Z$.
For $l\ne2,4,8$ we have $H[\iota_l,\iota_l]=\pm2$ is the generator
of $\im H$. Recall that $\pi_{2l-1}(S^l)\cong\Z\oplus(\text{a
finite group})$. By the above, the projection of
$[\iota_l,\iota_l]$ onto the $\Z$-summand is $\pm1$. Since
$\ker(\Sigma:\pi_{2l-1}(S^l)\to\pi^S_{l-1})=\left<[\iota_l,\iota_l]\right>$,
it follows that the above finite group is isomorphic to
$\pi_{l-1}^S$. Again by the algorithm of finding a quotient of
$\Z^t$ over linear relations we obtain that
$\dfrac{\pi_{2l-1}(S^l)\oplus\Z}{[\iota_l,\iota_l]\oplus2}\cong\pi_{l-1}^S\oplus\Z$.
\qed



\smallskip
{\bf Realization Theorem 2.4.b.} {\it Let CAT=PL or AD. 
For $1\le p\le l-2$ or $1=p=l-1$ there are homomorphisms
$$\tau^m_{p,q}=\tau_p=\tau:\pi_q(V_{m-q,p+1})\to KT^m_{p,q,DIFF}\quad\text{and}
\quad\omega_{p,l}=\omega_p=\omega:\Z_{(l)}\to KT^{3l+p}_{p,2l-1,DIFF}$$
such that for their compositions with the forgetful map 
$KT^m_{p,q,DIFF}\to KT^m_{p,q,CAT}$, which compositions are denoted by the same 
letters $\tau$ and $\omega$, the following holds:

$\tau^{3l+p}_{p,2l-1}$ is an isomorphism for $l\in\{3,7\}$;

$\tau^{3l+p}_{p,2l-1}\oplus\omega_{p,l}$ is an epimorphism with the kernel 
generated by $w_{l,p}\oplus(-2)$ when $l+1\ne2^s$ and CAT=AD;

$\tau^{3l+p}_{p,2l-1}\oplus\omega_{p,l}$ is an epimorphism with the kernel 
generated by $w_{l,p}\oplus(-2)$ and $0\oplus a$ for some $a\in\Z_{(l)}$ when 
either CAT=PL or $l+1=2^s\ge16$. }

\smallskip
{\it the Main Theorem 1.4} in the AD or PL category follows by 
the Realization Theorem 2.4.b in the AD or PL category, respectively. 

{\it The Realization Theorem 2.4.a} in the AD category follows by the 
case $l+1\ne2^s$ of the Realization Theorem 2.4.b in the AD category because 
$$\pi_{4k-1}(V_{2k+2,2})=\pi_{4k-1}(\partial TS^{2k+1})=
\pi_{4k-1}(S^{2k})\oplus\pi_{4k-1}(S^{2k+1})\quad\text{and so}\quad
\tau^{6k+1}_{1,4k-1}=\tau^1\oplus\tau^2.$$ 
Here $TS^{2k+1}$ is the space of vectors of length at most 1 tangent to 
$S^{2k+1}$.  

{\it The Realization Theorem 2.4.a in the PL category} follows by the 
Realization Theorem 2.4.a in the AD category and the Almost Smoothing Theorem 
2.3. 

\smallskip
{\it Proof of the Main Theorem 1.3.b in the PL category.}
By the Almost Smoothing Theorem 2.3 it suffices to prove the same result in the 
AD category.
Such a result follows by the Realization Theorem 2.4.b and the following exact 
sequence, in which $l=2k+1$ and $\Sigma\lambda''$ is multiplication by 2 [JW54].
$$\minCDarrowwidth{0pt}\CD
\pi_{2l}(S^{l+1}) @>> \lambda'' > \dfrac{\pi_{2l-1}(S^l)}{[\iota_l,\iota_l]}@>> \mu'' >
\dfrac{\pi_{2l-1}(V_{l+2,2})}{\mu''[\iota_l,\iota_l]} @>> \nu'' > 
\pi_{2l-1}(S^{l+1}) @>> \lambda'' > \pi_{2l-2}(S^l) \endCD.\quad\qed$$

{\it Proof of Remark (d) after the Main Theorem 1.3.b for $l=5$.}
We need to prove that $\pi_9(V_{7,2})/w_{5,1}\cong\Z_2$.
Consider the exact sequence [Pa56, To62]
$$\minCDarrowwidth{0pt}\CD
   \pi_9(S^5)                      @>> i > \pi_9(V_{7,2}) @>> j > \pi_9(S^6) \\
   @VV \cong V                              @VV \cong V           @VV \cong V \\
\left<[\iota_5,\iota_5]\right>=\Z_2 @>> > \Z_2\oplus\Z_2 @>> > \Z_{24} \endCD$$
Here $j$ could not be a monomorphism, so $i$ is non-zero (which also follows by [Os86])
and we are done.
\qed

\smallskip
The Realization Theorem 2.4.b is proved at the end of \S2 using the following 
construction and the results of the following subsections.  

\smallskip
{\it Construction of $\tau$.} 
We follow [Sk02, proof of Torus Lemma 6.1, Sk07, \S6]. 
Recall that $\pi_q(V_{m-q,p+1})$ is isomorphic to the
group of smooth maps $S^q\to V_{m-q,p+1}$ up to smooth homotopy. These
maps can be considered as smooth maps $\varphi:S^q\times
S^p\to\partial D^{m-q}$.  Define the smooth embedding $\tau(\varphi)$
as the composition $$S^p\times
S^q\overset{\varphi\times\pr_2}\to\to\partial D^{m-q}\times S^q
\subset D^{m-q}\times S^q\subset\R^m.$$

The map $\omega_{p,l}$ is defined by $\omega_{p,l}(s):=s\omega_{p,l,DIFF}$ 
using Whitehead Torus Theorem 2.5 and the 2-relation of the 
Relation Theorem 2.7 (both in the next subsection).

\subhead The Whitehead torus \endsubhead

\smallskip
{\bf The Whitehead Torus Theorem 2.5.} {\it For each $0\le p<l$ there is
a smooth embedding $\omega_p=\omega_{p,l,DIFF}:T^{p,2l-1}\to\R^{3l+p}$
(the Whitehead torus) such that

($\alpha$-triviality) $\alpha(\omega_p)=\alpha(0)$.

($\alpha_\infty$-triviality) $\alpha_\infty(\omega_p)=\alpha_\infty(0)$.

(non-triviality) If $l+1$ is not a power of 2, then $\omega_p$ is not 
almost smoothly isotopic to the standard embedding.} 

\smallskip
(Note that $\omega_p$ is PL isotopic to the standard embedding for
$1\le p<l\in\{3,7\}$ and is not such for $1=p<l\ne2^s-1$.) 

The PL Whitehead link $\omega_{0,l,PL}:T^{0,2l-1}\to\R^{3l}$ is constructed in 
\S1. The almost smooth Whitehead link is the same.

The smooth Whitehead link $\omega_{0,l,DIFF}:T^{0,2l-1}\to\R^{3l}$
is obtained from $\omega_{0,l,PL}$ by connected summation of the
second component of $\omega_{0,l,PL}$ with an embedding
$\varphi:S^{2l-1}\to S^{3l}$ that represents minus the element of
$C_{2l-1}^{l+1}$ obtained by connected summation of the components
of $\omega_{0,l,PL}$. 

(Note that $\omega_{0,l,PL}$ is smooth but is
not smoothly isotopic to $\omega_{0,l,DIFF}$.)

The Whitehead tori $\omega_p=\omega_{p,l,DIFF}$ are constructed inductively 
from $\omega_0=\omega_{0,l,DIFF}$ by the Extension Lemma 2.6.a below.

\smallskip
{\bf Extension Lemma 2.6.} {\it (a) Let $s\ge q+3$ and $f_0:T^{0,q}\to\R^s$ be a link.
In the smooth case assume that the connected sum of the components of $f_0$
is a smoothly trivial knot.
Then there is a sequence of embeddings $f_p:T^{p,q}\to\R^{s+p}$ such that for $p\ge1$
embedding $f_p$ is a mirror-symmetric extension of $f_{p-1}$.

(b) Suppose that $f_p,g_p:T^{p,q}\to\R^{s+p}$ are two sequences of
embeddings such that for $p\ge1$ embeddings $f_p,g_p$ are
mirror-symmetric extensions of $f_{p-1},g_{p-1}$, respectively. If
$f_0$ and $g_0$ are (almost) concordant and $p\le s-q-2$, then
$f_p$ and $g_p$ are (almost) concordant.} 

\smallskip
{\bf Relation Theorem 2.7.} {\it For $1\le p\le l-2$ or $1=p=l-1$  

(2-relation) $2\omega_p=0$ in the smooth category for $l$ odd $\ge3$.

($\tau$-relation) $\tau_p(w_{l,p})=2\omega_p$ in the almost smooth category for 
$l$ even and in the smooth category for $l$ odd.} 

\smallskip
Note that the 2-relation also holds for $p=0$.  

The properties of the Whitehead torus listed in the Whitehead Torus
Theorem 2.5 and the Relation Theorem 2.7 are proved in \S4 
(using the $\beta$-invariant Theorem 2.9 below and explicit constructions). 
Cf. Symmetry Lemma of \S4. 

\bigskip
{\bf Almost embeddings and their Haefliger-Wu invariant.} 

The {\it self-intersection set} of a map $F:N\to\R^m$ is
$$\Sigma(F)=\Cl\{x\in N\ |\ \#F^{-1}Fx\ge1\}.$$
A PL map $F:T^{p,q}\to\R^m$ is called a PL {\it almost embedding},
if $\Sigma(F)\subset T^{p,q}_-$. A map $F:T^{p,q}\to\R^m$ is
called a {\it smooth almost embedding}, if $\Sigma(F)\subset
T^{p,q}_-$ and $F|_{T^{p,q}_+}$ is a smooth embedding (i.e. a
smooth injective map such that $dF_x$ is non-degenerate at each
$x$). Almost smooth almost embedding is by definition the same as
smooth almost embedding.

A proper PL map $F:T^{p,q}\times I\to\R^m\times I$ is called a
{\it PL almost concordance} if $\Sigma(F)\subset T^{p,q}_-\times
I$. Analogously a {\it (almost) smooth almost concordance} is
defined. We say that a (almost) concordance $F$ is a (almost)
concordance {\it between} embeddings $F|_{T^{p,q}\times 0}$ and
$F|_{T^{p,q}\times 1}$ of $T^{p,q}$ into $\R^m$. The latter two
embeddings are called {\it almost concordant}. (Note that {\it
almost smooth} embedding or concordance is not the same as {\it
smooth almost} embedding or concordance.)

\smallskip
{\it Sketch of a construction of smooth almost isotopy 
$\Omega_0:T^{0,1}\times I\to\R^3\times I$ between the standard
link and the Whitehead link.} 
This interesting and well-known example introduces informally one of the main 
ingredients of our proof (see the details in the non-triviality of the 
$\beta$-invariant Theorem 2.9 below). 
Recall from \S1 that the
Whitehead link is obtained from the Borromean rings by joining two
of the rings with a tube. It is well-known that after making a
transversal self-intersection of the last two of the Borromean
rings we can move the first ring away from them. After that,
removing the self-intersection by the reverse move we obtain the
standard linking. This homotopy is not an almost isotopy. But
an almost isotopy is the corresponding homotopy $\Omega_0$ of
the Whitehead link, which link is obtained by joining with a tube those two of
the Borromean rings that intersect throughout the above homotopy.


\smallskip
Let $\overline{KT}^m_{p,q}$ be the set of almost embeddings
$T^{p,q}\to\R^m$ up to almost concordance. 
For an almost embedding $F:T^{p,q}\to\R^m$ and $m\ge2p+q+1$ we may assume by 
general position that $\Sigma(F)\cap T^{p,0}=\emptyset$.

The definition of {\it standardized almost embedding, standardized almost 
concordance} are analogous to those at the beginning of \S2. 
In Standardization Lemma 2.1 and Group Structure Theorem 2.2 (as well as in 
their proofs) we can replace 'embedding' and 'concordance' by 'almost 
embedding' and 'almost concordance'.

We have the following diagram:
$$\minCDarrowwidth{12pt}\CD
KT^m_{p,q} @>> \forg > \overline{KT}^m_{p,q} \\
@AA \tau A           @VV \bar\alpha V \\
\pi_q(V_{m-q,p+1}) @>> \alpha\tau=\bar\alpha\forg\tau >
\pi^{m-1}_{eq}(\t{T^{p,q}}) \endCD$$ 
The map $\forg$ is the obvious forgetful map.
The map $\alpha$ is the Haefliger-Wu invariant defined in \S1. 
The Haefliger-Wu invariant $\bar\alpha$ is well-defined by the Invariance Lemma 
and proof of $\alpha$-triviality of \S4, cf. [Sk02, Theorem 
5.2.$\alpha$, Sk05, Almost embeddings Theorem]. 
The map $\tau$ is defined above in \S2. 

\smallskip
{\bf Torus Theorem 2.8.} {\it Suppose that 
$$\text{either}\quad(m,p,q)=(7,1,3)\quad\text{or}\quad p\le q,\quad m\ge2p+q+3
\quad\text{and}\quad2m\ge3q+p+4.$$
The map $\bar\alpha$ is an isomorphism for 
$$2m\ge3q+2p+2\quad\text{and CAT=DIFF, \quad or }\quad 
2m\ge3q+2p+3\quad\text{and CAT=PL}.$$
The map $\alpha$ is an isomorphism for 
$$2m\ge3q+3p+4\quad\text{and CAT=DIFF, \quad or }\quad 
2m\ge3q+2p+4\quad\text{and CAT=PL}.$$
The map $\alpha\tau$ is an isomorphism for $2m\ge3q+2p+4$ and is
an epimorphism with the kernel $\left<w_{l,p}\right>$ for
$q=2l-1$, $m=3l+p$ and $p\ge1$. }


\smallskip
The assertion on $\alpha$ follows by the Classical Isotopy Theorem 1.5.a,b. 
The assertion on $\alpha\tau$ follows analogously to [HH62, Sk02, Torus 
Theorem 6.1], see the details in \S5. 
The assertion on $\bar\alpha$ in the PL category follows by [Sk02, Theorem 
2.2.q], cf. [Sk05, Almost Embeddings Theorem (b)]; there is also a direct
proof analogous to the proof of the smooth case.

{\it The new part of the Torus Theorem 2.8} is the assertion on
$\bar\alpha_{DIFF}$ for $2m\in\{3q+2p+2,3q+2p+3\}$, which is proved in \S5. 
This assertion implies the bijectivity of the forgetful homomorphism 
$\overline{KT}^m_{p,q,DIFF}\to\overline{KT}^m_{p,q,PL}$ for $2m=3q+2p+3$. 
This bijectivity does not follow from the Haefliger smoothing theory [Ha67, Ha] 
because the smoothing obstructions is in $C_q^{m-p-q}\ne0$.  

(For $2m=3q+2p+2$ the bijectivity of $\overline\alpha_{DIFF}$ is not used in 
the proof of the main results but could be used for an alternative proof of 
the the Triviality Lemma 6.2, cf. \S6.)

(Note that the Almost Smoothing Theorem 2.3 does not follow from the assertion 
on $\overline\alpha$ of the Torus Theorem 2.8: two almost smooth embeddings 
which are PL isotopic should be smoothly almost isotopic, but not necessarily 
almost smoothly isotopic.)

(Note that {\it $\alpha_{DIFF}^m(T^{p,q})$ is surjective for
$2m\ge3p+3q+3$} [Ha63]. 
The dimension assumption is sharp here because 
{\it $\alpha_{DIFF}^{3l+p+1}(T^{p,2l})$ is not surjective for 
$2\le p<l\in\{3,7\}$} (for general manifolds cf. [MRS03, Theorem 1.3]). 
Indeed, by the Non-triviality Lemma 6.3 below
there is $G\in\overline{KT}^{3l+p+1}_{p,2l,DIFF}$ such that
$\beta(G)=1$. If $\alpha(F)=\alpha(G)$ for some $F\in
KT^{3l+p+1}_{p,2l,DIFF}$, then by the assertion on
$\bar\alpha_{DIFF}$ of the Torus Theorem 2.8, $F$ is almost concordant to $G$. 
Hence $1=\beta(G)=\beta(F)=0$, which is a contradiction.)

(Note that $\forg:KT^m_{p,q,PL}\to\overline{KT}^m_{p,q,PL}$ 
is an isomorphism for $2m\ge3q+2p+4$ and is an epimorphism for $2m=3q+2p+3$.)

\bigskip
{\bf A new embedding invariant.} 

The assertion on $\overline\alpha$ of the Torus Theorem 2.8 for $2m=3q+2p+3$ 
suggests that in order to classify $KT^{3l+p}_{p,2l-1}$ we need to define 
for each almost concordance $F$ between embeddings $T^{p,2l-1}\to\R^{3l+p}$ 
an obstruction $\beta(F)$ to modification of $F$ to a concordance.

\smallskip
{\it Definition of the obstruction $\beta(F)$ for $0\le p<l$ and an almost embedding 
$F:T^{p,2l}\to\R^{3l+p+1}$.} 
Take any $x\in D^p_+$. 
Take a triangulation $T$ of $T^{p,2l}$ such that $F$ is linear on simplices of 
$T$ and $x\times S^{2l}$ is a subcomplex of $T$. 
Then $\Sigma(F)$ is a subcomplex of $T$. 
Since $F$ is an almost embedding, we have $\Sigma(F)\subset T^{p,2l}_-$. 
Denote by
$$[\Sigma(F)]\in C_{l+p-1}(T^{p,2l}_-;\Z_{(l)})$$
the sum of top-dimensional simplices of $\Sigma(F)$.

For $l$ even $[\Sigma(F)]$ is the sum of oriented simplices $\sigma$ whose 
$\pm$ coefficients are defined as follows. 
Fix in advance any orientation of $T^{p,2l}$ and of $\R^{3l+p+1}$. 
By general position there is a unique simplex $\sigma'$ of $T$ such that
$F(\sigma)=F(\sigma')$. 
The orientation on $\sigma$ induces an orientation on $F\sigma$ and then on 
$\sigma'$. 
The orientations on $\sigma$ and $\sigma'$ induce orientations on normal spaces 
in $T^{p,2l}$ to these simplices. 
The corresponding orientations in normal spaces in $F(T^{p,2l})$ to $F\sigma$ 
and to $F\sigma'$ (in this order) together with the orientation on $F\sigma$ 
induce an orientation on $\R^{3l+p+1}$. 
If the latter orientation agrees with the fixed orientation of $\R^{3l+p+1}$, then 
the coefficient of $\sigma$ is $+1$, otherwise $-1$. 
Change of orientation of $\sigma$ changes all the orientations in the above 
construction (except the fixed orientations of $T^{p,2l}$ and of $\R^{3l+p+1}$) 
and so changes the coefficient of $\sigma$ in $[\Sigma(F)]$. 
Thus the coefficient is well-defined.

(This definition of signs is equivalent to that of [Hu69, \S11, Hu70', cf. 
Sk07, \S2] given as follows. 
The orientation on $\sigma$ induces an orientation on $F\sigma$ and on 
$\sigma'$, hence it induces an orientation on their links. 
Consider the oriented $(2l+1)$-sphere $\lk F\sigma$, that is the link of 
$F\sigma$ in certain triangulation of $\R^{3l+p+1}$ 'compatible' with $T$. 
This sphere contains disjoint oriented $l$-spheres $F(\lk_T\sigma)$ and 
$F(\lk_T\sigma')$. 
Their linking coefficient 
$\link_{\lk F\sigma}(F(\lk_T\sigma),F(\lk_T\sigma'))\in\Z_{(l)}$ is the
coefficient of $\sigma$ in $[\Sigma(F)]$, which equals $\pm1$. 
For $l$ odd these signs can also be defined but are not used because 
$\partial[\Sigma(F)]=0$ only mod2.)

We have $\partial[\Sigma(F)]=0$ (only mod2 if $l$ is odd) [Hu69, Lemma 11.4,
Hu70', Lemma 1, cf. Sk07, \S2, the Whitney obstruction].
Since 
$$\Sigma(F)\subset T^{p,2l}_-\quad\text{and}\quad l+p-1<2l,\quad\text{we have}
\quad[\Sigma(F)]=\partial C\quad\text{for some}
\quad C\in C_{l+p}(T^{p,2l}_-;\Z_{(l)}).$$
We have $\partial FC=0$ [Hu70', Corollary 1.1].
Hence 
$$FC=\partial D\quad\text{for some}\quad D\in C_{l+p+1}(\R^{3l+p+1};\Z_{(l)}).$$
Since the support of $C$ is in $T^{p,2l}_-$, it follows that
$FC\cap F(x\times S^{2l})=\emptyset$.
So define 
$$\beta(F):=D\cap[F(x\times S^{2l})]\in\Z_{(l)},$$
where $\cap$ denotes algebraic intersection. 

Analogously one constructs {\it the obstruction $\beta(F)\in\Z_{(l)}$ for an 
almost concordance $F$ between embeddings $T^{p,2l-1}\to\R^{3l+p}$ almost 
concordant to the standard embedding.} 
Only replace
$$T^{p,2l}\quad\text{by}\quad T^{p,2l-1}\times I,\quad T^{p,2l}_\pm
\quad\text{by}\quad T^{p,2l-1}_\pm\times I\quad\text{and}\quad
\R^{3l+p+1}\quad\text{by}\quad \R^{3l+p}\times I.$$ 
(The condition that the 'boundary embeddings' of $F$ are almost concordant to the standard 
embedding is used not for the definition, but for the proof of the 
independence on $C$ in the case $p=l-1$, which case is only required for 
the Main Theorem 1.3.aPL for $k=1$ and Remark (h) after the Main Theorem 1.4.) 

\smallskip
{\it Proof that $\beta$-obstruction is well-defined, i.e. is
independent on the choices $C$ and $D$.} 
The independence on the choice of $D$ is standard. 

In order to prove the independence on the choice of $C$ assume that $F$ is 
an almost embedding and $\partial C_1=\partial C_2=[\Sigma(F)]$. 
Then 
$$\partial(C_1-C_2)=0\quad\text{and}
\quad C_1-C_2\in C_{l+p}(T^{p,2l}_-;\Z_{(l)}).$$
Since 
$$l+p<2l,\quad\text{we have}\quad C_1-C_2=\partial X\quad\text{for some}
\quad X\in C_{l+p+1}(T^{p,2l}_-;\Z_{(l)}).$$ 
Thus $FC_1-FC_2=\partial FX$. 
Hence we can take $D_1$ and $D_2$ so that 
$D_1-D_2=FX$ is disjoint with $[F(x\times S^{2l})]$. 
Therefore $D_1\cap[F(x\times S^{2l})]=D_2\cap[F(x\times S^{2l})]$. 

For an almost concordance $F$ the proof is analogous when $p\le l-2$.
When $p=l-1$ we have $\dim C=2l-1$.
So $C_1-C_2$ is homologous to $s[\{-x\}\times S^{2l-1}\times0]$ in $T^{p,2l-1}_-$ 
for some integer $s$.
Hence $F(C_1-C_2)$ is homologous to $sF[\{-x\}\times S^{2l-1}\times0]$.
Since $F$ is an almost concordance between embeddings almost concordant to the 
standard embedding,
it follows that $F$ extends to an almost embedding $\bar F:S^p\times D^{2l}_+\to D^{3l+p+1}_+$.
Since $C_1-C_2$ is null-homologous in $D^p_-\times D^{2l}_+$, it follows that
$F(C_1-C_2)$ is null-homologous in $D^{3l+p+1}_+-\bar F(x\times D^{2l}_+)$ and hence in
$S^{3l+p}\times I-F(x\times S^{2l-1}\times I)$.
\qed

\smallskip
The constructed $\beta$-obstruction is similar not only to [Hu70', Ha84] but 
also to the Fenn-Rolfsen-Koschorke-Kirk $\beta$-invariant of link maps [Ko88, 
\S4] and to the Sato-Levine invariant of semi-boundary links (for the 
Kirk-Livingston-Polyak-Viro-Akhmetiev generalized Sato-Levine invariant see 
[MR05]).

A map $f:T^{0,q}\to\R^m$ is an almost embedding if the images of
the components are disjoint and the restriction to the second
component is an embedding. Thus for $p=0$ the above
$\beta(F:T^{p,2l}\to\R^{3l+p+1})$ is a PL version of the Koschorke 
$\beta$-invariant [Ko88, \S4, Ko90, \S2]. 
For $p=0$ and $l=1$
the above $\beta(F:T^{p,2l-1}\times I\to\R^{3l+p}\times I)$ is the
Sato-Levine invariant. (It is independent on the choice of $C$ for
links almost concordant to the standard link,
but not for arbitrary links.)

\bigskip
{\bf Properties of the new embedding invariant.} 

\smallskip
{\bf $\beta$-invariant Theorem 2.9.} {\it Suppose that $0\le p<l$
and $F$ is an almost concordance between embeddings
$T^{p,2l-1}\to\R^{3l+p}$ almost concordant to the standard embedding.
We assume PL, AD or DIFF category, except
PL in the independence and DIFF in the completeness.

(obstruction) If $F$ is an embedding, then $\beta(F)=0$.

(invariance) $\beta(F)$ is invariant under almost concordance of
$F$ relative to the boundary.

(completeness) If $\beta(F)=0$ and $l\ge2$, then $F$ is almost
concordant relative to $T^{p,2l-1}_+\times I\cup T^{p,2l-1}\times\{0,1\}$
to a concordance.

(independence) If $l+1$ is not a power of 2, then for smooth
almost concordances $F$ between embeddings
$f_0,f_1:T^{p,2l-1}\to\R^{3l+p}$, the class $\beta(F)$ depends
only on $f_0$ and $f_1$.

(triviality) If $1\le p<l\in\{3,7\}$, then there is a concordance $F_1$ 
such that $F_1=F$ on $T^{p,2l-1}_+\times I$ and $F_1$ is a concordance between 
the same embeddings as the almost concordance $F$.

(non-triviality) There exists a smooth almost concordance
$\Omega_p$ between the Whitehead torus $\omega_p:T^{p,2l-1}\to
B^{3l+p}$ and the standard embedding such that
$\beta(\Omega_p)=1$.} 

\smallskip
The obstruction follows obviously by the definition of $\beta(F)$.
The invariance is simple and analogous to [Hu70', Lemma 2, cf. Hu69, Lemma 
11.6], see the details in [Sk05].
The completeness is a non-trivial property, but it is easily implied by
known results [Hu70', Theorem 2, Ha84, Theorem 4], see the details in \S6.

The remaining properties are new and non-trivial. 

The non-triviality is proved in \S4 using a higher-dimensional analogue of the 
above example of an almost concordance $\Omega_0$ and the Extension Lemma 2.6.b. 

The independence and the triviality are proved in \S6. 
In their proof we use an analogue of the Koschorke formula stating that for
link maps the $\beta$-invariant is the composition of the $\alpha$-invariant 
and the Hopf homomorphism. The triviality
implies that {\it the independence is false for $l\in\{3,7\}$},
thus the proof of the independence should not work for
$l\in\{3,7\}$ and so could not be easy. 
Note that the idea of the proof of a property analogous to the independence 
[HM87] does not work in our situation.

Note that the independence holds in the PL category for $p=1$ by the
Almost Smoothing Theorem 2.3. 
The condition that $l+1$ is not a power of 2 in the independence (also in its 
proof and applications) can be relaxed by $w_{l,p}\ne0$.

Let us state another simple properties of $\beta$-invariant. 
Denote by $\overline F$ the reversed $F$, i.e. $\overline
F(x,t)=F(x,1-t)$. For almost concordances $F$, $F'$ between $f_0$
and $f_1$, $f_1$ and $f_2$, respectively, define an almost
concordance $F\cup F'$ between $f_0$ and $f_2$ as 'the union' of
$$F:X\times[0,1]\to\R^m\times[0,1]\quad\text{and}\quad F':X\times[1,2]\to\R^m\times[1,2].$$
(Observe that the union of almost concordances is associative up to ambient isotopy.)

Analogously to the Group Structure Theorem 2.2 for $m\ge2p+q+3$ in 
the DIFF or AD category and for $m\ge\max\{2p+q+2,q+3\}$ in 
the PL category, we can construct a sum operation (not necessarily 
well-defined) on the set of almost concordances 
$T^{p,q}\times I\to\R^m\times I$, extending the sum operation on the boundary. 

(We can prove that the sum operation is well-defined only for $m\ge2p+q+4$ in the 
DIFF or AD category and for $m\ge2p+q+3$ in the PL category, because 
the analogue of the Standardization Lemma 2.1 for almost concordances between 
almost concordances would require such assumptions).

\smallskip
{\bf Additivity Theorem 2.10.} {\it Suppose that $0\le p<l$
and $F,F'$ are almost concordances between embeddings
$T^{p,2l-1}\to\R^{3l+p}$ almost concordant to the standard embedding.

(a) $\beta(F\cup F')=\beta(F)+\beta(F')$ and 
$\beta(\overline F)=-\beta(F)$.

(b) $\beta(F+F')=\beta(F)+\beta(F')$ (for $p=l-1$ only in the PL category).} 

\smallskip
Part (a) follows obviously by the definition of $\beta(F)$. 
Part (b) is simple, see the details in \S6.

\bigskip
{\bf Proof of the Realization Theorem 2.4.b.}

In this subsection we omit AD category from the notation (but in this
subsection this does not mean that the same holds for PL or DIFF category).
Denote $\Pi_{l,p}=\pi_{2l-1}(V_{l+p+1,p+1})$. 

\smallskip
{\it Proof of the Realization Theorem 2.4.b in the AD category.} 
By the assertion on $\alpha\tau$ of the Torus Theorem 2.8 we can identify the 
groups $\pi^{3l+p-1}_{eq}(\t{T^{p,2l-1}})\cong\Pi_{l,p}/w_{l,p}$ so that
$\alpha\tau$ is identified with
$\pr:\Pi_{l,p}\to\Pi_{l,p}/w_{l,p}$.

First assume that $l\in\{3,7\}$.
Then $\alpha\tau=\pm\id\Pi_{l,p}$, so it remains to prove that
$\pm\tau\alpha(f)=f$ for each $f\in KT^{3l+p}_{p,2l-1}$.
Since $\alpha\tau\alpha(f)=\pm\alpha(f)$, by the assertion on $\bar\alpha$ of
the Torus Theorem 2.8 there is an almost concordance $F$ between $\pm\tau\alpha(f)$
and $f$.
By the triviality of the $\beta$-invariant Theorem 2.9 we have 
$\pm\tau\alpha(f)=f$.

Now assume that $l\not\in\{3,7\}$.
Consider the following diagram.
We omit indices $p$ and $l$ from the notation of $\omega_{l,p}$ and $\tau_{l,p}$
in this proof.
$$\minCDarrowwidth{0pt}\CD
\ker\omega @>> > \Z_{(l)} @>>\omega> KT^{3l+p}_{p,2l-1} @>>\alpha> \Pi_{l,p}/w_{l,p} @>> >0\\
@ .      @VV i V                 @AA \widehat{\tau\oplus\omega} A             @AA j A @ .\\
  @ .       @>> i > \dfrac{\Pi_{l,p}\oplus\Z_{(l)}}{w_{l,p}\oplus(2\Z_{(l)}+\ker\omega)} @>> j > @ .
\endCD$$
Here $i(x):=[(0,x)]$ and $j[(a,b)]:=[a]$ (clearly, $j$ is
well-defined by this formula). 
Recall the construction of the maps $\tau$ and $\omega$ after the Realization 
Theorem 2.4.b. 
By the Relation Theorem 2.7 we can 
define the map $\widehat{\tau\oplus\omega}$ by
$[a\oplus b]\mapsto\tau(a)+b\omega$. Clearly, the left triangle is
commutative. By the $\alpha$-triviality of the Whitehead Torus
Theorem 2.5, $\alpha\omega=0$. 
This and $\alpha\tau=\pr$ imply that the right triangle is commutative. 
Since $\alpha\tau=\pr$, it follows that $\alpha$ is epimorphic, i.e.
the first line is exact at $\Pi_{l,p}/w_{l,p}$. 
The sequence is also exact at $\Z_{(l)}$. 
So if we prove that $\ker\alpha\subset\im\omega$, then the 5-lemma would 
imply that $\widehat{\tau\oplus\omega}$ is an isomorphism. 
Then the Theorem would follow because by the non-triviality of the 
Whitehead Torus Theorem 2.5 $\ker\omega=0$ for $l+1\ne2^s$. \qed


\smallskip
{\it Proof that $\ker\alpha\subset\im\omega$.} Suppose that
$\alpha(f)=0$ for some $f\in KT^{3l+p}_{p,2l-1}$. Then by the
assertion on $\bar\alpha$ of the Torus Theorem 2.8 there is an
almost concordance $F$ between $f$ and the standard embedding.

Take an almost concordance $\Omega$ given by the non-triviality of
the $\beta$-invariant Theorem 2.9. 
Let $\beta(F)\Omega$ be the sum
of $\beta(F)$ copies of $\Omega$ (recall that for $p=l-1$ the sum
is not necessarily well-defined, but we just take {\it any} of the
possible sums). Then $\beta(F)\Omega\cup\overline F$ is an almost
concordance between $\beta(F)\omega$ and $f$. 
By the Additivity Theorem 2.10 and the non-triviality of the 
$\beta$-invariant Theorem 2.9, 
$\beta(\beta(F)\Omega\cup\overline F)=\beta(F)\beta(\Omega)-\beta(F)=0$. 
Hence by the completeness of
the $\beta$-invariant Theorem 2.9, $\beta(F)\Omega\cup\overline F$
is almost concordant relative to the boundary to a concordance. So
$f=\beta(F)\omega\in\im\omega$. \qed

\smallskip
{\it Proof of the Realization Theorem 2.4.b in the PL category.} 
For $l\in\{3,7\}$ the proof is the same as in  the AD category. 
For $l\not\in\{3,7\}$ analogously to the AD category we obtain an isomorphism 
$$\widehat{\tau_{PL}\oplus\omega_{PL}}:
\frac{\Pi_{l,p}\oplus\Z_{(l)}}{w_{l,p}\oplus(2\Z_{(l)}+\ker\omega_{PL})}
\to KT^{3l+p}_{p,2l-1,PL}.\quad\qed$$

(For $p\ge2$ we do not know that $\ker\omega_{PL}=0$ for $l+1\ne2^s$. 
Note that the forgetful map $KT^{3l+p}_{p,2l-1,AD}\to KT^{3l+p}_{p,2l-1,PL}$
is surjective and its kernel is $\ker\omega_{PL}/\ker\omega_{AD}$.)

\head 3. Proof of the Group Structure Theorem 2.2 \endhead

Fix a point $y\in D^q_-\subset S^q$.

\smallskip
{\it Proof of the Standardization Lemma 2.1 for embeddings in the PL category.}
Since $m\ge2p+q+2\ge2p+2$, it follows that $g|_{S^p\times y}$ is unknotted in
$S^m$.
So there is an embedding $\hat g:D^{p+1}\subset S^m$ such that
$$(1)\qquad\hat g|_{\partial D^{p+1}}=g|_{S^p\times y}.$$
Since $m\ge2p+q+2$, by general position we may assume that
$$(2)\qquad\hat g\Int D^{p+1}\cap gT^{p,q}=\emptyset.$$
The regular neighborhood in $S^m$ of $\hat gD^{p+1}$ is homeomorphic to the
$m$-ball.
Take an isotopy moving this ball to $D^m_-$ and let $f':T^{p,q}\to S^m$ be the
PL embedding obtained from $g$.

Now we are done since the embedding $f'(S^p\times D^q_-)\subset D^m_-$ is PL 
isotopic to the standard embedding by the following result 
(because $m\ge p+q+3$, the pair $(S^p\times D^q,S^p\times S^{q-1})$
is $(q-1)$-connected and $q-1\ge 2(p+q)-m+1$).
\qed


\smallskip
{\bf Unknotting Theorem Moving the Boundary.}
{\it Let $N$ be a compact $n$-dimensional PL manifold and $f,g:N\to D^m$
proper PL embeddings.
If $m\ge n+3$ and $(N,\partial N)$ is $(2n-m+1)$-connected, then $f$ and $g$
are properly isotopic} [Hu69, Theorem 10.2]. 

\smallskip
{\it Proof of the Standardization Lemma 2.1 for embeddings in the smooth 
category.}
The bundle $\nu_{S^m}(g|_{S^p\times y})$ is stably trivial and $m-p-q\ge p$,
hence this bundle is trivial.
Take a $(m-p-q)$-framing $\xi$ of this bundle.  

Take the section formed by the first vectors of $\xi$.
As for the PL category we can 'extend' this section to a smooth embedding
$\widehat g:D^{p+1}\to\R^m$ satisfying to (1) and (2).
Take an isotopy moving the regular neighborhood in $S^m$ of $\widehat gD^{p+1}$
to $D^m_-$ and let $f':T^{p,q}\to S^m$ be the embedding obtained from $g$ by
this isotopy.

By deleting the first vector from $\xi$ we obtain a $(m-p-q-1)$-framing $\xi_1$ 
on $\widehat g(\partial D^{p+1})$ normal to $\widehat g(D^{p+1})$.  
Denote by $\eta$ the standard normal $q$-framing of 
$g(S^p\times y)$ in $gT^{p,q}$.
Then $(\xi_1,\eta)$ is a normal $(m-p-1)$-framing on 
$\widehat g(\partial D^{p+1})$ normal to $\widehat g(D^{p+1})$.  
Since $p<m-p-q-1$, the map $\pi_p(SO_{m-p-q-1})\to\pi_p(SO_{m-p-1})$ is
epimorphic.
Hence we can change $\xi_1$ (and thus $\xi$) so that $(\xi_1,\eta)$ 
extends to a normal framing on $\widehat g(D^{p+1})$.
Hence $f'(S^p\times D^q_-)\subset D^m_-$ is isotopic to the standard embedding.
\qed

\smallskip
{\it Proof of the Standardization Lemma 2.1 for concordances in the PL 
category.}
This is a relative version of the proof for embeddings with some additional 
efforts required for $m=2p+q+2$.
Take a concordance $G$ between standardized embeddings
$f_0,f_1:T^{p,q}\to S^m$.
There is a level-preserving  embedding
$\hat G:D^{p+1}\times\{0,1\}\to S^m\times\{0,1\}$ whose components satisfy to
(1) and (2).
Since $m\ge2p+q+2$ and $m\ge q+3$, it follows that $m+1\ge p+1+3$.
Hence $\hat G$ can be extended to an embedding
$\hat G:D^{p+1}\times I\subset S^m\times I$ such that
$$(1')\qquad\hat G|_{\partial D^{p+1}\times I}=G|_{S^p\times y\times I}.$$
Indeed, it suffices to see that any concordance $S^p\times I\to S^m\times I$
standard on the boundary is isotopic to the standard concordance, which follows
from the Haefliger-Zeeman Unknotting Theorem.

If $m\ge2p+q+3$, then by general position
$$(2')\qquad\hat G(\Int D^{p+1}\times I)\cap G(T^{p,q}\times I)=\emptyset.$$
If $m=2p+q+2$ (in particular, $(m,p,q)=(7,1,3)$), then $p\ge1$. 
So either $q=0$ and the Lemma is obvious or the property (2') can be achieved 
analogously to the Proposition below because $m+1=2(p+1)+q+1$, $p+1\ge2$ and 
$q\ge1$.

Take a regular neighborhood
$$B^m\times I\quad\text{in}\quad S^m\times I \quad\text{of}\quad
\hat GD^{p+1}\quad\text{such that}\quad
(B^m\times I)\cap(S^m\times\{0,1\})=D^m_-\times\{0,1\}.$$
Take an isotopy of $S^m\times I\rel S^m\times\{0,1\}$ moving $B^m\times I$ to
$D^m_-\times I$.
Let $F'$ be the concordance obtained from $G$ by this isotopy.

The embedding $F'(S^p\times D^q_-\times I)\subset D^m_-\times I$ is isotopic
$\rel D^m_-\times\{0,1\}$) to the identical concordance by the
following Unknotting Theorem Moving Part of the Boundary (which is proved
analogously to [Hu69, Theorem 10.2 on p. 199]).

{\it Let $N$ be a compact $n$-dimensional PL manifold, $A$ a codimension
zero submanifold of $\partial N$ and $f,g:N\to D^m$ proper PL embeddings.
If $m\ge n+3$ and $(N,A)$ is $(2n-m+1)$-connected, then $f$ and $g$
are properly isotopic $\rel\partial N-A$.}
\qed


\smallskip
{\bf Proposition.}
{\it If $m=2p+q+1$, $p\ge2$, $q\ge1$ and $g:T^{p,q}\to S^m$ is a PL 
embedding, then there is a PL embedding $\hat g:D^{p+1}\to S^m$ such that (1) 
and (2) hold.} 

\smallskip
{\it Proof.}
Since $m=2p+q+1\ge p+q+3$, it follows that $2m\ge3p+4$.
Hence $g(S^p\times y)$ is unknotted in $S^m$.
So there is an embedding $g':D^{p+1}\subset S^m$ such that
$g'|_{\partial D^{p+1}}=g|_{S^p\times y}$.
We shall achieve (2) by modification of $D^{p+1}$ modulo the boundary
analogously to [Hu63, proof of Lemma 1 for $M$ = a point, cf.  Sk98].

By general position $g'\Int D^{p+1}\cap gT^{p,q}$ is a finite number of
points.
These points can be joined by an arc in $g'D^{p+1}$ meeting
$\partial D^{p+1}$ only in endpoint $x$.
These points can also be joined by an arc in $gT^{p,q}$ meeting
$g(S^p\times y)$ in the endpoint $x$ only and passing the double ponts in
the same order as the previous arc.
The union of these two arcs is contained in a union
$C^2\subset S^m$ of 2-disks.
Since $m\ge p+4$ and $m\ge p+q+3$, by general position it follows that
$C^2\cap D^{p+1}$ and $C^2\cap gT^{p,q}$
are exactly the first and the second of the above arcs, respectively.
A small regular neighborhood of $C^2$ in $\R^m$ is a ball $B^m$.
Therefore

$B^m\cap gT^{p,q}$ is a proper $(p+q)$-ball in $B^m$;

$B^m\cap \partial D^{p+1}$ is a proper $p$-ball in the ball
$B^m\cap gT^{p,q}$;

$B^m\cap D^{p+1}$ is a $(p+1)$-ball in $B^m$ whose boundary
is contained in $\partial D^{p+1}\cup\partial B^m$; and

$\partial B^m\cap D^{p+1}\cap gT^{p,q}=\partial B^m\cap \partial D^{p+1}$.

Since $m\ge p+q+3$, by the Unknotting Balls Theorem there is a PL homeomorphism
$h:\con\partial B^m\to B^m$ identical on the boundary and such that
$h\con(\partial B^m\cap gT^{p,q})=B^m\cap gT^{p,q}$.
Then replacing $B^m\cap g'D^{p+1}$ with $h\con(\partial B^m\cap g'D^{p+1})$
we obtain the required embedding $\hat g$.
\qed

\smallskip
{\it The Standardization Lemma 2.1 for concordances in the smooth category}
is proved using the ideas of proof of the Standardization Lemma 2.1 
for concordances in the PL category and for embeddings in the smooth category.

{\it The Standardization Lemma 2.1 in the AD category} is proved analogously to 
the one in the smooth category (all the modifications are performed outside a 
ball on which the embedding or the concordance is not smooth). 

\smallskip
{\it Proof of the Group Structure Theorem 2.2.}
By the Standardization Lemma 2.1 for concordances, the concordance (= the 
isotopy) class of $f_0+f_1$ does not depend on the choice of $f_0'$ and $f_1'$.

Clearly, $f+0=f$. 

Let $R_{k,t}$ be the rotation of $\R^k$ whose restriction to the
plane $\R^2$ generated by $e_1$ and $e_2$ is a rotation through the angle
$\pi t$ and which leaves the orthogonal complement $\R^{k-2}$ fixed.
For each embedding $f:T^{p,q}\to\R^m$ the embeddings $f$ and
$R_{m,-t}\circ f\circ(\id_{S^p}\times R_{q,t})$ are isotopic.
The sum operation is {\it commutative} because 
$$f_0'+f_1'\quad\text{is isotopic to}
\quad R_{m,1}\circ(f_0'+f_1')\circ(\id\phantom{}_{S^p}\times R_{q,1})=f_1'+f_0'.$$ 
In order to prove {\it the associativity}, denote
$D^q_{++}=\{x\in D^q\ |\ x_1\ge0,\ x_2\ge0\}$ and define $\R^m_{++}$,
$D^q_{+-}$, $D^q_{-+}$ and $D^q_{--}$ analogously.
By the Standardization Lemma 2.1 {\it any embedding $T^{p,q}\to\R^m$ is 
isotopic to an embedding $f:T^{p,q}\to\R^m$ such that

$f:S^p\times(S^q-\Int D^q_{++})\to\R^m-\Int\R^m_{++}$ is the standard
embedding and 

$f(S^p\times\Int D^q_{++})\subset\Int\R^m_{++}$.}

If $f$, $g$ and $h$ are such embeddings, then both $f+(g+h)$ and $(f+g)+h$
are isotopic to the embedding which is standard on $S^p\times D^q_{--}$ and
whose restrictions to $S^p\times\Int D^q_{+-}$,
to $S^p\times\Int D^q_{++}$ and to $S^p\times\Int D^q_{-+}$ are obtained by
rotation from the restrictions of $f$, $g$ and $h$ onto
$S^p\times\delet D^q_{++}$.

Clearly, the embedding $-f_0'$ is standardized. 
The embedding $f_0'+(-f_0')$ can be extended to an embedding 
$S^p\times D^{q+1}_+\to\R^{m+1}_+$ by mapping linearly the segment 
$x\times[y,\sigma^1_qy]$ to the segment $fx\times[fy,\sigma^1_mfy]$.
Hence the embedding $f_0'+(-f_0')$ is isotopic to the standard one by 
the Triviality Criterion below. \qed

\smallskip
{\bf Triviality Criterion.} {\it For $m\ge\max\{2p+q+2,q+3\}$ an embedding 
$f:T^{p,q}\to S^m$ is isotopic to the standard embedding if and only if it 
extends to an embedding $\overline f:S^p\times D^{q+1}_+\to D^{m+1}_+$.}

\smallskip
{\it Proof.}
The 'only if' part is easy and does not require the dimension assumption.
In order to prove the 'if' part observe that analogously to the Standardization 
Lemma 2.1

{\it for $m\ge\max\{2p+q+2,q+3\}$ any proper embedding 
$S^p\times D^{q+1}_+\to D^{m+1}_+$ is concordant relative
to the boundary to an embedding $f:S^p\times D^{q+1}_+\to S^{m+1}$ such that

$f(S^p\times\frac12 D^{q+1}_+)\subset\frac12 D^{m+1}_+$ is the standard 
embedding and 

$f(S^p\times(D^{q+1}_+-\frac12D^{q+1}_+))\subset D^{m+1}_+-\frac12D^{m+1}_+$.} 


Apply this assertion to $\overline f$. 
Then the restriction of the obtained embedding to 
$S^p\times(D^{q+1}_+-\frac12\Int D^{q+1}_+)\to D^{m+1}_+-\frac12\Int D^{m+1}_+$
is a concordance from $f$ to the standard embedding.
\qed

\smallskip
{\bf Remarks.}
(a) {\it For $p+1=q=2k$ and $m=2p+q+2=6k$ there are no group structures on
$KT^m_{p,q,DIFF}$ such that the Whitney invariant 
$W:KT^m_{p,q,DIFF}\to\Z$ [Sk07, Sk06] is a homomorphism.} 
Indeed, $W$-preimages of distinct elements consist of different number of 
elements [Sk06, Classification Theorem and Higher-dimensional Classification 
Theorem].
(Note that $W=\alpha$.)

{\it For $p=q=2k$ and $m=2p+q+1=6k+1$ there are no group structures on
$KT^m_{p,q,DIFF}$ such that the Whitney invariant
$W:KT^m_{p,q,DIFF}\to\Z\oplus\Z$ [BH70, KS05] is a homomorphism}   
(even a sum operation which is not well-defined does not exist). 
Indeed, $W(KT^{6k+1}_{2k,2k,DIFF})=2\Z\oplus0\cup0\oplus2\Z$ [Bo71] which 
is not a subgroup of $\Z\oplus\Z$.
(Note that $BH=2W=2\alpha$.)

(b) In the Standardization Lemma 2.1 for embeddings in the DIFF and AD categories 
we only need $m\ge\max\{2p+q+2,q+3\}$. 
Hence there is a sum operation '$+$' on $KT^m_{p,q,DIFF}$ which is not 
necessarily well-defined. 
It is indeed not well-defined on $KT^{6k}_{2k-1,2k,DIFF}$ because 
for the Whitney invariant $W$ [Sk07, Sk06] we have $W(f+g)=W(f)+W(g)$ and 
$W$-preimages of distinct elements consist of different number of 
elements

(c) {\it The analogues for $m=2p+q+1$ of the Standardization Lemma 2.1
for embeddings and of the Triviality Criterion are false in the PL category}, 
in spite of the Proposition.
Indeed, a proper embedding $S^p\times D^q\subset D^{2p+q+1}$ is not
necessarily properly isotopic to the standard embedding;
counterexamples to the analogue for $m=2p+q+1\ge2q+p+2$ of the Triviality Criterion
are obtained by taking a non-trivial embedding $T^{p,q}\to \partial B^{2p+q+1}$ and
then extending them to an embedding $S^p\times D^{q+1}\to B^{2p+q+2}$ 
analogously to \S4 or to [Hu63].

(d) The Standardization Lemma 2.1 in the smooth category can also be proved
analogously to the PL category using the smooth analogue of Unknotting Theorem
Moving the Boundary (for which no additonal metastable dimension restriction is
required).
The PL case can in turn be proved analogously to the smooth case, we only
replace $SO$ by $PL$, use the Fibering Lemma 7.1 below and
$\pi_p(V^{PL}_{m-q,p+1})=0$ for $m\ge2p+q+2$ [HW66]; instead of $k$-frames
over $X$ we consider PL embeddings of $X\times D^k$.

(e) {\it An alternative proof of the if part of the Triviality Criterion in the 
PL category without use of the analogue of the Standardization Lemma 2.1.}
Take an extension $\overline f_0:S^p\times D^{q+1}_+\to B^{m+1}_+$ of the
standard embedding $f_0:T^{p,q}\to S^m$.
Since $m\ge2p+q+2$, by the Unknotting Theorem Moving the Boundary any two
proper embeddings $S^p\times D^q_-\to B^m_-$ are properly isotopic.
Hence $\overline f$ is properly isotopic to $\overline f_0$, so $f$ is
isotopic to $f_0$.

(f) Under conditions of the Standardization Lemma 2.1 for embeddings the following 
interesting 
situation occur: any proper embedding $S^p\times D^q_\pm\to D^m_\pm$ is 
isotopic to the standard embedding (moving the boundary), but an embedding 
$S^p\times S^q\to S^m$ formed by two such proper embeddings can be non-isotopic 
to the standard embedding. 

Also note that any embedding symmetric w.r.t. $S^p\times S^{q-1}$ is isotopic to 
the standard embedding, but there are embeddings symmetric w.r.t. 
$S^{p-1}\times S^q$ and non-isotopic to the standard embedding.


\head 4. The Whitehead Torus and the Whitehead almost concordance \endhead

{\bf Proof of the Extension Lemma 2.6.}

{\it Proof of (a).} Add a strip to $f_0$, i.e. extend it to an
embedding 
$$f_0':T^{0,q}\bigcup\limits_{S^0\times D^q_+=\partial D^1_+\times D^q_+} 
D^1_+\times D^q_+\to\R^s\quad\text{which is a smooth embedding on}\quad D^1_+\times D^q_+.$$ 
This embedding contains connected sum of the
components of $f_0$. The union of $f_0'$ and the cone over the
connected sum forms an embedding $D^1_+\times S^q\to\R^{s+1}_+$.
This latter embedding can clearly be shifted to a proper embedding
$f_{1,+}$.

If $f_0$ is almost smooth, then $f_{1,+}$ can be made almost
smooth. If $f_0$ is smooth, then the obstruction to smoothing the
above almost smooth $f_{1,+}$ equals to the class in $C_q^{s-q}$
of the connected sum of the components of $f_0$. So this
obstruction is zero and we may assume that $f_{1,+}$ is smooth.

Given a CAT embedding $f_{p-1,+}:T^{p-1,q}_+\to\R^{s+p-1}_+$ we rotate it in
$\R^{s+p}_+$ with respect to $\R^{s+p-2}=\partial\R^{s+p-1}_+$ to obtain
an embedding $f_{p,+}:T^{p,q}_+\to\R^{s+p}_+$.
Clearly, $f_{p,+}$ is a CAT embedding for CAT=PL, AD, DIFF.

So define inductively $f_{p,+}$ for $p\ge2$ (but not for $p=1$!)
starting with $f_{1,+}=f_+$.
Then set $f_p=\partial f_{p+1,+}$.
Clearly, $f_p$ is mirror-symmetric for $p\ge1$.
\qed



\smallskip
{\it Proof of (b).}
Let $W_0$ be a (almost) concordance between $f_0$ and $g_0$.
We have 
$$\R^{s+1}_+\cup\R^s\times I\cup\R^{s+1}_-\cong\R^{s+1}\quad\text{and}
\quad T^{1,q}_+\cup T^{0,q}\times I\cup T^{1,q}_-\cong T^{1,q}.$$ 
Let $F_0:T^{1,q}\to\R^{s+1}$ be an (almost) embedding obtained from 
$f_{1+}\cup g_{1+}\cup W_0$ by these homeomorphisms.

Now analogously to the proof of (a) extend $F_0$ to an (almost)
embedding
$$T^{1,q}\bigcup\limits_{S^1\times D^q_+=\partial D^2_+\times D^q_+}
D^2_+\times D^q_+\to\R^{s+1}\quad\text{which is a smooth embedding on}\quad D^2_+\times D^q_+$$ 
and then to a map $F_+:T^{2,q}_+\to\R^{s+2}_+$.
If $F_0$ is an embedding, then $F_+$ is also an embedding.
If $F_0$ is an almost embedding, then
$\Sigma(F_0)\cap x\times S^{q+1}=\emptyset$. Hence by the above
construction $\Sigma(F_+)\cap x\times S^{q+1}=\emptyset$, i.e.
$F_+$ is also an almost embedding.

This (almost) embedding $F_+$ can clearly be shifted to a proper (almost) concordance $W_{1+}$ between $f_{1+}$ and $g_{1+}$.
If $F_0$ is almost smooth, then $W_{1+}$ can be made almost smooth.

If $W_0$ is a smooth concordance, then the complete obstruction to smoothing
$W_{1+}$ is in  $C_{q+1}^{s-q}$ [BH70, Bo71]. If we change a
smooth concordance $W_0$ between $f_0$ and $g_0$ by connected sum
with a smooth sphere $h:S^{q+1}\to\R^s\times I$, then $W_{1+}$
changes by a connected sum with the cone over $s$. Hence the
obstruction to smoothing $W_{1+}$ changes by $[h]\in
C_{q+1}^{s-q}$ [BH70, Bo71]. Therefore by changing $W_0$ inside
its almost smooth isotopy class (and modulo $T^{1,q}$) we can make
$W_{1+}$ smooth. So we may assume that $W_{1+}$ is smooth.

In this paragraph assume that $f_p$ and $g_p$ not only are mirror-symmetric but 
also are obtained from the restrictions $f_{1+},g_{1+}:T^{1,q}_+\to\R^{s+1}_+$ 
of $f_1$ and $g_1$ as in the proof of the Extension Lemma 2.6.a.
(This simpler case is sufficient to the proof of most parts of
the Whitehead Torus Theorem 2.5.)
In this case we can define $W_p$ from $W_{1+}$ analogously to the 'rotation' construction
of $f_p$ from $f_{1+}$ in the proof of the Extension Lemma 2.6.a.

The general case is proved as follows.
Define $W_{1-}$ by symmetry to $W_{1+}$.
The two maps $W_{1+}$ and $W_{1-}$ fit together to give the required 
concordance $W_1$ between $f_1$ and $g_1$ (in the smooth category we can indeed 
make $W_1$ smooth).
So we define $W_p$ inductively using the following assertion.

{\it If $f_p,g_p:T^{p,q}\to\R^m$ are mirror-symmetric
embeddings whose restrictions $f_{p-1},g_{p-1}:T^{p-1,q}\to\R^{m-1}$ are
(almost) concordant and $m\ge2p+q+2$, then $f_p$ and $g_p$ are mirror-symmetrically
(almost) concordant.}

This assertion was proved in the first four paragraphs of this proof for $p=1$, and the proof for arbitrary $p$ is analogous.
\qed

\smallskip
Note that $f_+$ can be constructed from $f_0$ not explicitly but by
the Haefliger-Irwin-Zeeman Embedding Theorem, cf. [Hu63].
Using the corresponding Unknotting Theorem we obtain that for $p\ge1$, $m\ge2p+q+2$
and PL or AD category a map $\hat\mu'_p:KT^{m-1}_{p-1,q}\to KT^m_{p,q}$ is uniquely
defined by the condition that $\hat\mu'_p(f)$ is a mirror-symmetric extension of $f$.
Clearly, $\tau_p\hat\mu''_p=\hat\mu'_p\tau_{p-1}$,
where the homomorphism $\hat\mu''_p$ is induced by 'adding one vector' map
$V_{m-q-1,p}\to V_{m-q,p+1}$.

\bigskip
{\bf Proof of the non-triviality of $\beta$-invariant.}

{\it Construction of the Whitehead almost concordance $\Omega_p$ between
$\omega_p$ and the standard embedding.}
We start with the case $p=0$, for which we use the idea of [Ha62, \S4].
Consider the half-space $\R^{3l+1}_+$ of coordinates
$$(x,y,z,t)=(x_1,\dots,x_l,y_1,\dots,y_l,z_1,\dots,z_l,t),\quad t\ge0.$$
Take three $2l$-disks $D_x,D_y,D_z\subset\R^{3l+1}_+$ given by
the equations 
$$\cases x=0\\ y^2+2z^2+t^2=1\endcases,\qquad \cases y=0\\
z^2+2x^2+2t^2=1\endcases\qquad \text{and}\quad\cases z=0\\
x^2+2y^2+1.5t^2=1\endcases,$$ 
respectively. 
These disks bound Borromean rings.

Join $D_x$ and $D_y$ by a half-tube
$D^{2l-1}_+\times I$ such that $\partial D^{2l-1}_+\times I\subset\R^{3l}$.
We obtain a self-intersecting disk bounding the second component of the
Whitehead link. 
The union of this disk and $D_z$ form a smooth map 
$\Omega'_0:S^0\times D^{2l}\to\R^{3l+1}_+$ whose restriction
to the boundary is $\omega_0$. 

On the intersection $D_x\cap D_z$ we have $y^2+t^2=2y^2+1.5t^2=1$, hence 
$D_x\cap D_z=\emptyset$. 
Analogously $D_y\cap D_z=\emptyset$. 
So $\Sigma(\Omega'_0)\subset D_x\cap D_y$ and hence $\Omega'_0$ is an almost 
embedding. 
Deleting from $\R^{3l+1}_+$ a small regular neighborhood of an arc joining a
point from $D_z$ to a point from $D_x-D_y$ we obtain an almost concordance 
$\Omega_0$ between the Whitehead link and the standard link.

For the smooth category we modify the above $\Omega_0$ by an almost smooth isotopy 
to obtain a smooth almost concordance $\Omega_0$ between $\omega_{0,DIFF}$ and 
the standard link.

The Whitehead almost concordances $\Omega_p$ are constructed inductively by
the Extension Lemma 2.6.b starting from $\Omega_0$.
\qed

\smallskip
For the proof of the non-triviality (and thus of the Whitehead Torus Example 
1.7)
calculation of $\beta(\Omega_0)$ below can be omitted, because we
only need that $\beta(\Omega_0)=\beta(\Omega_1)\ne0$, but not that
$\beta(\Omega_0)=1$. We prove $\beta(\Omega_0)\ne0$ simpler than
$\beta(\Omega_0)=1$ as follows: if $\beta(\Omega_0)=0$, then by
the completeness of $\beta$-invariant $\omega_0$ is concordant to the standard link, which contradicts to $\lambda_{21}(\omega_0)\ne0$.

\smallskip
{\it Proof that $\beta(\Omega_p)=1$.}
First we prove the case $p=0$.
The set $D_x\cap D_y$ is an $(l-1)$-sphere in $\R^{3l+1}_+$ given by the equations
$$\cases x=y=0\\ z^2=t^2=1/3\endcases.$$
Clearly, $[\Sigma(\Omega_0)]$ is represented by $\Omega_0^{-1}(D_x\cap D_y)$.
The sphere $D_x\cap D_y$ bounds in $D_x$ and in $D_y$ two $l$-disks $C_x$ and $C_y$ given by the equations
$$\cases x=y=0\\ z^2\le1/3\\ 2z^2+t^2=1\endcases\quad\text{and}\quad
\cases x=y=0\\ z^2\le1/3\\ z^2+2t^2=1\endcases.$$
Then $\Omega _0C$ is represented by the $l$-sphere $C_x\cup C_y$.
This sphere bounds in $\R^{3l+1}$ an $(l+1)$-disk $D$ given by equations
$$\cases x=y=0\\ z^2\le1/3\\ \frac{1-z^2}2\le t^2\le1-2z^2\endcases.$$
By the Alexander-Pontryagin duality $\beta(\Omega_0)$ is the
algebraic sum of the intersection points of $D\cap D_z$.
It is easy to check that $D\cap D_z=\{(0,0,0,\sqrt{2/3})\}$.
Thus $\beta(\Omega_0)=\pm1$, and we can choose orientations so that
$\beta(\Omega_0)=1$ for the almost smooth $\Omega_0$.
Since adding an isotopy to an almost concordance does not change $\beta$-invariant, this
holds also for the smooth $\Omega_0$.

We have $\beta(\Omega_p)=\beta(\Omega_0)=1$ by the part (b) of the following
Proposition.
\qed

\smallskip
{\bf Proposition.} {\it (a) If $W:T^{p,2l}\to\R^{3l+p+1}$ is a
mirror-symmetric almost embedding and $W_0:T^{p-1,2l}\to\R^{3l+p}$
is the restriction of $W$, then $\beta(W)=\beta(W_0)$.

(b) If $W:T^{p,2l-1}\times I\to\R^{3l+p}\times I$ is a
mirror-symmetric almost concordance between embeddings and
$W_0:T^{p-1,2l-1}\times I\to\R^{3l+p-1}\times I$ is the
restriction of $W$, then $\beta(W)=\beta(W_0)$.}

\smallskip
{\it Proof.}
First we prove (a).
As in the definition of $\beta$ we take chains
$$C_0\quad\text{and}\quad C_\pm\quad\text{so that}\quad\partial C_0=[\Sigma(W_0)],
\quad\partial C_\pm=[\Sigma(W_\pm)]\pm C_0\quad\text{and}
\quad C_\pm\cap x\times S^{2l-1}=\emptyset.$$
Then $\partial(C_++C_-)=[\Sigma(W)]$.
Let $D\subset\R^{3l+p}$ be a $(2l+1)$-chain represented by a ball whose boundary is $W(x\times S^{2l})=W_0(x\times S^{2l})\subset\R^{3l+p}$.
If $l$ is even, then the chain $D$ is with integer coefficients and so the  representing disk should be oriented.
Then by the duality we have
$$\beta(W)=D\cap W(C_++C_-)=D\cap W_0(C_0)=\beta(W_0).$$

The proof of (b) is analogous.
We only take as $D\subset\R^{3l+p-1}\times I$ a $(2l+1)$-ball whose boundary
is the union of
$$W(x\times S^{2l-1}\times I)=W_0(x\times S^{2l-1}\times I)\subset\R^{m-1}\times I$$
and two $2l$-balls in $\R^{3l+p-1}\times\{0,1\}$ spanned by
$W(x\times S^{2l-1}\times i)$, $i=0,1$.
\qed

\bigskip
{\bf Proof of the Whitehead Torus Theorem 2.5.}

\smallskip
{\bf Invariance Lemma.} {\it If $f,g:N\to\R^m$ are 
embeddings of a closed $n$-manifold $N$ such that there is a codimension 0 ball 
$B\subset N$ and a homotopy $N\times I\to\R^m\times I$ between $f$ and $g$ 
whose self-intersection set is contained in $B\times I$, then 
$\alpha_\infty(f)=\alpha_\infty(g)$.} 

\smallskip
We present a proof of this folklore result for completeness. 
Cf. [Sk02, Theorem 5.2.$\alpha$, Sk05, Almost Embeddings Theorem (a)].  

\smallskip
{\it Proof.} 
It suffices to prove that for each $i$ there is an equivariant deformation 
retraction 
$$R:\t N^i\to\overline N^i:=
\{(x_1,\dots,x_i)\in \t N^i\ |\ \#(\{x_1,\dots,x_i\}\cap B)\le1\}.$$
Take $a\in N$ and a metric on $N$ such that $B$ is the closed 1-neighborhood of 
$a$ and the close 2-neighborhood $B_2$ of $a$ is a ball. 
Identify $B_2$ with the ball of radius 2 in $\R^n$. 
Now, given $(x_1,\dots,x_i)\in\t N^i-\overline N^i$, there exist $u\ne v$ 
such that $|a-x_u|$ and $|a-x_v|$ are the smallest among all $|a-x_s|$ (possibly 
$|a-x_u|=|a-x_v|$). 
Then $x_u,x_v\in B$.   
Assume that $|a-x_u|\ge|a-x_v|$. 
Then $|a-x_u|\ne0$. 
Take a piecewise-linear function 
$$M:[0,2]\to[0,2]\quad\text{such that}\quad M(0)=0,\quad M(|a-x_u|)=1
\quad\text{and}\quad M(2)=2.$$  
Define a map 
$$R:\t N^i-\overline N^i\to\overline N^i\quad\text{by}\quad R_s(x_1,\dots,x_i)
:=\cases x_sM(|a-x_s|), &x_s\in B_2\\ x_s &x_s\not\in B_2\endcases.$$
Extend this map identically over $\t N^i$. 
We obtain the required equivariant deformation retraction. 
\qed


\smallskip
{\it Proof of the $\alpha$- and the $\alpha_\infty$-triviality.} 
By the non-triviality of $\beta$-invariant there is an almost concordance 
$\Omega_p$ between $\omega_p$ and the standard embedding (only the 
existence of $\Omega_p$ is necessary here). 
We may assume that in fact $\Omega_p$ is an {\it almost isotopy} (either 
checking that the construction could give in fact an almost isotopy, cf. 
sketch of construction of $\Omega_0$ in \S2, or applying the Melikhov 'almost 
concordance implies almost isotopy in codimension at least 3' theorem). 
Since $2l+p+p+1<3l+p+1$, by general position we may assume that 
$\Sigma(\Omega_p)\cap S^p\times I=\emptyset$. 
Therefore the $\alpha$- and the $\alpha_\infty$-triviality follow by 
the Invariance Lemma.  
\qed

\smallskip
An alternative proof of the $\alpha$-triviality (without use of the Invariance 
Lemma) can be obtained observing that $\alpha$ changes as a suspension under 
the construction of the Extension Lemma 2.6.a (cf. [Sk02, Decomposition Lemma 7.1, 
Ke59, Lemma 5.1]) and proving that $\alpha(\omega_0)=0$ (either because
$\alpha(\omega_0)=\pm\Sigma\lambda(\omega_0)$ or because $\omega_0$ is almost 
concordant to the standard link).

\smallskip
{\it Proof of the non-triviality.}
Suppose to the contrary that there is an almost smooth isotopy $F$ between 
$\omega_p$ and the standard embedding.
Then by the obstruction, the independence and the non-triviality of
$\beta$-invariant we have $0=\beta(F)=\beta(\Omega_p)\ne0$, which is a
contradiction.
\qed

\bigskip
{\bf Proof of the Relation Theorem 2.7.}

\smallskip
{\bf Symmetry Lemma.} {\it For 
$$1\le p\le l-2\quad\text{or}\quad 1=p=l-1\quad\text{we have}\quad
\omega_p=\sigma_p\omega_p=\sigma_{3l+p}\omega_p=-\sigma_{2l-1}\omega_p$$
in the smooth category, except the first equality for
$l$ even, which holds only in the almost smooth category. }

\smallskip
{\it Proof.}
By construction $\omega_p$ is mirror symmetric for $p>0$, i.e.
$\sigma_p\omega_p=\sigma_{3l+p}\omega_p$.
By definition of the inverse element $\sigma_{3l+p}\omega_p=-\sigma_{2l-1}\omega_p$.

Recall that $\omega_p$ and $\sigma_p\omega_p$ are mirror-symmetric extensions
of $\omega_{p-1}$ and of $\sigma_{p-1}\omega_{p-1}$.
Hence analogously to the Extension Lemma 2.6.b (by the assertion at the end of its proof)
it suffices to prove that there exists a concordance (not necessarily mirror-symmetric)
between $\omega_1$ and $\sigma_1\omega_1$.

By the non-triviality of $\beta$-invariant $\Omega_1$ and $\sigma_1\Omega_1$
are smooth almost concordances from $\omega_1$ and from $\sigma_1\omega_1$ to the standard embedding.
Since $\sigma_1$ reverses the orientation of $S^1\times S^{2l-1}\times I$
but not of $x\times S^{2l-1}\times I$ and not of $\R^{3l+1}\times I$,
it follows that $\sigma_1$ preserves the orientation of chains
$[\Sigma(\Omega_1)]$, $C$ and $\Omega_1C$ from the definition of $\beta$.
Therefore $\beta(\sigma_1\Omega_1)=\beta(\Omega_1)$, so
$\beta(\sigma_1\Omega_1\cup\overline\Omega_1)=0$.
Hence by the completeness of $\beta$-invariant
$\sigma_1\omega_1$ is almost smoothly concordant to $\omega_1$.

If $l\ge3$ is odd, then the complete obstruction to smoothing this almost smooth concordance is in $C^{l+1}_{2l}=0$ [Ha66', Mi72], so $\sigma_1\omega_1$ is
smoothly concordant to $\omega_1$.
\qed

\smallskip
The relation $\omega_p=\sigma_p\omega_p$ (and the $\tau$-relation) in the 
smooth category holds also for $p\ge3$ and $l\equiv\pm2\ mod12$ by the above 
(and below) proof and $C_{2l+2}^{l+1}=0$ [Mi72, Corollary G].

Note that $\omega_{0,1}=\sigma_0\omega_{0,1}$ but $\omega_{0,l}\ne\sigma_0\omega_{0,l}$
for $l\not\in\{1,3,7\}$ by the Haefliger Theorem below.
It would be interesting to know if $\omega_{0,l}=\sigma_0\omega_{0,l}$ for $l\in\{3,7\}$.

\smallskip
{\bf The Haefliger Theorem.} [Ha62', \S3]
{\it For $l\not\in\{1,3,7\}$ two smooth links $f,g:T^{0,2l-1}\to\R^{3l}$ are

PL isotopic if and only if $\lambda_{12}(f)=\lambda_{12}(g)$ and
$\lambda_{21}(f)=\lambda_{21}(g)$

smoothly isotopic if and only if they are PL isotopic and their restrictions
to components are smoothly isotopic.}

\smallskip
{\it Definition of the linking coefficients $\lambda_{ij}$.}
Suppose that $p,q\le m-3$ and $f:S^p\sqcup S^q\to S^m$ is an embedding.
Take orientations of $S^p$, $S^q$, $S^m$ and $D^{m-p}$.
Take an embedding $g:D^{m-p}\to S^m$ such that $gD^{m-p}$
intersects $fS^p$ transversely at exactly one point of sign $+1$.
Then $g|_{\partial D^{m-p}}:\partial D^{m-p}\to S^m-fS^p$ is a homotopy equivalence.
Let $h_f:S^m-fS^p\to\partial D^{m-p}$ be a homotopy inverse.
Clearly, the homotopy class of $h_f$ does not depend on the choice of $g$.
  The linking coefficient is
$$\lambda_{12}(f)\ =\ [f|_{S^q}:S^q\to S^m-fS^p\overset{h_f}\to\to\partial D^{m-p}]
\ \in\ \pi_q(S^{m-p-1}).$$
Analogously is defined $\lambda_{21}(f)\in\pi_p(S^{m-q-1})$.

\smallskip
For odd $l\not\in\{1,3,7\}$ we have $2\omega_0=0$ by the Haefliger Theorem
because $2[\iota_l,\iota_l]=0$.
The following proof works for {\it each} odd $l\ge3$ and, for the almost smooth 
category, is independent on the Haefliger Theorem.

\smallskip
{\it Proof of the 2-relation.}
By the Extension Lemma 2.6.b it suffices to prove that $2\omega_0=0$.
Take a smooth almost concordance $\Omega_0$ given by the non-triviality of
$\beta$-invariant.
Then $\sigma_{3l}\sigma_{2l-1}\Omega_0\cup\overline\Omega_0$ is a smooth almost concordance
between $\sigma_{3l}\sigma_{2l-1}\omega_0=-\omega_0$ and $\omega_0$.
Since $l$ is odd, we have
$$\beta(\sigma_{3l}\sigma_{2l-1}\Omega_0\cup\overline\Omega_0)=
\pm\beta(\Omega_0)-\beta(\Omega_0)=0\in\Z_2.$$
Therefore by the completeness of $\beta$-invariant $\omega_0$ and $-\omega_0$ are
almost smoothly concordant, and hence PL isotopic.

In the smooth category the restrictions of both $\omega_0$ and $-\omega_0$ to 
their components are smoothly isotopic because
$2\varphi=0\in C_{2l-1}^{l+1}\cong\Z_2$.
Hence $\omega_0$ and $-\omega_0$ are smoothly isotopic by the Haefliger Theorem.
\qed

\smallskip
{\it Proof of the $\tau$-relation.}
For $l\in\{3,7\}$ we have $\tau(w_{l,p})=0$ because $[\iota_l,\iota_l]=0$.
So the $\tau$-relation follows from the 2-relation.

For $l\not\in\{1,3,7\}$ let us prove that

{\it (a) $\tau_0(w_{l,0})=\omega_0+\sigma_0\omega_0$ in the almost smooth category.}

Recall that $\sigma_0\psi$ is the link obtained from a link $\psi:T^{0,q}\to\R^m$ by
interchanging the components.
Recall that
$$\lambda_{12}(\omega_0)=[\iota_l,\iota_l]\quad\text{and}\quad
\lambda_{21}(\omega_0)=0,\quad\text{hence}\quad
\lambda_{12}(\sigma_0\omega_0)=0\quad\text{and}\quad
\lambda_{21}(\sigma_0\omega_0)=[\iota_l,\iota_l].$$
Since linking coefficients are additive under connected sum, it follows that
both linking coefficients of $\omega_0+\sigma_0\omega_0$ are $[\iota_l,\iota_l]$.

Identify $\pi_{2l-1}(S^l)$ and $\pi_{2l-1}(V_{l+1,1})$ by the
isomorphism $\mu''_0$ (which we thus omit from the notation).
Clearly, for each $\varphi\in\pi_{2l-1}(S^l)$ we have
$$\lambda_{12}(\tau\varphi)=\varphi\quad\text{and}\quad
\lambda_{21}(\tau\varphi)=
\lambda_{12}(\sigma_0\tau\varphi)=\lambda_{12}(\tau(a_l\circ\varphi)).$$
Here $a_l$ is the antipodal map of $S^l$.
The last equality follows because $\sigma_0\tau\varphi=\tau(a_l\circ\varphi)$.
If $l$ is odd, then $a_l\circ\varphi=\varphi$.
If $l$ is even, then
$a_l\circ\varphi=-\varphi+[\iota_l,\iota_l]\circ h_0(\varphi)$
[Po85, complement to Lecture 6, (10)].
Since $h_0[\iota_l,\iota_l]=2$ [Po85, Lecture 6, (7)], it follows that
$a_l\circ[\iota_l,\iota_l]=[\iota_l,\iota_l]$.
Therefore both linking coefficients of $\tau[\iota_l,\iota_l]$ are $[\iota_l,\iota_l]$.

So by the Haefliger Theorem $\tau_0(w_{l,0})$ is smoothly isotopic to
$\omega_0+\sigma_0\omega_0$ because the restrictions of both links to components
are standard.

Recall that $\omega_p$, $\sigma_p\omega_p$ and $\tau_p(w_{l,p})$ are mirror-symmetric
extensions of $\omega_{p-1}$, $\sigma_{p-1}\omega_{p-1}$ and $\tau_{p-1}(w_{l,p-1})$, 
respectively.
Hence by the Extension Lemma 2.6.b $\tau_p(w_{l,p})=\omega_p+\sigma_p\omega_p$
in the almost smooth category.

Recall that for $l$ odd the complete obstruction to smoothing an almost smooth concordance
$T^{1,2l-1}\times I\to\R^{3l+1}\times I$ is in $C^{l+1}_{2l}=0$ [Ha66', Mi72].
Hence $\tau_1(w_{l,1})$ is smoothly concordant to $\omega_1+\sigma_1\omega_1$ (the
concordance is not necessarily mirror-symmetric).
Hence analogously to the Extension Lemma 2.6.b
$\tau_p(w_{l,p})=\omega_p+\sigma_p\omega_p$ in the smooth category for $p\ge1$.

Now the $\tau$-relation follows by the relation $\sigma_p\omega_p=\omega_p$ of 
the Symmetry Lemma.
\qed

\smallskip
Note that the smooth version of (a) is false
(because the smooth $\omega_0$ is constructed by connected summation with $\varphi$).

\head 5. Proof of the Torus Theorem 2.8 \endhead

Denote $KT^m_{p,q,+,CAT}:=\Emb\phantom{}^m_{CAT}(D^p_+\times S^q)$. 
For $m\ge\max\{2p+q+2,q+3\}$ a group structure on $KT^m_{p,q,+}$ is defined 
analogously to that on $KT^m_{p,q}$ above. 
For the smooth category the sum is connected sum of $q$-spheres together with 
fields of $p$ normal vectors.

\smallskip
{\bf Lemma 5.1.} {\it For $p\le q$, $m\ge p+q+3$ and $m\ge 2p+q+2$
there are homomorphisms
$$\pi_q(V_{m-q,p+1})@>> \tau_+ > KT^m_{p+1,q,+} @>> \alpha >
\pi^{m-1}_{eq}(\t{T^{p+1,q}_+})@>> r > \pi^{m-1}_{eq}(\t{T^{p,q}})
@<< \sigma' < \pi_q(V^{eq}_{m-q,p+1}).$$

Here $r$ is the isomorphism induced by restriction.

The map $\alpha r$ is an isomorphism for $2m\ge3q+2p+4$ in the
smooth category and for $2m\ge3q+2p+5$ in the PL category.

The map $\tau_+$ is defined as follows (analogously to $\tau$ in \S2). 
Represent $\varphi\in\pi_q(V_{m-q,p+1})$ as a mapping 
$D^p\times S^q\to D^{m-q}$. 
Define $\tau_+(\varphi)$ to be the composition 
$D^p\times S^q\overset{\varphi\times\pr_2}\to\to D^{m-q}\times S^q\subset\R^m$. 

The map $\tau_+$ is an isomorphism for $2m\ge3q+4$ in the smooth category and 
for $2m\ge3q+2p+5$ in the PL category.

The} equivariant Stiefel manifold $V^{eq}_{mn}$ {\it is the space
of maps $S^{n-1}\to S^{m-1}$ which are equivariant with respect to
the antipodal involutions. 
The map $\sigma'$ is defined (as $\gamma^{-1}i_2\sigma$) and for $2m\ge3q+p+4$ 
is proved to be an isomorphism in 
[Sk02, Proof of Torus Lemma 6.1, Sk07, Proof of Torus Lemma 6.1].} 

\smallskip
{\it Proof.} Analogously to [Sk02, Torus Lemma 6.1, Sk07, Torus Lemma 6.1]
there is an equivariant deformation retraction
$$\t{T^{p+1,q}_+}\to\adiag(\partial D^{p+1})\times S^q\times S^q
\bigcup\limits_{\adiag(\partial D^{p+1})\times\adiag S^q}
D^{p+1}\times D^{p+1}\times\adiag S^q.$$ Hence by general position
$r$ is an isomorphism for $m\ge p+q+3$ and $m\ge 2p+q+2$.

Then the assertion on $\alpha r$ follows by [Ha63, 6.4, Sk02,
Theorems 1.1.$\alpha\partial$ and 1.3.$\alpha\partial$].

The assertion on $\tau_+$ in the PL category follows the assertion on
$\alpha r$ and the assertion on $\alpha\tau$ of the Torus Theorem 2.8.

The assertion on $\tau_+$ in the smooth category follows because for
$2m\ge3q+4$ every smooth embedding $S^q\to\R^m$ is smoothly
isotopic to the standard one and by the (trivial) smooth version
of the Fibering Lemma 7.1 below, cf. the Filled-Tori Theorem 7.3.a below.
\qed

\proclaim{Restriction Lemma 5.2} For $p\ge1$, $m\ge2p+q+2$ and $2m\ge3q+p+4$
there is the following commutative (up to sign) diagram with exact lines:
$$\minCDarrowwidth{0pt}\CD
@>> \lambda''_{q+1} > \pi_q(S^{m-p-q-1}) @>> \mu'' > \pi_q(V_{m-q,p+1}) @>> \nu'' >
\pi_q(V_{m-q,p}) @>> \lambda''_q > \pi_{q-1}(S^{m-p-q-1}) \\
@. @VV \Sigma^p V @VV \tau V @VV \tau_+ V @VV \Sigma^p V \\
@>> \overline\lambda_{q+1} > \pi_{p+q}(S^{m-q-1}) @>> \overline\mu > \overline{KT}_{p,q}^m
@>>\overline\nu> KT^m_{p,q,+} @>> \overline\lambda_q > \pi_{p+q-1}(S^{m-q-1})\\
@. @VV = V @VV \bar\alpha' V @VV \alpha' r V @VV = V \\
@>> \lambda'_{q+1} > \pi_{p+q}(S^{m-q-1})@>> \mu' > \pi_q(V^{eq}_{m-q,p+1})  @>> \nu' >
\pi_q(V^{eq}_{m-q,p}) @>> \lambda'_q > \pi_{p+q-1}(S^{m-q-1})\endCD.$$
Here the the upper and the bottom lines are the exact sequences of the 
'restriction' Serre fibrations $S^{m-p-q-1}\to V_{m-q,p+1}\to V_{m-q,p}$ and
$\Omega_p S^{m-p-q-1}\to V^{eq}_{m-q,p+1}\to  V^{eq}_{m-q,p}$.
The homomorphism $\overline\nu$ is restriction-induced.
The homomorphisms $\overline\lambda_q$ and $\overline\mu$ are defined below.

By $\bar\alpha'$ and $\alpha'r$ we denote the compositions of
$\bar\alpha$ and $\alpha r$ with the isomorphism $\sigma'$ of Lemma 5.1. 
Denote by $\rho=\rho_p:\pi_q(V_{m-q,p})\to\pi_q(V^{eq}_{m-q,p})$ the
inclusion-induced homomorphism. 
Then $\bar\alpha'\tau=\pm\rho_{p+1}$ and $\alpha'r\tau_+=\pm\rho_p$.

For $p\ge1$ and $2m=3q+p+3\ge4p+2q+4$ all the diagram except $\bar\alpha'$ is still
defined and commutative (up to sign), and the lines are exact.
\endproclaim


Recall that $T^{p,q}=S^p\times S^q$ and $T^{p,q}_+=D^p_+\times S^q$.
Denote $$B^{p+q}:=T^{p,q}-(\Int D^p_+\times S^q_+\cup\Int S^p\times D^q_+).$$ 
Take $x\in\Int D^p_-$. 
Denote by $G_0:T^{p,q}_+\to\R^m$ the standard embedding.
Recall the definition of $h_f$ from the definition of linking coefficients in \S4. 

\smallskip
{\it Definition of $\overline\lambda_q$.}
Take an embedding $G:T^{p,q}_+\to S^m$.
Since $m\ge 2p+q+2$, by general position $G$ has a unique up to isotopy 
extension to an embedding $G':T^{p,q}-\Int B^{p+q}\to S^m$.
Define $\overline\lambda_q(G)$ to be the homotopy class of the map 
$$G'|_{\partial B^{p+q}}:
\partial B^{p+q}\to S^m-G(x\times S^q)\overset{h_G}\to\simeq S^{m-q-1}.$$ 

\smallskip
{\it Definition of $\overline\mu$.}
For $x\in\pi_{p+q}(S^{m-q-1})$ take a map 
$S^{p+q}\to S^m-G_0(T^{p+1,q}_+)\overset{h_{G_0}}\to\simeq S^{m-q-1}$
representing the class $x$.
Define $\overline\mu(x)$ to be the connected sum of the embedding 
$G_0|_{T^{p,q}}$ with this map.  

\smallskip
{\it Proof of the exactness at $KT^m_{p,q,+}$ in the Restriction Lemma 5.2.}
By definition $\overline\lambda_q(G)$ is the obstruction to extending an embedding 
$G:T^{p,q}_+\to S^m$ to an almost embedding $T^{p,q}\to\R^m$, so
the sequence is exact at $KT^m_{p,q,+}$.
\qed

\smallskip
{\it Proof of the exactness at $\overline{KT}^m_{p,q}$ in the Restriction
Lemma 5.2.}
Clearly, $\overline\nu\overline\mu=0$.
Let $G:T^{p,q}\to S^m$ be an almost embedding such that $\overline\nu(G)=0$. 
Then $G$ is isotopic to an almost embedding standard on $T^{p,q}_+$. 
Thus we may assume that $G$ itself is standard on $T^{p,q}_+$. 
Since $m\ge 2p+q+2$, by general position $G$ is isotopic relative to 
$T^{p,q}_+$ to an embedding standard outside $B^{p+q}$.
Thus we may assume that $G$ itself is standard outside $B^{p+q}$.

Hence $G|_{B^{p+q}}$ and $G_0|_{B^{p+q}}$ form together a map
$S^{p+q}\to S^m-G_0(\Int T^{p,q}_+)\simeq S^{m-q-1}$.
Let $x$ be the homotopy class of this map. 
Then $G=\overline\mu(x)$
(because $[G(B^{p+q})\cup G_0(B^{p+q})]\cap GT^{p,q}_+=\emptyset$, so there 
is an isotopy of $S^m\rel G(T^{p,q}_+)$ moving $G_0(T^{p+1,q}_+)$ to a 
neighborhood in $S^m$ of $GT^{p,q}_+\rel G\partial T^{p,q}_+$
that misses $G(B^{p+q})\cup G_0(B^{p+q})$).

\smallskip
{\it Definition of $\lambda(f)\in\pi_{p+q}(S^{m-q-1})$ for an almost embedding 
$f:T^{p,q}\to S^m$ coinciding with $G_0$ on $T^{p,q}_+$.}
The restrictions of $f$ and $G_0$ onto $T^{p,q}_-$ coincide on the boundary and 
so form together a map 
$$\lambda'(f):S^p\times S^q\to S^m-f(x\times S^q)\simeq S^{m-q-1}.$$ 
The map $S^p\times S^q\to(S^p\times S^q)/(S^p\vee S^q)\cong S^{p+q}$ induces
 a 1--1 correspondence between homotopy classes of maps
$S^p\times S^q\to S^{m-q-1}$ and $S^{p+q}\to S^{m-q-1}$.  
Let $\lambda(f)$ be the homotopy class of the map corresponding to 
$\lambda'(f)$. 
(Note that $x=\lambda(G)$ in the previous proof.)

\smallskip
{\it Proof of the exactness at $\pi_{p+q}(S^{m-q-1})$ in the Restriction Lemma 
5.2.}
First we prove that $\overline\mu\overline\lambda_{q+1}=0$ (note that only this part of
the exactness is required for the proof of the main results).
Take an embedding $\psi:T^{p,q+1}_+\to\R^{m+1}$.
Represent  $S^{q+1}=D^{q+1}_<\cup D^q\times I\cup D^{q+1}_>$.
Analogously to the Standardization Lemma 2.1 for
$m\ge2p+q+2$ in the PL category and for $m\ge2p+q+3$ in the smooth
category making an isotopy we can assume that 

$\psi|_{D^p_+\times D^{q+1}_>}$ is the standard proper embedding into
$\R^m\times[1,+\infty)$, 

$\psi|_{D^p_+\times D^{q+1}_<}$ is the standard proper embedding into 
$\R^m\times(-\infty,0]$, and

$\psi|_{D^p_+\times D^{q+1}\times I}$ is a concordance between
standard embeddings. 

The latter concordance is ambient, so there is a homeomorphism 
$$\Psi:\R^m\times I\to\R^m\times I\quad\text{such that}\quad \Psi(\psi(y,0),t)=
\psi(y,t)\quad\text{for each }y\in D^p_+\times D^{q+1},\ t\in I.$$ 
So the standard embedding $G_0:T^{p,q}\to\R^m$ is concordant to an embedding 
$\overline\psi_1:T^{p,q}\to\R^m\times1$ defined by 
$\overline\psi_1(y):=\Psi(G_0(y),1)$. 
The concordance $\Psi|_{G_0(T^{p,q})\times I}$ together with the standard 
extensions of $\psi|_{D^p_+\times D^{q+1}_{>,<}}$ to $S^p\times D^{q+1}_{>,<}$ 
form an extension of $\psi$ to $T^{p,q+1}-B^{p+q+1}$. 
Hence $\lambda(\overline\psi_1)=\overline\lambda_{q+1}(\psi)=
\lambda(\overline\mu\overline\lambda_{q+1}(\psi))$. 
Clearly, $\lambda(\overline\mu x)=x$, and 

{\it almost embeddings $f,g:T^{p,q}\to S^m$ coinciding with $G_0$ on 
$T^{p,q}_+$ are almost concordant relative to $T^{p,q}_+$ if and only if 
$\lambda(f)=\lambda(g)$.} 

Therefore $\overline\mu\overline\lambda_{q+1}(\psi)=\overline\psi_1$ is 
concordant to the standard embedding. 

Now we prove that $\ker\overline\mu\subset\im\overline\lambda_{q+1}$.
If $\overline\mu(x)=0$, then take the restriction 
$T^{p,q}_+\times I\to\R^m\times I$ of
an almost concordance between $\overline\mu(x)$ and the standard embedding.
This restriction is a concordance between standard embeddings and so can be 
completed to an embedding $\psi:T^{p,q+1}_+\to\R^{m+1}$.
Analogously to the above it is proved that $\overline\lambda_{q+1}(\psi)=x$.
\qed

\comment

As in the definition of $\overline\mu$ extend the standard embedding 
$$\psi|_{D^p\times\partial D^{q+1}_>}: D^p\times\partial D^{q+1}_>
\to\R^m\times 1\quad\text{to an embedding}\quad \overline\psi_1:=
\overline\mu\overline\lambda_{q+1}(\psi):S^p\times S^q\to\R^m\times 1.$$ 
Define an embedding 
$$\overline\psi_0:T^{p,q}\to\R^m\times 0\quad\text{by}\quad
\overline\psi_0(y):=\Psi(\overline\psi_1(y),0).$$ 
Then $\Psi$ is a concordance between $\overline\psi_0$ and $\overline\psi_1$. 

So it suffices to prove that $\lambda(\overline\psi_0,f)=0$. 
Since $\overline\psi_1=\overline\mu\overline\lambda_{q+1}(\psi)$, we have 
$\lambda(\overline\psi_1,f)=\overline\lambda_{q+1}(\psi)$. 
Thus it suffices to prove that
$\lambda(\overline\psi_1,f)-\lambda(\overline\psi_0,f)=
\overline\lambda_{q+1}(\psi)$.
Since $2m\ge3q+p+4$, by the Freudenthal Suspension Theorem 
the group $\pi_{p+q}(S^{m-q-1})$ is stable. 
Hence it suffices to prove that
$\Sigma^\infty(\lambda(\overline\psi_1,f)-\lambda(\overline\psi_0,f)-
\overline\lambda_{q+1}(\psi))=0$.

This relation is proved using the Pontryagin-Thom construction and
intersection interpretation of homotopy classes [Ko88, \S1]. 
In this argument we omit constructions of framings because they are evident 
(and analogous to [Ko88]). 
We may assume that $\overline\psi_1$ is a framed immersion for which there 
exists a framed immersion  
$$F_1:S^p\times D^{q+1}\to\R^m\times1\quad\text{such that}\quad F_1|_
{S^p\times\partial D^{q+1}}=\overline\psi_1|_{T^{p,q}_+}\cup f|_{T^{p,q}_+}.$$ 
Define a framed immersion  
$$F:S^p\times D^{q+1}\times I\to\R^m\times I\quad\text{by}\quad F(y,t):=
\Psi(F_1(y),t).$$
Then 
$F|_{S^p\times\partial D^{q+1}\times0}=
\overline\psi_0|_{T^{p,q}_+}\cup f|_{T^{p,q}_+}$.  
Denote $s^q:=f(x\times S^q)$. 
For $j\in\{0,1\}$ the class $\Sigma^\infty\lambda(\overline\psi_j,f)$ is 
represented by the immersed framed manifold 
$$\{y\in S^p\times D^{q+1}\ |\ F(y,j)\in s^q\times j\}
\subset S^p\times D^{q+1}\subset S^{p+q+1}.$$ 
The class $\Sigma\overline\lambda_{q+1}(\psi)$ is represented by the 
immersed framed manifold 
$$\{(y,t)\in T^{p,q}_+\times I\ |\ \Psi(\overline\psi_0(y),t)\in
\Psi(s^q\times I)\} \subset T^{p,q}_+\times I\subset S^{p+q+1}.$$ 
The framed immersion $\partial F:=F|_{S^p\times\partial(D^{q+1}\times I)}$ is 
formed by framed immersions $\Psi|_{\overline\psi_0T^{p,q}_+\times I}$, $F_0$, 
$F_1$ and $\Psi|_{fT^{p,q}_+\times I}$. 
Since $\Psi(fT^{p,q}_+\times I)\cap\Psi(s^q\times I)=\emptyset$, 
it follows that the class 
$\Sigma^\infty(\lambda(\overline\psi_1,f)-\lambda(\overline\psi_0,f)-
\overline\lambda_{q+1}(\psi))$
is represented by the 
immersed framed manifold 
$(\partial F)^{-1}\Psi(s^q\times I)\subset S^p\times S^{q+1}$. 
The 
immersed framed manifold 
$F^{-1}\Psi(s^q\times I)\subset S^p\times D^{q+1}\times I$ is a 
framed bordism from the latter to zero. 
This proves the required relation. 

\endcomment

\smallskip
For the sequel we need $\Pi^{m-1}_{p,q}:=\pi^{m-1}_{eq}(S^p\times S^{2q})$, 
where the involution on $S^p\times S^{2q}$ is the product of the antipodal
involution on $S^p$ and the symmetry with respect to $S^q\subset S^{2q}$. 
The group structure on $\Pi^{m-1}_{pq}$ is defined as $S^p$-parametric version 
of the group structure on $[S^{2q},S^{m-1}]$, see the details in [Sk02, Torus 
Lemma 6.1, Sk07, Torus Lemma 6.1].
An element $\varphi\in\pi_q(V^{eq}_{m-q,p+1})$ can be considered
as a map $\varphi:S^p\times S^q\to S^{m-q-1}$ such that
$\varphi(-x,y)=-\varphi(x,y)$ for each $x\in S^p$.
Define a homomorphism $\sigma:\pi_q(V^{eq}_{m-q,p+1})\to\Pi^{m-1}_{p,q}$ 
by setting $\sigma(\varphi)$ to be the $q$-fold $S^p$-fiberwise suspension
of such a map $\varphi$, i.e.\  $\sigma(\varphi)|_{x\times
S^{2q}}=\Sigma^q(\varphi|_{x\times S^q})$.  
By [Sk02, Proof of Torus Lemma 6.1, Sk07, Proof of Torus Lemma 6.1], 
$\sigma$ is an isomomorphism for $2m\ge3p+q+4$.

\smallskip
{\it Proof of the commutativity in the Restriction Lemma 5.2.}
In this proof 'commutativity' means 'commutativity up to sign'.
The commutativity of the middle squares is trivial.

It is clear and well-known [HH62] that the big vertical rectangulars are 
commutative (because the composition 
$\pi_q(S^{m-p-q-1})\overset{\Sigma^p}\to\to\pi_{p+q}(S^{m-q-1})\cong
\pi_q(\Omega_pS^{m-q-1})$ is induced by $\rho_{p+1}$).

Let us prove the commutativity of the left upper square.
The case $p=0$ is clear.
The case of arbitrary $p$ follows by the Extension Lemma 2.6.b because $m\ge2p+q+2$ and
$\tau_p\mu''_px:T^{p,q}\to\R^m$ is a mirror-symmetric extension of
$\tau_{p-1}\mu''_{p-1}x:T^{p-1,q}\to\R^{m-1}$ and
$\overline\mu_p\Sigma^px:T^{p,q}\to\R^m$ is a mirror-symmetric extension of
$\overline\mu_{p-1}\Sigma^{p-1}x:T^{p-1,q}\to\R^{m-1}$.

The commutativity of the right upper square for $2m\ge3q+p+3$ follows by the
commutativity of the right bottom square (proved below) and the commutativity 
of a right big vertical rectangular.

By [Sk02, Torus Lemma 6.1] the bottom line can be identified with the 
following sequence:
$$\CD @>> \lambda_{q+1} > \pi_{p+2q}(S^{m-1})@>> \mu > \Pi^{m-1}_{p,q}
@>> \nu > \Pi^{m-1}_{p-1,q} @>> \lambda_q > \pi_{p+2q-1}(S^{m-1})\endCD.$$ 
Under this identification the vertical arrows marked with the equality are 
identified with the suspension homomorphisms $\Sigma^q$ (which are isomorphisms 
because $2m\ge3q+p+3$).
Hence the commutativity of the left bottom square follows by the commutativity 
of the left bottom square from [Sk02, Decomposition Lemma 7.1].

Under this identification the map $\lambda_q'$ is identified with the map 
$\lambda_q$ defined as the obstruction to extension over $T^{p,2q}_+$. 
Hence the commutativity of the right bottom square for $2m\ge3q+p+3$
follows from the relation $\lambda_q\alpha'r=\Sigma^q\overline\lambda_q$ which 
is clear by the proof of [Sk02, proof of Torus Lemma 6.1, Sk07, proof of
Torus Lemma 6.1]. 
\qed

\smallskip
{\it Proof of the relations $\bar\alpha'\tau=\pm\rho_{p+1}$ and 
$\alpha'r\tau_+=\pm\rho_p$ in the Restriction Lemma 5.2.}
It suffices to prove the first relation. 
It reduces to $\alpha\tau=\pm\sigma\rho_{p+1}$. 
Represent $\varphi\in\pi_q(V_{m-q,p+1})$ as a map $S^p\times S^q\to S^{m-q-1}$. 
The relation is proved for $p=0$ using the representation 
$S^{m-1}\cong S^{m-q-1}*S^{q-1}$ and deforming $\alpha\tau(\varphi)$ 
to the ($S^p$-fiberwise) suspension $\sigma(\varphi)$ of $\varphi$ [Ke59].
The same relation for $p>0$ is proved by applying this deformation for each 
$x\in S^p$ independently.
\qed

\smallskip
{\it Proof of the assertion on $\bar\alpha$ in the smooth category 
in the Torus Lemma 2.7.}
Consider the bottom two lines of the diagram from the Restriction Lemma 5.2.
By Lemma 5.1 and [Ha63, 6.4, Sk02, Theorem 1.1.$\alpha\partial$] the map 
$\alpha'r$ is an isomorphism for $2m\ge3q+2p+2$ and an epimorphism for 
$2m=3q+2p+1$. 
Hence by the 5-lemma it follows that $\bar\alpha$ is an isomorphism for 
$2m\ge3q+2p+2$.
\qed

\smallskip
{\it Proof of the assertions on $\alpha\tau$ and the 'moreover' part 
in the Torus Theorem 2.8.}
Let $n=p+1$ and $l=m-p-q$.
Identify the sets of Lemma 5.1 by the isomorphisms  of Lemma 5.1.

In the following two paragraphs we reproduce the argument from [HH62, (1.1)].
The proof is by induction on $n$.
For $n=1$ the map $\rho_n$ is an isomorphism because 
$V_{l+1,1}\cong V_{l+1,1}^{eq}\cong S^l$. 
Suppose now that $n\ge2$.
Consider the following diagram formed by the upper and the bottom line of
the diagram from the Restriction Lemma 5.2:
$$\minCDarrowwidth{5pt}\CD
\pi_{q+1}(V_{l+n,n-1}) @>> >
\pi_q(S^l) @>> \mu'' >
\pi_q(V_{l+n,n}) @>> \nu'' >
\pi_q(V_{l+n,n-1}) @>> \lambda'' >
\pi_{q-1}(S^l) \\
@VV \rho_{n-1} V @VV \Sigma^{n-1} V @VV \rho_n V @VV \rho_{n-1} V @VV \Sigma^{n-1} V \\
\pi_{q+1}(V^{eq}_{l+n,n-1}) @>> >
\pi_{q+n-1}(S^{l+n-1}) @>> \mu' >
\pi_q(V^{eq}_{l+n,n}) @>> \nu' >
\pi_q(V^{eq}_{l+n,n-1}) @>> \lambda' >
\pi_{q+n-2}(S^{l+n-1}) \endCD.$$
Suppose that $q\ge1$ (the argument for $q=0$ is the same, only the right-hand 
terms in the above diagram should be replaced by zeros, since both restrictions 
inducing $\nu''$ and $\nu'$ are surjective).

Let us prove that $\rho_n$ is an isomorphism for 
$2m\ge3q+2p+4\Leftrightarrow q\le2l-2$.
By the Freudenthal Suspension Theorem $\Sigma^{n-1}$ are isomorphisms.
By the inductive hypothesis $\rho_{n-1}$ are isomorphisms.
So by the 5-lemma $\rho_n$ is an isomorphism.

Let us prove the 'moreover' part.
We have $q=2l-1$ and $n\ge2$.
Hence by the Freudenthal Suspension Theorem the right-hand
$\Sigma^{n-1}$ is an isomorphism and the left-hand $\Sigma^{n-1}$
is an epimorphism whose kernel is generated by $[\iota_l,\iota_l]$.
By the inductive hypothesis $\rho_{n-1}$ are isomorphisms.
Therefore if we factorize $[\iota_l,\iota_l]\in\pi_q(S^l)$ and 
$w_{l,n-1}\in\pi_q(V_{l+n,n})$, by the 5-lemma we obtain that $\rho_n$ an 
epimorphism whose kernel is generated by $w_{l,n-1}$.
\qed

\smallskip
{\bf Remarks.} (a) {\it An alternative proof that the map $\tau$ of the Torus 
Theorem 2.8 is an isomorphism for $m\ge\max\{2p+q+2,3q/2+p+2\}$ in the PL 
category and for $m\ge\max\{2p+q+3,3(q+p)/2+2\}$ in the smooth category.}
Recall that the forgetful map $KT_{p,q}^m\to\overline{KT}_{p,q}^m$
is an isomorphism by [Sk02, Theorem 2.2.$\alpha$q] and general
position (and, in the smooth category, by smoothing). 
Consider the upper two lines of the diagram from the Restriction Lemma 5.2. 
The map $\tau_+$ is an isomorphisms for $2m\ge3q+2p+5$ by the assertion on
$\tau_+$ of Lemma 5.1, and is an epimorphism for $2m=3q+2p+4$. 
So by the induction on $q$ using the 5-lemma we obtain that $\tau$ is an 
isomorphism. 


(b) {\it An alternative proof of the commutativity of the upper right square 
from the Restriction Lemma 5.2 for $p=1$ and $q\le2(m-q-1)-2$.}
We omit index $q$. 
Recall that $\lambda''\Sigma x=(1-(-1)^{m-q})x$ for $x\in\pi_{q-1}(S^{m-q-2})$ 
[JW54] (in that paper $\lambda''$ was denoted by $\Delta$).
On the other hand, $\bar\lambda y$ is the signed sum of linking coefficients of 
the link $y|_{T^{0,q}}$, thus $\bar\lambda y=(1-(-1)^{m-q})y$.
By the Freudenthal Suspension Theorem it follows that for each $y\in \pi_q(S^{m-q-1})$ there
is $x\in\pi_{q-1}(S^{m-q-2})$ such that $y=\Sigma x$.
Hence $\Sigma\lambda''y=\Sigma\lambda''\Sigma x=\Sigma(1-(-1)^{m-q})x=\bar\lambda y$.

(c) {\it An alternative proof of the commutativity of the upper left square 
from the Restriction Lemma 5.2 for $2m\ge3q+2p+2$.}
Follows by the commutativity of the left bottom square and the commutativity
of a left square of a diagram from the proof of the assertions on $\alpha\tau$,
because $\bar\alpha$ is injective for $2m\ge3q+2p+2$ (which is proved using
only the two bottom lines of the diagram).


(d) {\it The map $\alpha\tau$ is an epimorphism for $q=2l$, $m=3l+p+1$, $p\ge2$ 
and $l+1$ not a power of 2.} 
(The condition that $l+1$ is not a power of 2 can be replaced by $w_{l,p}\ne0$.)
This can be proved using the diagram from the proof of the
assertions on $\alpha\tau$ and the diagram choice as if in the
proof of the Triviality Lemma 6.2. 
This fact together with the assertion on $\bar\alpha$ of the Torus Theorem 2.8 for 
the smooth category implies the Triviality Lemma 6.2 below.

\comment


By the Freudenthal Suspension Theorem the right $\Sigma^{n-1}$ is an
epimorphism whose kernel is $\left<[\iota_l,\iota_l]\right>$.
If $l$ is even, then $\left<[\iota_l,\iota_l]\right>$ is infinite while
$\pi_{2l}(V_{l+n,n-1})$ is finite (the latter is proved in the proof of
the Triviality Lemma 6.2).
Therefore $\im\lambda''\cap\left<[\iota_l,\iota_l]\right>=0$.
If $l$ odd and $l+1$ is not a power of 2, then
$\left<[\iota_l,\iota_l]\right>=\{[\iota_l,\iota_l],0\}$ and
$w_{l,n-1}\ne0$ [Os86].
Hence by exactness $[\iota_l,\iota_l]\not\in\im\lambda''$,
so again $\im\lambda''\cap\left<[\iota_l,\iota_l]\right>=0$.
Now $\rho_n$ is epimorphic because for any $x\in\pi_q(V^{eq}_{l+n,n})$
we have
$$\lambda'\nu'(x)=0\quad\Rightarrow
\quad\lambda''\nu'(x)\in\left<[\iota_l,\iota_l]\right>\quad\Rightarrow
\quad\lambda''\nu'(x)=0\quad\Rightarrow
\quad\nu'(x)=\nu''x'\text{ for some }x'\quad\Rightarrow$$
$$\Rightarrow\quad \nu'(x-\rho_n(x'))=0 \quad\Rightarrow
\quad x-\rho_n(x')=\mu'(y)\text{ for some }y\quad\overset{n\ge3}\to
\Rightarrow$$
$$\Rightarrow\quad y=\Sigma^{n-1}y'\text{ for some }y'\quad\Rightarrow
\quad x=\rho_n(x'+\mu''y').$$
Here the implication marked with '$n\ge3$' holds for $n\ge3$
by the James Double Suspension Theorem [Jam54].

{\it An alternative proof of the Triviality Lemma 6.2 for $p\ge2$.}
Consider the composition
$$\pi_{2l}(V_{l+p+1,p+1})\overset{r\tau}\to\to\overline{KT}^{3l+p+1}_{p,2l,DIFF}
\overset{\bar\alpha}\to\cong\pi_{eq}^{3l+p}(\t{T^{p,2l}}).$$ Since
$2(3l+p+1)\ge3\cdot2l+p+4$ and $2(3l+p+1)=3\cdot2l+2p+2$, by the
assertion on $\bar\alpha$ the Torus Theorem 2.8 for the smooth category 
$\bar\alpha$ is an isomorphism. Since $p\ge2$ and $l+1$ is not a
power of 2, by the above the map $\bar\alpha\tau=\alpha\tau$ is
onto. Hence $r\tau$ is onto, so $\beta=0$. (Instead of using the
assertion on $\bar\alpha$ for the smooth category, since $p\ge2$ and
$2(3l+p+1)=3(2l+p)+1-(p-1)$ we can use [Sk02, Theorem 2.3.q] to
see that the map $\bar\alpha$ defined on the set
$\t{KT}^{3l+p+1}_{p,2l}$ of {\it smooth} almost embeddings
$T^{p,2l}\to S^{3l+p+1}$ up to {\it PL} almost concordance is
injective.) (This proof is in fact not essentially different from
the previous one because the surjectivity of $\alpha\tau$ is
proved by the same diagram search as in the previous proof.)




\endcomment

\head 6. Proof of the properties of $\beta$-invariant \endhead


{\it Proof of the completeness in the $\beta$-invariant Theorem 2.9.}
By general position, since $3l+p+1\ge 2(p+1)+(2l-1)+1$, there is a homotopy of
$F\rel T^{p,2l-1}_+\times I\cup T^{p,2l-1}\times\{0,1\}$ to an almost
concordance (denoted by the same letter $F$) for which
$\Sigma(F)\subset B^{2l+p}:=D^p_-\times D^{2l-1}_-\times[\frac13,\frac23]$.
In the smooth category since $2(3l+p+1)\ge 3(p+1)+2(2l-1)+4$, we may assume 
that $F$ is a smooth embedding outside $B^{2l+p}$.
Let
$$M=\R^{3l+p}\times I-\Int R_{\R^{3l+p}\times I}
(F(T^{p,2l}\times I-\Int B^{2l+p}),F\partial B^{2l+p}).$$
Observe that $F|_{B^{2l+p}}:B^{2l+p}\to M$ is a proper map whose restriction to
the boundary is an embedding.
Since $T^{p,2l}\times I$ is homologically $(p-1)$-connected, by
Alexander duality $M$ is homologically $(l+p-1)$-connected.
Since $M$ is simply-connected, it follows that $M$ is $(l+p-1)$-connected.

Since $\beta(F)=0$, by Alexander duality it follows that 
$[FC]=0\in H_{l+p}(M;\Z_{(l)})$.
Since $[FC]=0$, $M$ is $(l+p-1)$-connected and $3l+p+1\ge2l+p+3$,
by the completeness theorem [Hu70', Theorem 2, Ha84, Theorem 4]
it follows that $F|_{B^{2l+p}}$ is homotopic $\rel\partial
B^{2l+p}$ to an embedding $F'$. Extending $F'$ over $N$ by $F$ we
obtain an embedding $F'$. Thus we obtain the required almost
concordance from $F$ to a concordance. \qed


\proclaim{Reduction Lemma 6.1}
For $q=2l-1$, $m=3l+p$ and $p<l\ge2$ the following conditions are equivalent:

(G) $\beta(G)=0$ for every almost embedding $G:T^{p,q+1}\to S^{m+1}$;

(F) for almost concordances $F$ between embeddings $f_0,f_1:T^{p,q}\to\R^m$
almost concordant to the standard embedding,
the invariant $\beta(F)$ depends only on $f_0$ and $f_1$.
\endproclaim

\smallskip
{\it Proof that (F) implies (G).}
By the Standardization Lemma 2.1 and the assertion in the proof of the 
Triviality Criterion in \S3 (or just analogously to the Standardization Lemma 
2.1) 
we may assume that $G$ is standardized, 
$G(S^p\times\frac12 D^{q+1}_+)\subset\frac12 D^{m+1}_+$ is the standard 
embedding and 
$G(S^p\times(D^{q+1}_+-\frac12D^q_+))\subset D^{m+1}_+-\frac12D^{m+1}_+$.  
Thus the restriction 
$F:S^p\times(D^{q+1}_+-\frac12D^{q+1}_+)\to D^{m+1}_+-\frac12D^{m+1}_+$ of $G$ 
is an almost concordance between standard embeddings. So
$\beta(G)=\beta(F)=0$. \qed


\smallskip
{\it Proof that (G) implies (F).}
For $s=0,1$ let $F_s$ be any almost concordances from $f_s$ to the standard
embedding.
Given an almost concordance $F$ between $f_0$ and $f_1$ there is an
almost concordance $F'=\overline F_0\cup F\cup F_1$ between standard
embeddings.
This $F'$ can be 'capped' to obtain an almost embedding $G:T^{p,q+1}\to S^{m+1}$
without new self-intersections.
Hence by the additivity of the $\beta$-invariant Theorem 2.9 we have
$0=\beta(G)=\beta(F')=\beta(F)+\beta(F_0)-\beta(F_1)$.
Thus $\beta(F)=\beta(F_1)-\beta(F_0)$ is independent on $F$.
\qed

\smallskip
{\it The Additivity Theorem 2.10.b} is proved analogously to the following 
absolute version. 

\smallskip
{\it Proof that $\beta(F+F')=\beta(F)+\beta(F')$ for almost embeddings
$F,F':T^{p,2l}\to\R^{3l+p+1}$ and $0<p<l$.} 
(Recall that for $p=l-1$ in the DIFF or AD category the sum is not necessarily 
well-defined and we denote by $F+F'$ any sum.) 
Since $p<l$, we have $3l+p+1\ge2p+2l+2$.
Hence we may assume that $F$ and $F'$ are standardized.
Then we may assume that the supports of $C_F$ and $C_{F'}$ are in 
the ball $S^p\times D^{2l}_+\cap T^{p,2l}_-$.

Denote $m:=3l+p+1$. 
Since $R_m$ and $R_{2l}$ are isotopic to the identity maps of $\R^m$ and of 
$S^{2l}$, they do not change orientations. 
Hence $[\Sigma(F+F')]=[\Sigma(F)]+(\id S^p\times R_{2l})[\Sigma(F')]$. 
So we can take 
$$C_{F+F'}:=C_F+(\id S^p\times R_{2l})C_{F'}\quad\text{and}
\quad D_{F+F'}:=D_F+R_mD_{F'}.$$ 
We may assume that the supports of $D_F$ and $D_{F'}$ are in $\R^m_+$. 
Then  
$$\beta(F+F')=D_{F+F'}\cap[(F+F')(x\times S^{2l})]=
(D_F+R_mD_{F'})\cap[F(x\times D^{2l}_+)+R_mF'(x\times D^{2l}_+)]=$$
$$=D_F\cap[F(x\times D^{2l}_+)]+R_mD_{F'}\cap[R_mF'(x\times D^{2l}_+)]=
\beta(F)+\beta(F').\quad\qed$$



Note that for $p\le l-2$ the equality $\beta(F+F')=\beta(F)+\beta(F')$ holds 
for {\it general} almost concordances $F$ and $F'$ (for which $f_0$ and $f_1$ 
are not supposed to be almost concordant to the standard embedding).



\smallskip
{\it An alternative proof that for $p\le l-2$ (G) implies that for almost concordances $F$
between {\it arbitrary} embeddings $f_0,f_1:T^{p,q}\to\R^m$ the invariant $\beta(F)$ depends
only on $f_0$ and $f_1$.}
Let $F$ and $F'$ be almost concordances between embeddings
$f_0,f_1:T^{p,q}\to\R^m$.
Then $F-F'$ is an almost concordance between embeddings isotopic to the
standard embedding.
Hence $F-F'$ can be 'capped' to obtain an almost embedding
$G:T^{p,q+1}\to S^{m+1}$ without new self-intersections.
So by the Additivity Theorem 2.10.b we have 
$\beta(F)-\beta(F')=\beta(F-F')=\beta(G)=0$.
\qed

\smallskip
The independence on $C$ for $p=l-1$ and arbitrary $f_0,f_1$ cannot be reduced
to the case when $f_0$ and $f_1$ are almost concordant to the standard embedding
(as above) because $\beta(F+F')=\beta(F)+\beta(F')$ is not proved for $p=l-1$ and {\it general} almost concordances $F$ and $F'$ .

\proclaim{Triviality Lemma 6.2} Suppose that either $1\le p\le l-2$ and
$l+1$ is not a power of 2, or $1=p=l-1$. 
Then $\beta(G)=0$ for every
smooth almost embedding $G:T^{p,2l}\to\R^{3l+p+1}$.
\endproclaim

\proclaim{Non-triviality Lemma 6.3} For $1\le p<l\in\{3,7\}$ there
exists a smooth almost embedding $G:T^{p,2l}\to\R^{3l+p+1}$ such
that $\beta(G)=1$.
\endproclaim

We postpone the proof of these lemmas.
Note that the Triviality Lemma 6.2 is true for $p=0\le l-2$. 

{\it The independence in the $\beta$-invariant Theorem 2.9} follows from
$(G)\Rightarrow(F)$ of the Reduction Lemma 6.1 and the Triviality Lemma 6.2.

{\it The triviality in the $\beta$-invariant Theorem 2.9} follows from
$(F)\Rightarrow(G)$ of the Reduction Lemma 6.1, the Non-Triviality Lemma 6.3 
and the completeness of $\beta$-invariant.

For the proof of the Triviality Lemma 6.2 and the Non-triviality Lemma 6.3 we 
need the Restriction Lemma 5.2 of \S5 and so the condition $3l+p+1\ge2p+2l+2$, i.e.
$p<l$.
More specifically, for the Non-triviality Lemma 6.3 we only need the definition of
$\overline\mu$, for the Triviality Lemma 6.2 in the case when $l$ is even we only need the
definitions of $\overline\mu$ and $\overline\nu$ together with the exactness at
$\overline{KT}_{p,q}^m$, and for the Triviality Lemma 6.2 in the general case we need
the upper two lines of the Restriction Lemma 5.2 except the exactness at
$\pi_{p+q}(S^{m-q-1})$.

From now until the end of this section we work in the smooth category which we 
omit from the notation.

\demo{Proof the Non-triviality Lemma 6.3}
Consider the following diagram, in which $H$ is the Hopf invariant.
$$\minCDarrowwidth{0pt}\CD
\pi_{2l+p}(S^{l+p})  @>>\overline\mu> \overline{KT}^{3l+p+1}_{p,2l}\\
@VV \Sigma^\infty V @VV \beta V \\
\pi^S_l @>> H > \Z_2 \endCD$$
The diagram is (anti)commutative analogously to [Ko88, Theorem 4.8].
Since $p\ge1$, it follows that the group $\pi_{2l+p}(S^{l+p})$ is either
stable or metastable, so $\Sigma^\infty$ is epimorphic.
Since $l\in\{3,7\}$, it follows that $H$ is epimorphic.
Hence $\beta$ is epimorphic.
\qed\enddemo

\demo{Proof of the Triviality Lemma 6.2}
Since $2(3l+p+1)\ge3\cdot2l+4$, we can
identify $\pi_{2l}(V_{l+p+1,p})$ with $KT^{3l+p+1}_{p,2l,+}$ by
the isomorphism $\tau_+$ of Lemma 5.1 which we omit from notation.

The set $\overline{KT}_{p,2l}^{3l+p+1}$ is a group for $p\le l-2$ or 
$1=p=l-1$ by Group Structure Lemma 2.2 and the Almost Smoothing Theorem 2.3. 

In this paragraph assume that {\it $l$ is even}.
By the Additivity Theorem 2.10.b $\beta:\overline{KT}^{3l+p+1}_{p,2l}\to\Z$ is 
a homomorphism. 
So it suffices to prove that its domain is finite.
Consider the following exact sequence given by the second line of the diagram 
from the Restriction Lemma 5.2 for $m=3l+p+1$ and $q=2l$:
$$\pi_{p+2l}(S^{l+p}) @>>\overline\mu>\overline{KT}_{p,2l}^{3l+p+1}
@>>\overline\nu> \pi_{2l}(V_{l+p+1,p}).$$
The domain of $\overline\mu$ is finite for $l$ even.
The range of $\overline\lambda$ is finite for $l\ge2$.
(The finiteness is proved by the induction on $p$: for $p=1$ and $l\ge2$ the 
group $\pi_{2l}(S^{l+1})\cong\pi_{l-1}^S$ is finite; the inductive step 
follows by applying exact sequence 
$\pi_{2l}(V_{l+p,p-1})\to\pi_{2l}(V_{l+p+1,p})\to\pi_{2l}(S^{l+p})$.)
Hence $\overline{KT}^{3l+p+1}_{p,2l}$ is finite.


In this paragraph assume that {\it $l$ is odd and $l+1$ is not a power of 2.}
Consider the first two lines of the diagram from the Restriction Lemma 5.2 for 
$m=3l+p+1$ and $q=2l$ (except of the very left column).
By $\Sigma^p$ we denote {\it the right} $\Sigma^p$ of the diagram, which maps 
the metastable group to the stable group.
By $\mu''$ we denote the map from this metastable group to the right
(which is not shown on the diagram).
Then $\left<[\iota_l,\iota_l]\right>=\{[\iota_l,\iota_l],0\}$ and
$w_{l,p}=\mu''[\iota_l,\iota_l]\ne0$ [Os86].
Hence by exactness $[\iota_l,\iota_l]\not\in\im\lambda''$,
so $\im\lambda''\cap\ker\Sigma^p=0$.
Hence for any $x\in\overline{KT}^{3l+p+1}_{p,2l}$ we have
$$\overline\lambda\overline\nu(x)=0\quad \Rightarrow
\quad \lambda''\overline\nu(x)\in\ker\Sigma^p\quad \Rightarrow
\quad \lambda''\overline\nu(x)=0\quad \Rightarrow
\quad \overline\nu(x)=\nu''x'\text{ for some }x'\quad \Rightarrow$$
$$\Rightarrow\quad \overline\nu(x-\tau(x'))=0 \quad\Rightarrow
\quad x-\tau(x')=\overline\mu(y)\text{ for some }y.$$
So $\overline{KT}^{3l+p+1}_{p,2l}=\im\tau+\im\overline\mu$.
Since $\beta\tau=0$ and $\beta\overline\mu=H\Sigma^\infty=0$, the Lemma follows.
\qed\enddemo


\head 7. Proof of the Almost Smoothing Theorem 2.3 \endhead

We need the following definitions and results.
Recall from the beginning of \S2 the convention concerning PL and PD 
categories.  
Recall that $T^{p,q}_+=D^p_+\times S^q$. 
Let $i:T^{m-q,q}_+\to\R^m$ and $D^p\subset D^{m-q}$ be the
standard inclusions.
An embedding $g:T^{p,q}_+\to\R^m$ is called {\it fiberwise} if 
$$g(0,x)=i(0,x)\quad\text{and}\quad g(D^p\times x)\subset i(D^{m-q}\times x)
\quad\text{for each}\quad x\in S^q.$$
Analogously a fiberwise concordance is defined.

\smallskip
{\bf Fibering Lemma 7.1.} {\it For $m\ge p+q+3$ every PL embedding
$f:T^{p,q}_+\to\R^m$ is PL isotopic to a fiberwise PL embedding
and every PL concordance between fiberwise PL embeddings $T^{p,q}_+\to\R^m$
is PL isotopic relative to the ends to a fiberwise PL concordance.}

\smallskip
Note that codimension 0 analogue of the Fibering Lemma 7.1 is false 
[Hi66, Corollary B].

\smallskip
{\bf Slicing Lemma 7.2.} {\it Let $X$ be a finite simplicial $n$-complex, 
$m-n\ge3$ and $s:X\to\R^m$ a simplicial map for some triangulation $T$ of 
$\R^m$.  
Then for each $\varepsilon>0$ there exists $\delta>0$ such that if a PL 
embedding $f:X\to\R^m$ is $\delta$-close to $s$, then $f$ is PL 
$\varepsilon$-ambient isotopic to a PL embedding $h:X\to\R^m$ such that
$h^{-1}(\sigma)=s^{-1}(\sigma)$ for each cell $\sigma$ of the
dual to $T$ cell-subdivision of $\R^m$} [Me02, Lemma 4.1].

\smallskip
{\it Proof of the Fibering Lemma 7.1.}
(This proof appeared in a discussion of the author with S. Melikhov.)
We prove only the assertion on embedding, the assertion on concordance is 
proved analogously.
By the Zeeman Unknotting Spheres Theorem [Ze62] we may assume that
$f(0,x)=i(0,x)$ for each $x\in S^q$.
By making an isotopy of $f$ which 'shrinks' $f(D^p\times x)$ we may assume that
$f(T^{p,q}_+)\subset i(T^{m-q,q}_+)$ and $f(a,x)=i(b,y)$ implies that 
$|x,y|<\delta$. 

Apply the Slicing Lemma 7.2 to the standard embedding $s$ and the embedding $f$ 
to obtain an embedding $h$. 
Then for some cell-subdivision $T_1$ of $S^q$ we have 
$h(D^p\times\sigma_1)\subset i(D^{m-q}\times\sigma_1)$ for each dual cell 
$\sigma_1$ of $T_1$ (because $D^{m-q}\times\sigma_1$ is a dual cell for 
certain triangulation of $\R^m$).

By induction on $k$ we may assume that $h(D^p\times x)\subset D^{m-q}\times x$
for points $x$ of the $k$-skeleton of $T_1$. 
Recall the relative $k$-Concordance Implies $k$-Isotopy Theorem: 

{\it if an embedding $g:D^p\times D^k\to D^{m-q}\times D^k$ on 
$D^p\times\partial D^k$ commutes with the projection onto the second factor, 
and $g(0,x)=i(0,x)$ for each $x\in D^k$, 
then the embedding is isotopic relative to 
$D^p\times\partial D^k\cup 0\times D^k$ to an embedding commuting with 
the projection onto the second factor.} 

Applying this result for each $(k+1)$-cell of $T_1$ independently we obtain an 
isotopy of $h$ relative to the $k$-skeleton to an embedding $f'$ such that 
$f'(D^p\times x)\subset D^{m-q}\times x$ for points $x$ of the 
$(k+1)$-skeleton of $T_1$. 
This proves the inductive step.  
For $k=n$ we obtain a fiberwise PL embedding.
\qed

\smallskip
{\it Proof of the Almost Smoothing Theorem 2.3.} Since any embedding
$S^0\to S^{m-q-1}$ is canonically isotopic to an embedding whose
image consists of two antipodal points, by the Fibering Lemma 7.1 it
follows that

{\it every PL embedding $T^{1,q}_+\to\R^m$ is PL isotopic to a
smooth embedding, and every PL concordance between smooth embeddings 
$T^{1,q}_+\to\R^m$ is PL isotopic to a smooth concordance.}

Now the Almost Smoothing Theorem 2.3 follows because

{\it for $m\ge q+\frac{3p+5}2$ the group $KT^m_{p,q,AD}$ is isomorphic to that of
PL embeddings $T^{p,q}\to\R^m$ which are smooth embeddings on $T^{p,q}_+$,
up to PL isotopy which is a smooth isotopy on $T^{p,q}_+\times I$.}

Let us prove the latter statement.
By [Ha67, Ha] for such an embedding $T^{p,q}\to\R^m$ the obstructions to
extension of smoothing from $T^{p,q}_+$ to
$T^{p,q}_+\bigcup\limits_{\partial D^p_-\times D^q_+} D^p_-\times D^q_+$ are in
$$H^{i+1}(D^p_-\times D^q_+,\partial D^p_-\times D^q_+;C^{m-q}_i)
\cong
H^{i+1}(S^p;C^{m-q}_i)\cong0\quad\text{for}\quad2(m-q)\ge3p+4.$$
Therefore every embedding as above is {\it almost smoothable},
i.e. is PL concordant to an almost smooth embedding.

Analogously, for $2(m+1-q)\ge 3(p+1)+4$ every PL isotopy between
smooth embeddings, which is a smooth isotopy on $T^{p,q}_+\times
I$, is {\it almost smoothable}, i.e is PL concordant relative to
the ends to an almost smooth concordance.
\qed

\subhead Appendix: Manifolds with boundary \endsubhead

We show that the dimension restriction is sharp
in the completeness results for $\alpha$-invariant of manifolds with boundary
[Ha63, 6.4, Sk02, Theorems 1.1.$\alpha\partial$ and 1.3.$\alpha\partial$]. 
This is (much simpler) non-closed analogue of the Whitehead Torus Example 1.7.

\proclaim{Filled-tori Theorem 7.3}
(a) $KT^m_{p,q,+,PL}\cong\pi_q(V^{PL}_{m-q,p})$ for $m\ge p+q+3$, where the PL
Stiefel manifold $V^{PL}_{m-q,p}$ is the space of PL embeddings $S^{p-1}\to S^{m-q-1}$.

(b) The map $\alpha^m(T^{1,q}_+)$ is not surjective
if $m\ge q+4$ and $\Sigma^\infty:\pi_q(S^{m-q-1})\to\pi_{2q+1-m}^S$ is not
epimorphic (concrete examples can be found from the table in [Sk02,
Example 1.4.s]).

(c) The map $\alpha^m(T^{1,q}_+)$ is not injective if
$m\ge q+4$ and
$\Sigma^\infty:\pi_q(S^{m-q-1})\to\pi_{2q+1-m}^S$ is not monomorphic, e.g.
for $q=2l-1$, $m=3l$, $l\ge4$ and $l\ne7$.

(d) The map $\alpha_{DIFF}^{3l+p-1}(T^{p,2l-1}_+)$ is
not injective for each even $l$ and $1\le p<l$.
\endproclaim

{\it Proof.} Part (a) follows by the Fibering Lemma 7.1.

Since $V^{CAT}_{m-q,1}\simeq S^{m-q-1}$, by part (a) and the assertions on 
$\tau$ and $\alpha\tau$ of the Torus Theorem 2.8 it follows that
$$KT^m_{1,q,+}\cong\pi_q(S^{m-q-1})\quad\text{and}\quad 
\pi^{m-1}_{eq}(\t{T^{1,q}_+})\cong\pi^S_{2q+1-m}.$$
Analogously to the proof of the Torus Theorem 2.8, $\alpha$ corresponds to the
suspension under the above isomorphisms.
This implies parts (b) and (c).

Part (d) follows because
$$KT^{3l+p-1}_{p,2l-1,+,DIFF}\cong\pi_{2l-1}(V_{l+p,p})\quad
\text{is infinite and}\quad\pi_{eq}^{3l+p-2}(\t{T^{p,2l-1}_+})\cong
\Pi^{3l+p-2}_{p-1,2l-1}\quad\text{is finite}$$
by the asserion on $\tau$ of the Torus Theorem 2.8 and [Sk02, Lemma 7.3.a].
The infiniteness of $\pi_{2l-1}(V_{l+p,p})$ for $l$ even and $1\le p<l$ is
proved by the induction on $p$ [cf. Sk02, Proof of Lemma 7.3.a].
The base $p=1$ is due to Serre.
The inductive step is proved using the following exact sequence for $p>1$:
$$\dots\to\pi_{2l}(S^{l+p-1})\to\pi_{2l-1}(V_{l+p-1,p-1})\to
\pi_{2l-1}(V_{l+p,p})\to\pi_{2l-1}(S^{l+p-1})\to\dots\qed$$

\head 8. Proofs of the Main Theorems 1.3 and 1.4 in the smooth category 
\endhead

{\bf Smoothing Theorem 8.1.} 
{\it For $m\ge2p+q+3$ there is an exact sequence of groups:
$$\dots\to KT^{m+1}_{p,q+1,AD}\overset{\sigma_{q+1}}\to\to C^{m-p-q}_{p+q}
\overset\zeta\to\to KT^m_{p,q,DIFF}\overset\forg\to\to KT^m_{p,q,AD}
\overset{\sigma_q}\to\to C^{m-p-q}_{p+q-1}\to\dots$$} 

For $m\ge p+q+3$ there is exact sequence of 
{\it sets} as in the Smoothing Theorem 8.1. 

\smallskip
{\it Definition of $\sigma_q$.}   
Take the codimension zero ball $B^{p+q}\subset T^{p,q}$ from the definition of 
an almost smooth embedding, and an almost smooth embedding $f:T^{p,q}\to\R^m$. 
By making a PL isotopy fixed outside $\Int B^{p+q}$ we may assume that $f$ is 
smooth outside a fixed single point [Ha67], cf. beginning of [BH70, Bo71].
Consider a small smooth oriented $(m-1)$-sphere $\Sigma$ with the center at 
the image of this point.
Take the natural orientation on the $(p+q-1)$-sphere $f^{-1}\Sigma$.
Let $\sigma(f)$ be the isotopy class of the abbreviation 
$f^{-1}\Sigma\to\Sigma$ of $f$. 
By [BH70, Bo71] $\sigma(f)$ is well-defined, i.e. independent on choices in the 
definition. 

\smallskip
For all applications of the Smoothing Theorem 8.1 but Theorem 1.2DIFF the exactness 
at $KT^m_{p,q,DIFF}$ is sufficient. 

\smallskip
{\it Proof of the exactness at $KT^m_{p,q,DIFF}$.}
Clearly, $\zeta$ is a homomorphism.
By the PL Unknotting Spheres Theorem $\forg\circ\zeta=0$.

Take a smooth embedding $f:T^{p,q}\to S^m$ such that $\forg f=0$.
Let $F$ be an almost smooth isotopy from $f$ to the standard embedding. 
Analogously to the above definition of $\sigma_q$ it is defined a complete
obstruction $\sigma(F)\in C^{m-p-q}_{p+q}$ to {\it smoothability}
of $F$, i.e. to the existence of an almost smooth concordance from
$F$ to a smooth embedding [Ha67, Ha, BH70, Bo71]. 

Take a smooth embedding $g:S^{p+q}\to S^m$ representing $-\sigma(F)$. 
From $\Sigma g:S^{p+q+1}\to S^{m+1}$ we can easily construct an almost smooth 
isotopy $G:S^{p+q}\times I\to S^m\times I$ between $g$ and the standard
embedding such that $\sigma(G)=[g]=-\sigma(F)$. 
Then $F\#G$ is an
almost smooth concordance from $f\#g=f-\zeta(\sigma(F))$ to the
standard embedding. We have $\sigma(F\#G)=\sigma(F)+\sigma(G)=0$,
so $F\#G$ is smoothable. Therefore
$f=\zeta(\sigma(F))\in\im\zeta$. \qed

\smallskip
{\it The exactness at $KT^m_{p,q,AD}$} is clear. 

\smallskip
{\it Proof of the exactness at $C^{m-p-q}_{p+q}$.} 
The proof appeared in a discussion with M. Skopenkov. 
Let $X$ be the set of proper almost smooth embeddings 
$S^p\times D^{q+1}\to D^m$ 
up to proper almost smooth isotopy (not necessarily fixed on the boundary).  
Consider the sequence 
$$KT^{m+1}_{p,q+1,AD}\overset u\to\to X\overset r\to\to KT^m_{p,q,DIFF}.$$ 
Here $r$ is the restriction to the boundary. 
In order to define $u$, take an almost smooth embedding 
$f:T^{p,q+1}\to\R^{m+1}$. 
By the Standardization Lemma 2.1 on embeddings we may assume that $f$ is 
standardized. 
Let $u(f)=f|_{T^{p,q+1}_+}:T^{p,q+1}_+\to\R^{m+1}_+$. 
By the Standardization Lemma 2.1 on concordances this is well-defined. 

By the Triviality Criterion of \S3 $ru=0$. 
Take an almost smooth embedding $f_+:T^{p,q+1}_+\to\R^{m+1}_+$ such that 
$rf=0$.
Capping a smooth isotopy from $rf$ to the standard embedding we obtain a smooth 
embedding $f_-:T^{p,q+1}_-\to\R^{m+1}_-$ agreeing with $f$ on the boundary. 
The embedding $f_-$ is smoothly isotopic to the standard embedding by 
the Uniqueness Lemma 8.2 below. 
Since this isotopy is ambient, $f_+\cup f_-$ is isotopic to a standardized 
embedding $f:T^{p,q+1}\to\R^{m+1}$ whose restriction to $T^{p,q+1}_+$ is 
properly smoothly isotopic to $f_+$. 
Hence $f_+=uf$. 
Thus the $(r,u)$-sequence is exact. 

There is a 'boundary connected sum with cone' action 
$\zeta':C^{m-p-q}_{p+q}\to X$. 
Analogously to the above definition of $\sigma_q$ it is defined a complete
obstruction $\sigma:X\to C^{m-p-q}_{p+q}$ to smoothing. 
Clearly, $\zeta'\sigma=\id_X$ and $\sigma\zeta'=\id_{C^{m-p-q}_{p+q}}$. 
Thus $\zeta'$ and $\sigma$ are 1--1 correspondences. 
Clearly, $r=\zeta'\zeta$ and $\sigma_{q+1}=\zeta'u$. 
Hence the sequence of the Smoothing Theorem 8.1 is exact at $C^{m-p-q}_{p+q}$. 

\smallskip
{\bf Uniqueness Lemma 8.2.} {\it If $m\ge2p+q+3$, then each two smooth proper 
embeddings $S^p\times D^q\to D^m$ are smoothly proper isotopic.} 

\smallskip
{\it Proof.} 
Use the assertion in the proof of the Triviality Criterion in \S3. 
Since $m\ge p+q+3$, the restriction 
$f|_{S^p\times(D^{q+1}_+-\frac12D^{q+1}_+)}$ is smoothly isotopic relative to 
the boundary to a smooth 
isotopy $F:S^p\times\partial D^{q+1}_+\times I\to\R^m\times I$ between 
$F_0=f|_{S^p\times\partial D^{q+1}_+}$ and the standard embedding $F_1$. 
For 
$$\tau\in[0,1]\quad\text{define}\quad\tau_*:\R^{m+1}\to\R^{m+1}\quad\text{by}
\quad\tau(x_0,x_1,\dots,x_m):=(x_0-\tau,x_1,\dots,x_m).$$ 
Take an isotopy $\Phi$ from $F\cup f|_{S^p\times\frac12D^{q+1}_+}$ to the 
standard embedding such that the image of 
$\Phi_\tau$ is the union of 
$\tau_*f(S^p\times\frac12D^{q+1}_+)$ and $\tau_*F_t$ for $t\ge\tau$. 
\qed


\smallskip
{\it Proof of the Main Theorem 1.3.b in the smooth case.} 
Follows by the AD case (proved in \S2 together with the PL case), the exactness 
at $KT^m_{p,q,DIFF}$ of the Smoothing Theorem 8.1 
and the equality $C_{4k+2}^{2k+2}=0$ [Ha66, 8.15, Mi72, Theorem F] because 
{\it the map $\forg$ is surjective for $p=1$, $q=2k+1$ and $m=6k+4$}. 
The surjectivity of $\forg$ follows because the composition of $\forg$ 
with $\tau\oplus\omega$ is surjective by the Realization Theorem 2.4.b.

\smallskip
We denote by $\widehat{h}$ certain quotient homomorphism of a homomorphism $h$.

\smallskip
{\it Proof of the Main Theorem 1.4.DIFF.}
The existence of the exact sequence follows from the Smoothing Theorem 8.1. 
The surjectivity of the map $\forg$ follows because the composition of $\forg$ 
with $\tau\oplus\omega$ is surjective by the Realization Theorem 2.3.b.

Let us prove that {\it for $l$ odd the exact sequence splits, i.e. the above
map $\forg$ has a right inverse}.
By the Main Theorem 1.4.AD and the Realization Theorem 2.4.b we can identify
$$KT^{3l+p}_{p,2l-1,AD}\quad\text{with}\quad 
\dfrac{\pi_{2l-1}(V_{l+p+1,p+1})}{w_{l,p}}\oplus\widehat{\Z_2}\quad
\text{via the isomorphism}\quad\widehat\tau_{AD}\oplus\widehat\omega_{AD}.$$
By the Relation Theorem 2.7 $\tau_{DIFF}(w_{l,p})=0$.
Then $\tau_{DIFF}\oplus\omega_{DIFF}$ factors through a right inverse of
$\forg$.
\qed

\smallskip
{\it Proof of the Main Theorem 1.3.aDIFF.}
Follows from the Main Theorem 1.4.AD,DIFF and the Smoothing Theorem 8.1 because {\it 
for $p=1$, $q=4k\ge8$ and $m=6k+1$ the above map $\forg$ has a right inverse.}

Let us prove the latter statement.
We omit the DIFF category from the notation.
Let $l=2k$.
By the Realization Theorem 2.4.a in the AD category we can identify 
$$KT^{3l+1}_{1,2l-1,AD}\quad\text{with}\quad
\pi_{l-2}^S\oplus\dfrac{\pi_{2l-1}(S^l)\oplus\Z}{[\iota_l,\iota_l]\oplus2}\quad
\text{via the isomorphism}\quad 
\tau^1_{AD}\oplus\widehat{(\tau^2_{AD}\oplus\omega_{AD})}.$$
By the $\tau$-relation of the Relation Theorem 2.7 we have
$\tau^2_{AD}[\iota_l,\iota_l]=2\omega_{AD}$.
Then by smoothing theory [BH70, Bo71] there is an element 
$\psi\in C_{2l}^{l+1}$ such that $\tau^2[\iota_l,\iota_l]-2\omega=\zeta(\psi)$.
Since $[\iota_l,\iota_l]$ is not divisible by 2 for $l$ even, $l\ge4$ [Co95],
it follows that $[\iota_l,\iota_l]\oplus2$ is primitive.
Hence there is a homomorphism
$$\varphi:\pi_{2l-1}(S^l)\oplus\Z\to C_{2l}^{l+1}\quad\text{such that}\quad
\varphi([\iota_l,\iota_l]\oplus2)=\psi.$$
Then $\tau^1\oplus[(\tau^2\oplus\omega)-\zeta\varphi]$
factors through a right inverse of $\forg$.
\qed

\smallskip
{\it Proof of Theorem 1.2DIFF.} 
For $2m\ge3q+2p+4$ we have $C^{m-p-q}_q=C^{m-p-q}_p=0$ [Ha66], so
$KT^m_{p,q,PL}=KT^m_{p,q,AD}$ by smoothing theory [BH70, Bo71], cf. [Ha67, Ha]. 

For $2m\ge3q+2p+4$ by the assertion on $\alpha\tau$ of the Torus Theorem 2.8 
the map $\tau_{DIFF}(\alpha\tau_{DIFF})^{-1}\alpha_{PL}$ is a right inverse of 
$\forg$ (hence $\sigma_q=0$).   
If $2m=3q+2p+3$, then $m\ge2p+q+3$ implies that $q\ge2p+3$, so 
$\sigma_q=0$ by the Main Theorem 1.4DIFF. 
Thus Theorem 1.2DIFF follows from Theorem 1.2PL and the Smoothing Theorem 8.1.

\head Appendix: Remarks, conjectures and open problems \endhead

See some open problems in [Sk07, \S3].

(0) Explicit constructions of \S2 allow to prove the following remark (which does not 
follow immediately from the definition of inverse elements).

\smallskip
{\bf Symmetry Remark.} {\it On $KT^7_{1,3,PL}\cong KT^7_{1,3,AD}$ 
we have $\id=\sigma_1=-\sigma_3=\sigma_7$.} 


\smallskip
{\it Proof.} Let $a_k:S^k\to S^k$ be the antipodal involution. 
We have
$$\sigma_1\tau^1=\tau^1(a_3\circ\iota_3)=\tau^1\quad\text{and}\quad
\sigma_3\tau^1=\sigma^3\tau^1(\Sigma\iota_2)=\tau^1(\sigma_3\circ\iota_3)=
\tau^1(-\iota_3)=-\tau^1.$$ 
Hence $\sigma_7\tau^1=-\sigma_3\tau^1=\tau^1$. 

Let $\eta\in\pi_3(S^2)\cong\Z$ be the generator (i.e. the class of
the Hopf map) and $h_0:\pi_3(S^2)\to\Z$ the Hopf isomorphism. We
have by [Po85, Complement to Lecture 6, (10)]
$$\sigma_1\tau^2=\sigma_1^2\tau^2=\tau^2(a_2\circ\eta)=
\tau^2(-\eta+[\iota_2,\iota_2]h_0(\eta))=\tau^2(-\eta+2\eta)=\tau^2.$$
Hence 
$$\sigma_7\tau^2=\sigma_1^1\tau^2=\sigma_1\tau^2=\tau^2
\quad\text{and}\quad\sigma_3\tau^2=-\sigma_7\tau^2=-\tau^2.$$
The relation for $\omega$ follows from the Symmetry Lemma of \S4. 
\qed

\smallskip
Note that $\omega_{0,AD}$ cannot be extended to a {\it smooth} embedding
$T^{1,2l-1}_+\to\R^{3l+1}_+$ without connected summation with $\varphi$.
Indeed, if there existed such an extension, then we could modify it by boundary surgery to a smoothly embedded $2l$-manifold bounding the connected sum of
Borromean rings.
Then surger certain circle in this manifold so as to obtain a smoothly embedded
$2l$-disk bounding the connected sum of Borromean rings.
The latter is impossible by [Ha62'].

\smallskip
(1) Describe possible homotopy types of complements to knotted tori.
E. g. for an embedding $f:T^{1,3}\to S^7$ the complement is
determined by an element of $\pi_4(S^2\vee S^3)\cong\Z_2\oplus\Z_2\oplus\Z$.
Describe possible normal bundles of knotted tori, cf. [MR71].
E. g. prove that
$$\nu(\tau(\varphi))\cong1\oplus\kappa^*\pi^*\varphi^*t,\quad\text{where}
\quad S^q\overset{\varphi}\to\to V_{m-q,p+1}\overset{\pi}\to\to G_{m-q,p+1}
\overset{\kappa}\to\cong G_{m-q,m-q-p-1}\to G_{m,m-p-1}$$
and $t$ is the tautological bundle. 
Describe possible immersion classes of knotted tori.

(2) Which values can assume the linking coefficient of an embedding
$T^{0,q}\to\R^m$, extendable to an embedding $T^{p,q}\to\R^m$?
More generally, using the classification of links [Ha66', Hab86] and
knotted tori, describe natural maps
$$KT^m_{p,q}\to KT^m_{p,q,+}\to KT^m_{p-1,q}\to KT^m_{0,q}\quad\text{and}\quad 
KT^m_{p,q}\to KT^{m+1}_{p,q}.$$
See [CR05].
Observe that the restriction
$KT^m_{p,q,CAT}\to KT^m_{p,q,CAT,+}\cong\pi_q(V^{CAT}_{m-q,p})$
is a complete invariant for regular homotopy.
For $p=1$ this invariant is the linking coefficient $\nu$ defined in \S5
[Sk02, Decomposition Lemma 7.1].

From the Filled-tori Theorem 7.3.a and Theorem 1.2 it follows that

{\it The restriction map $r^m_{p,q}:KT^m_{p,q,+}\to KT^m_{p-1,q}$
is a 1--1 correspondence if

$2m\ge 3q+2p+4$ in the PL category, or

$p-1\le q$ and $m\ge\max\{(3p+3q+1)/2,2p+q+2\}$ in the smooth category.}

Let us sketch an idea of a possible alternative direct proof of the first of these facts.
Let us prove the surjectivity of $r$ (the injectivity can be proved analogously).
Take an embedding $f:T^{p-1,q}\to\R^m$.
To each section of the normal bundle of $f$ there corresponds an embedding
$f':T^{p-1,q}\times I\to\R^m$.
The latter can be extended to an embedding $f'':(T^{p,q}_+-B^{p+q})\to\R^m$ by general
position.
We can always find a section such that the latter embedding extends to a map
$f''':T^{p,q}_+\to\R^m$ whose self-intersection set is contained in $B^{p+q}$.
The latter map is homotopic to an embedding relative to $T^{p,q}_+-B^{p+q}$ 
analogously to [Sk02, Theorem 2.2.q].

A {\it fiberwise embedding} $T^{p,q}\to\R^m$ is defined analogously
to that of $T^{p,q}_+\to\R^m$.
It would be interesting to find a direct proof of the following corollary 
of Filled-tori Theorem 7.3.d and Theorem 1.2. 

{\it For $m\ge\max\{2p+q+2,3q/2+p+2\}$ every PL embedding
$f:T^{p,q}\to\R^m$ is PL isotopic to a fiberwise PL embedding
and every PL isotopy between fiberwise PL embeddings is PL isotopic
relative to the ends to a fiberwise PL isotopy.}


(3) We conjecture that if
$$f,g:T^{p,q}\to S^m\quad\text{are embeddings},\quad M=S^m-R(f(S^p\vee S^q)),
\quad D^{p+q}=T^{p,q}-f^{-1}(\Int M),$$
$$\quad f=g\text{ on }f^{-1}R(f(S^p\vee S^q))\quad\text{and}
\quad g(D^{p+q})\subset M,\quad\text{then}\quad
[f|_{D^{p+q}}]=[g|_{D^{p+q}}]\in\pi_{p+q}(M,\partial M).$$

(4) Construct examples of non-surjectivity of
$\alpha^{3l+p+1}(T^{p,2l})$.

We conjecture that $\beta^{3l+p-1}(T^{p,2l-1}_+)$ is not injective
and so the dimenson restiction in the injectivity of [Sk02, Theorem
1.1.$\beta\partial$] is sharp.

(5) We conjecture that the exact sequence of the Main Theorem 1.4 in
the DIFF case splits also for $l$ even. The proof of this
conjecture requires the smooth case of the $\tau$-relation for $l$
even. A possible approach to prove this is to follow the proof of
the completeness of $\beta$-invariant [Ha84] and prove that for
our particular almost concordance
$\sigma_1\Omega_1\cup\overline\Omega_1$ we can obtain a {\it
smooth} concordance from $\sigma_1\omega_1$ to $\omega_1$. Or else
the required assertion is implied by the fact that every almost
smooth embedding $T^{1,2l}\to\R^{3l+2}$ is smoothable. We also
need to smoothen the relation
$\tau_p(w_{l,p})=\omega_p+\sigma_p\omega_p$.

Let us prove that {\it the map $\forg$ from the Smoothing Theorem 8.1
has a right inverse if $2m=3q+2p+3$, $q=2l-1$, $p<l$,
$l$ is even and $w_{l,p}$ is not divisible by 2.}


Indeed, by the AD case of the Main Theorem 1.4 we can identify $KT^{3l+p}_{p,2l-1,AD}$ with
\linebreak 
$\frac{\pi_{2l-1}(V_{l+p+1,p+1})\oplus\Z}{w_{l,p}\oplus2}$ via
the isomorphism $\widehat{\tau_{AD}\oplus\omega_{AD}}$.
Recall that $\tau(w_{l,p})=2\omega_p$ in the almost smooth category.
Then there is an element $x\in C_{2l+p-1}^{l+1}$ such that
$\tau(w_{l,p})-2\omega_p=\zeta(x)$ in the smooth category.
Since $w_{l,p}$ is not divisible by 2,
it follows that $w_{l,p}\oplus2$ is a primitive element.
Hence there is a homomorphism
$\varphi:\pi_{2l-1}(V_{l+p+1,p+1})\oplus\Z\to C_{2l+p-1}^{l+1}$ such that
$\varphi(w_{l,p}\oplus2)=x$.
Then $(\tau\oplus\omega)-\zeta\varphi$
factors through a right inverse of $\forg$.

(6) We conjecture that the map $\bar\alpha'\tau$ from Retsriction Lemma 5.2 
is not epimorphic and $\beta\equiv0$ when $w_{l,p}=0$ (in particular, for $l=15$).

(7) We conjecture that 
$\pi^{m-1}_{eq}((S^p\times S^q-\Int D^{p+q})\t{\ })\cong\Pi^{m-1}_{p-1,q}$ for 
$m\ge 2p+q+2$.

(8) We conjecture that {\it the forgetful map
$KT^{3l+p}_{p,2l-1,AD}\to KT^{3l+p}_{p,2l-1,PL}$ is an isomorphism for $0<p<l$}, or,
equivalently, that {\it the independence in the $\beta$-invariant Theorem 2.9 
hold in the PL
case}, or, equivalently, that {\it the Triviality Lemma 6.2 hold in the PL case}.
If this conjecture is true, then there is less indeterminancy in the Main Theorem 1.4 for
the PL case, and the Whitehead Torus Example 1.7 is true for embeddings 
$T^{p,2l-1}\to\R^{3l+p}$ and $l+1$ 
not a power of 2.


If $m\ge\frac{3q}2+p+2$, then by smoothing theory the forgetful map
$KT^m_{p,q,AD}\to KT^m_{p,q,PL}$ is an isomorphism  (moreover, every PL embedding
$T^{p,q}\to\R^m$ is almost smoothable and every PL isotopy between such embeddings is
almost smoothable {\it relative to the boundary}).

The forgetful homomorphism
$KT^{3l+p}_{p,2l-1,+,DIFF}\to KT^{3l+p}_{p,2l-1,+,PL}$ is an isomorphism by the assertion
on $\tau$ of the Torus Theorem 2.8 because $2(3l+p)=3(2l-1)+2(p-1)+5$.
But this does not imply the conjecture because passage from PL isotopy to a smooth one is 
made by a homotopy which could have self-intersections.
This also does not imply that
every PL (almost) concordance $F$ between smooth embeddings $T^{p,2l-1}_+\to\R^{3l+p}$ is
smoothable relative to the boundary, thus changing such an almost concordance to
a smooth almost concordance can possibly change $\beta$-invariant.
The obstruction to such a smoothing can be non-zero because this is an obstruction to
smoothing but not just to an existence of a smooth concordance between the same ends
(or, in other words, to extension of a smooth embedding to $T^{p,2l}_+$ from the boundary
not just existence of a smooth embedding).

The smoothing obstruction is that to extend smoothing given on a neighborhood of
$x\times\partial D^{2l}$, to $x\times D^{2l}$.
This obstruction splits into two parts.
The first one is that to extend smoothing from $x\times\partial D^{2l}$ to $x\times D^{2l}$
and is in $C_{2l-1}^{l+p+1}=0$ for $p>0$.
The second one is to extend a smooth $p$-frame from $x\times\partial D^{2l}$ to
$x\times D^{2l}$ and is in $\pi_{2l-1}(V^{PL}_{l+p+1,p},V_{l+p+1,p})$.
The stabilization of the second obstruction in $\pi_{2l-1}(V^{PL}_{l+M+1,M},V_{l+M+1,M})
\cong C_{2l-1}^{l+1}$ is perhaps the Haefliger smoothing obstruction [Ha].

The set of invariants of $\ker(KT^{3l+p}_{p,2l-1,PL}\to KT^{3l+p}_{p,2l-1,+,PL})$ depend on
$KT^{3l+p+1}_{p,2l,+,PL}$ and so it is hard to compare it with the same kernel in the almost
smooth category.

It would be interesting to understand how non-smoothable embeddings
$T^{p,2l,+}\to\R^{3l+p+1}$ are constructed from a knot $S^{2l-1}\to S^{3l}$ (presently they
are coming by the Miller isomorphism from homotopy groups of spheres),
and to prove that they extend to a PL embedding $T^{p,2l}\to\R^{3l+p+1}$.
In other words, it would be interesting to construct a map $\zeta'$ such that the sequence
$\Z_{(l)}\cong C^{l+1}_{2l-1}\overset{\zeta'}\to\to KT^{3l+p}_{p,2l,AD}
\to KT^{3l+p}_{p,2l,PL}$ is exact, cf. [CRS], and then to use this map to prove the
independence of $\beta$-invariant in the PL case for $p>1$.
Note that multiplication of such a knot with $S^p$ for $p<<l$ gives an obstruction in
$C^{l+1}_{2l-1+p}$ not in $C^{l+1}_{2l-1}$, while multiplication of such a knot with
$D^p$ for $p\ge l-2$ gives an embedding isotopic to the standard embedding.
For $p\ne1$ this map is defined using not $S^p$ but $p$-manifolds (or even $p$-cycles).


The forgetful homomorphism
$\overline{KT}^{3l+p}_{p,2l-1,AD}\to\overline{KT}^{3l+p}_{p,2l-1,PL}$ is an 
isomorphism by the assertion on $\overline\alpha$ of the Torus Theorem 2.8 because 
$2(3l+p)=3(2l-1)+2p+3$.
But this does not imply the conjecture because we cannot apply 5-lemma to 
the forgetful map of exact sequences 
$\overline{KT}^{3l+p+1}_{p,2l}\to\Z_{(l)}\to KT^{3l+p}_{p,2l-1}
\to\overline{KT}^{3l+p}_{p,2l-1}\to0$, 
since 
$\forg:\overline{KT}^{3l+p+1}_{p,2l,AD}\to\overline{KT}^{3l+p+1}_{p,2l,PL}$ 
is not proved to be epimorphic. 
We know that $\overline\alpha^{3l+p+1}_{p,2l,AD}$ is an isomorphism, but this 
does not imply that $(\overline\alpha_{AD})^{-1}\overline\alpha_{PL}$ is an 
inverse of $\forg$, so again we cannot conclude that $\beta_{PL}=0$.  

The generalization of the above conjecture for arbitrary $m,p,q$ is true if and only if the
inclusion-induced homomorphism $\pi_q(V_{m-q,p})\to\pi_q(V^{PL}_{m-q,p})$ is an isomorphism
(as in the two Haefliger's examples below).

Recall some known facts on the PL Stiefel manifolds.
For $p>q$ and $l>2$ by [Ha66, 4.8, Ha, 10.2, 11.2] there is an exact sequence
$$\dots\to\pi_{q+1}(V^{PL}_{p+l,p})\to C^l_q\to\pi_q(V_{p+l,p})
\overset{\rho^{PL}_q}\to\to\pi_q(V^{PL}_{p+l,p})\to\dots.$$ 
Note that such a sequence is not exact for $p=1$.
If $p\ge q$ and $l>2$ this and [Ha66] imply that {\it the inclusion-induced map
$\pi_q(V_{p+l,p})\to\pi_q(V^{PL}_{p+l,p})$ is an isomorphism for
$q\le2l-4$ and an epimorphism for $q=2l-3$.}

Denote $V=V_{l+p+1,p}$ and $V^{PL}=V_{l+p+1,p}^{PL}$.
Consider the following part of the above exact sequence for $p\ge2l$:
$$\pi_{2l+1}(V)\to\pi_{2l+1}(V^{PL})\to C_{2l}^{l+1}\to\pi_{2l}(V)\to
\pi_{2l}(V^{PL})\to C_{2l-1}^{l+1}\to\pi_{2l-1}(V)\to\pi_{2l-1}(V^{PL})\to0.$$
Then 

for $l$ odd $\rho^{PL}_{2l}$ is monomorphic and $\rho^{PL}_{2l+1}$ is 
epimorphic (because $C_{2l}^{l+1}=0$),

$\rho^{PL}_{2l}$ is not epimorphic for $l$ even and $p=2l$ [Ha67, \S5], and

$\rho^{PL}_{2l+1}$ is not monomorphic for $l$ odd and $p=2l+1$ [Ha67, \S5].


(9) Is it correct that $\beta(h)=\beta(f)+\beta(g,f(x\times S^{2l}))$ if the almost
embedding $h:T^{p,2l}\to\R^{3l+p}$ is obtained from an almost embedding
$f:T^{p,2l}\to\R^{3l+p}$ by connected summation with a map
$g:S^{p+2l}\to\R^{3l+p}-f(x\times S^{2l})$, when the images
of $f$ and $g$ are not assumed to be contained in disjoint balls?
U. Koschorke kindly informed me that the difference
$\beta(f)+\beta(g,f(x\times S^{2l}))-\beta(h)$ is formed by intersection in which we can
choose the first and the second sheet, so the classification space is $D^0$ not
$\R P^\infty$, thus such a difference is a simple invariant (but still it is not clear
whether it is trivial or not).

\smallskip
(10) Note that the Fibering Lemma 7.1 does not immediately follow from the 
'$n$-concordance implies $n$-isotopy' theorem  
because it is not clear why there exist balls 
$D^{m-q}\times D^q_\pm\subset\R^m$ such that 
$$f(D^p\times D^q_\pm)\subset D^{m-q}\times D^q_\pm\quad\text{and}\quad 
D^{m-q}\times D^q_+\cap D^{m-q}\times D^q_-=D^{m-q}\times \partial D^q_+=
D^{m-q}\times\partial D^q_-.$$
Such balls are constructed using the Slicing Lemma 7.2.

\smallskip
(11) {\bf Proof of the Slicing Lemma 7.2.} 

We present the proof from [Me02] with S. Melikhov's kind permission, because
the proof is less technical in our situation.
Consider the following statement. 

(k,l)\quad {\it Let $X$ be a finite simplicial $n$-complex, $m-n\ge3$ and 
$s:X\to\R^m$ a simplicial map for some triangulation $T$ of $\R^m$. 
Let $k=\dim s(X)$, let $C$ be the union of some top-dimensional dual cells 
in $\R^m$ and let $l$ be the number of top-dimensional dual cells in $\Cl(\R^m-C)$.
Then for each $\varepsilon>0$ there exists $\delta>0$ such that each PL embedding 
$f:X\to\R^m$ which is $\delta$-close to $s$ and such that 
$f^{-1}(\sigma)=s^{-1}(\sigma)$ for each dual cell $\sigma\subset C$, 
is PL $\varepsilon$-ambient isotopic 
keeping $C$ fixed to a PL embedding $h:X\to\R^m$ such that 
$h^{-1}(\sigma)=s^{-1}(\sigma)$ for each cell $\sigma$ of the dual to $T$ 
cell-subdivision of $\R^m$.} 

Then (i,0) and (0,j) are trivial for any $i$, $j$.
Assuming that (i,j) is proved for $i<k$ and arbitrary $j$, and also for $i=k$ 
and $j<l$, let us prove (k,l).
We shall reduce it to the following Edwards' lemma.
Consider the projections
$$\Pi_1:Q\times\R\to Q,\qquad\Pi_2:Q\times\R\to\R\qquad\text{and}\qquad
\pi:X\times[-1,1]\to[-1,1]\subset\R.$$

\smallskip
{\bf Edwards' Slicing Lemma.}
{\it For each $\varepsilon>0$ and a positive integer
$n$ there exists $\delta>0$ such that the following holds.

Let $X$ be a compact $n$-polyhedron and $Y$ its $(n-1)$-subpolyhedron,
$Q$ a PL $m$-manifold with boundary, $m-n\ge3$, and
$f:(X,Y)\times[-1,1]\to(Q,\partial Q)\times\R$ a PL embedding such that
$\Pi_2\circ f$ is $\delta$-close to $\pi$.
Then for each $\gamma>0$ there is a PL ambient isotopy $H_t$ with support in
$Q\times[-\varepsilon,\varepsilon]$ such that 

$H_t$ moves points less than $\varepsilon$, 

$\Pi_1\circ H_t$ moves points less than $\gamma$, 

$H_t$ takes $f$ onto a PL embedding $g$ such that 
$g^{-1}(Q\times J)=X\times(J\cap[-1,1])$ for each
$J=(-\infty,0]$, $\{0\}$, $[0,+\infty)$.

Moreover, if $f^{-1}(\partial Q\times J)=Y\times J$ for each $J$ as above, then
$H_t$ can be chosen to fix $\partial Q\times\R$} [Ed75, Lemma 4.1].

\smallskip
{\it Proof of (k,l).}
Choose any vertex $v\in\R^m-C$.
Let $D=\st(v,T')$ be its dual cell, and denote $E=\Cl(\partial D-\partial C)$.
Notice that the pair $(E,\partial E)$ is bi-collared in
$(\Cl(\R^m-C),\partial C)$.  By Edwards' Slicing Lemma, for
any $\varepsilon_1>0$ the number $\delta<\varepsilon_1$ can be chosen
so that $f$ is PL $\varepsilon_1$-ambient isotopic, keeping $C$ fixed,
to a PL embedding $h_1:X\to\R^m$ such that $h_1^{-1}(E)=s^{-1}(E)$.  It
follows that, in addition, $h_1^{-1}(D)=s^{-1}(D)$.  Note that $h_1$ is
$(\varepsilon_1+\delta)$-close to $s$ and
$\varepsilon_1+\delta<2\varepsilon_1$.

Consider the triangulation of $D$ defined from $T$ by a pseudo-radial
projection [RS72] $\partial D\to\partial\st(v,T)$
(note that in general $\partial\st(v,T)\ne\lk(v,T)$).
Apply (k-1,l') in $\partial D$ equipped with the above triangulation,
where $l'$ is the number of dual cells in $C\cap\partial D$.
Using collaring, we obtain that for each $\varepsilon_2>0$ we can choose
$\varepsilon_1<\varepsilon_2$ so that $h_1$ (which is $2\varepsilon_1$-close to $s$)
is PL $\varepsilon_2$-ambient isotopic, keeping $C$ fixed, to a PL embedding
$h_2:X\to\R^m$ such that $h_2^{-1}(\sigma)=s^{-1}(\sigma)$ for each
dual cell $\sigma$ of $C\cup D$.
Note that $h_2$ is $(\varepsilon_2+2\varepsilon_1)$-close to $s$ and
$\varepsilon_2+2\varepsilon_1<3\varepsilon_2$.

By (k,l-1), the number $\varepsilon_2<\varepsilon/3$ can be chosen so that
$h_2$ (which is $3\varepsilon_2$-close to $s$) is PL $\varepsilon/3$-ambient
isotopic, keeping $C\cup D$ fixed, to a PL embedding $h:X\to\R^m$ such
that $h^{-1}(\sigma)=s^{-1}(\sigma)$ for each dual cell $\sigma$ of
$T$.  Thus $f$ is $\varepsilon$-ambient isotopic to $h$, keeping $C$
fixed, because $\varepsilon/3+2\varepsilon_2<\varepsilon$.  \qed

\smallskip
Although the statement of the Fibering Lemma 7.1 does not involve
$\varepsilon$-control as in the Slicing Lemma 7.2, this control is used
in the proof.

We conjecture that the Fibering Lemma 7.1 is true even for TOP locally flat
embeddings.
Analogously to the Fibering Lemma 7.1 it is proved that if a PL manifold $N$ is
a $D^p$-bundle over another manifold, then every embedding $N\to\R^m$ in
codimension at least 3 is isotopic to a fiberwise embedding (which is defined
analogously).

\smallskip
(12) {\it An idea how to define $\sigma(f)$.} 
Let $N$ be a compact $n$-manifold.
 Fix an $n$-ball $B$ in the interior of $N$.
 Take a proper PL embedding $f:N\to D^m$ such that $f_{N-\Int B^n}$
 is a smooth embedding.
If we define a smooth regular neighborhood $B^m$ of $f(B^n)$ relative to
 $f(\partial B^n)$, we can define $\sigma(f)$ be the isotopy class of the 
abbreviation $\partial B^n\to\partial B^m$ of $f$ (and try to prove that 
this is independent on the choice of $B^m$). 


(13) {\it Sketch of a new proof that 
$KT^{3l}_{0,2l-1,PL}\cong\pi_{2l-1}(S^l)\oplus\Z_{(l)}$ for $l\ge2$.}
Analogously to the proof of Realization Theorem 2.3.b considering the following 
diagram: 
$$\minCDarrowwidth{0pt}\CD
0@>> > \Z_{(l)} @>>\omega> KT^{3l}_{0,2l-1} @>>\lambda> \pi_{2l-1}(S^l) @>> >0\\
@ .   @VV i V                 @AA\tau\oplus\omega A             @AA j A @ .\\
  @ .       @>> i > \pi_{2l-1}(S^l)\oplus\Z_{(l)} @>> j > @ .
\endCD$$
Here $i(x):=[(0,x)]$, $j[(a,b)]:=[a]$ (clearly, $j$ is
well-defined by this formula) and $\lambda$ is the linking coefficient.


\enddocument
\end

Idea of a possible proof of the commutativity of the right-upper square for $p\ge1$.
Recall that any $\varphi\in\pi_q(V_{m-q,p})$ can be represented as a trivial $p$-subbundle
of a trivial $(m-q)$-bundle over $S^q$.
In these terms $\lambda''(\varphi)$ is the first obstruction to existence of a non-zero
section in the orthogonal complement $\varphi^\perp$ to $\varphi$ (i.e. the framed Euler
class of $\varphi^\perp$).
Take a framing of $\varphi^\perp$ on $D^q_-$ and a general position section $g$ of
$\varphi^\perp$ which is non-zero on $D^q_+$.
Then $\pm\Sigma\lambda''(\varphi)$ goes under the Pontryagin isomorphism to the framed
submanifold $D^q_-\cap_{\varphi^\perp}gD^q_-\subset S^q$, where $D^q$ is the zero section
of $\varphi^\perp$.
Now take any general position extension $S^p\times S^q\to\R^m$ of the
embedding $\varphi:D^p_-\times S^q\to\R^m$.
The element $\pm\Sigma\overline\lambda\varphi$ goes under the Pontryagin isomorphism to the
framed submanifold
$D^p_-\times S^q\cap_{D^{m-q}\times S^q}S^q \subset D^p_-\times S^q\subset S^{p+q}$
(why $D^p_-\times S^q$ is framed???)
This is proved analogously to the above.
Hence by transversality we can 'cancel $D^p_-$' and obtain
$\pm\Sigma\overline\lambda\varphi=\Sigma^{p+1}\lambda''\varphi$.
Since $2m\ge3q+p+4$, it follows that $\pm\overline\lambda\varphi=\Sigma^p\lambda''\varphi$.

The idea of an alternative proof of the Triviality Lemma 6.2 for $l$ odd and $p\ge2$
(together with the application of the Hopf invariant) can perhaps
work also for $p=1$ because $\rho_n\oplus\mu'$ is an epimorphism for $q=2l$
and $l+1$ not a power of 2, even when $n=2$ (this is proved in the proof
of 2.7.$\rho''$ before the implication marked with '$n\ge3$').
The problem is that for $p=1$ the range of $\overline\alpha^{3l+p+1}_{p,2l}$ is not
known to be isomorphic to $\pi_{2l}(V^{eq}_{l+p+1,p+1})$.
However, we conjecture that $\alpha:KT^m_{p,q}\to\pi_q(V_{m-q,p+1})$ is
well-defined and bijective when $p=1$ and $m\le\frac{3q}2+2$
(i.e. when we do not know that $\sigma$ is an isomorphism).

Note that the analogous obstruction $\beta'(F|_{B^n})\in
H_{2n-m+1}(M,\partial M)$ to almost concordance (not $\rel\partial
B^n$) of $F|_{B^n}$ to an embedding into $M$ [Hu70', Ha84] is
zero. Indeed, consider the following segment of the exact sequence
of the pair $(M,\partial M)$:
$$H_{2n-m+1}(\partial M)\overset{i}\to\to H_{2n-m+1}(M)\overset{j}\to
\to H_{2n-m+1}(M,\partial M).$$
Since $F|_{\partial B^n}$ is an embedding, by definition of
$\beta'(F|_{B^n})$ it follows that $\beta'(F|_{B^n})\in\im j$.
From the Mayer-Vietoris sequence for $B^m=M\cup(B^m-\Int M)$ we obtain that the sum of
inclusions $H_k(\partial M)\to H_k(M)\oplus H_k(B^m-\Int M)$ is an isomorphism.
Hence the above map $i$ is a projection onto a direct summand.
So $j$ is the zero homomorphism and $\beta'(F|_{B^n})=0$.

\subhead Remarks to the proof of the Whitehead Torus Theorem \endsubhead

\proclaim{The Haefliger-Irwin-Zeeman Embedding Theorem 9.0} [Hae61, Irw65]
A map $f:N\to M$ of a $(2n-m)$-connected $n$-manifold $N$ to a
$(2n-m+1)$-connected $m$-manifold $M$ (with or without boundary) such that
$f|_{\partial N}$ is an embedding into $\partial M$,
is homotopic to an embedding, provided $m\ge n+3$ and $2m\ge3n+3$
in the PL and DIFF category, respectively.
\endproclaim

\proclaim{Mirror-Symmetry Lemma 9.2} Suppose that a smooth $n$-manifold $N$
has a smooth involution $\sigma$ whose fixed set is a smooth submanifold
$N_0\subset N$.

(a) If a PL embedding $f:N\to\R^m$ is mirror-symmetric with respect to
$\R^{m-1}\subset\R^m$ (i.e. $\sigma f=\sigma_mf$) and is smooth
outside the union of two $\sigma$-symmetric $n$-balls in $N-N_0$,
then $f$ is smoothable.

(b) If a PL isotopy $F:N\times I\to\R^m\times I$ is
mirror-symmetric with respect to $\R^{m-1}\times I\subset\R^m\times I$
and is smooth outside the union of two mirror-symmetric $(n+1)$-balls in
$(N-N_0)\times(0,1)$, then $F$ is smoothable.
\endproclaim

The inductive step in the construction of $\omega_{p,DIFF}$ follows by the
AD case of the Extension Lemma 9.1.b and Mirror-Symmetry Lemma 9.2 below.
This $\omega_{p,DIFF}$ is not necessarily unique or mirror-symmetric:
the smoothing of a PL embedding or isotopy obtained by
Mirror-Symmetry Lemma 9.2 is not necessarily mirror-symmetric
(and the non-uniqueness can appear in smoothing of $\omega_{p-1,DIFF}$,
not in isotopy of $\omega_{p,\pm,DIFF}$).

\demo{Proof of Mirror-Symmetry Lemma 9.2}
We prove only (a) because the proof of (b) is analogous.
There is a union $D_+^m\sqcup D_-^m\subset\R^m$ of two disjoint balls such that
$f|_{N-\Int D_+^n-\Int D_-^n}$ is
a proper smooth embedding into $\R^m-\Int D_+^m-\Int D_-^m$.
If we fix an orientation on $N$, then we fix orientations on $\partial D_\pm^n$,
i.e. diffeomorphisms from them to the standard sphere $S^{n-1}$.
Any orientation on $\R^m$ gives orientations on $\partial D_\pm^m$, i.e.
diffeomorphisms from them to the standard sphere $S^{m-1}$.
In this way the class $\eta_\pm\in C_{n-1}^{m-n}$ of
$f|_{\partial D_\pm^n}:\partial D_\pm^n\to\partial D_\pm^m$ is well-defined.
By [BH70, Hae], there is a complete obstruction
$\eta\in C_{n-1}^{m-n}$ to smoothing of $f$ and $\eta=\eta_++\eta_-$.
Since $f$ is mirror-symmetric with respect to $\R^{m-1}$,
by [Ha66, \S1] we have $\eta_-=-\eta_+$, so $\eta=\eta_++\eta_-=0$.
Therefore $f$ is smoothable.
\qed\enddemo

The equality $\eta_+=-\eta_-$ of Mirror-Symmetry Lemma 9.2
can be illustrated by the following example.
Since for $m-n\ge3$ every PL embedding $S^n\to\R^m$ is smoothable,
it follows that the obstruction to smoothing a suspension over a smooth knot
$\eta_+:S^{n-1}\to\R^{m-1}$  equals 0 but not $2\eta_+$. For $p=1$ and $l$
even the smoothability of $\omega_p$ also follows by [BH70, Bo71].

In the above proof the equality $\eta_+=\pm\eta_-$ does not
require orientations and is is simpler to understand than
$\eta_+=-\eta_-$. This weaker equality implies that
$\eta(\omega_p)$ is either 0 or $2\eta_+$. For $p=1$ and $l$ odd
this is sufficient because $C_{2l+p-2}^{l+1}\cong\Z_2$.

By Mirror-Symmetry Lemma 9.2 the elements of $\im\mu'_{p,AD}$ are
smoothable (the smoothing is not necessarily unique).
Note that $\mu'_{p,DIFF}$ is not surjective for $p=1$, $q=4k-1$ and $m=6k+1$
by the diagram in open problem (3).

\bigskip
The PL and AD cases of the 2-relation for $p>0$ can be proved
exactly as in the above proof for $p=0$, replacing $\omega_0$ to
$\omega_p$ and $\Omega_0$ to $\Omega_p$. Since the almost
concordance $(-\Omega_p)\cup\overline\Omega_p$ is smooth, the
obtained PL concordance between $-\omega_p$ and $\omega_p$ is
almost smooth. (Here we do not use the fact that $\omega_p$ is
mirror-symmetric, but even if we use it, the obtained almost
smooth concordance is not necessarily mirror-symmetric and so
cannot be smoothed by Mirror-Symmetry Lemma 9.2.)

Let us prove independently that $\omega_1=\sigma_{3l+1}\omega_1$.
Since $\sigma_{3l+1}$ reverses the orientation of $\R^{3l+1}\times I$
but not of $S^1\times S^{2l-1}\times I$,
it follows that $\sigma_{3l+1}$ reverses the orientation of
$[\Sigma(\Omega_1)]$ and $C$.
Since $\sigma_{3l+1}$ reverses the orientation of $\R^{3l+1}\times I$
and of $C$ but not of $x\times S^{2l-1}\times I$,
it follows that $\beta(\sigma_{3l+1}\Omega _1)=\beta(\Omega _1)$.
Hence by the completeness of $\beta$-invariant
$\sigma_{3l+1}\omega_1$ is almost smoothly isotopic to $\omega_1$.

Let us prove independently that $\omega_1=-\sigma_{2l-1}\omega_1$.
Since $\sigma_{2l-1}$ reverses the orientation of
$S^1\times S^{2l-1}\times I$ but not of $\R^{3l+1}\times I$,
it follows that $\sigma_{2l-1}$ preserves the orientation of
$[\Sigma(\Omega_1)]$ and $C$.
Since $\sigma_{2l-1}$ preserves the orientation of $\R^{3l+1}\times I$ and of $C$
but reverses orientation of $x\times S^{2l-1}\times I$,
it follows that $\beta(\sigma_{2l-1}\Omega _1)=-\beta(\Omega _1)$.
Hence by the completeness of $\beta$-invariant
$-\sigma_{2l-1}\omega_1$ is almost smoothly isotopic to $\omega_1$.

Since $\sigma_1$ reverses the orientation of $S^1\times S^{2l-1}\times I$
but not of $x\times S^{2l-1}\times I$ and not of $\R^{3l+1}\times I$,
it follows that $\beta(\sigma_1\Omega _1)=-\beta(\Omega _1)$.

Note that $\sigma_4\psi=\psi$ for each $\psi\in C_4^3$ can be proved directly
analogously to $\sigma_7\psi=-\psi$ (but not just deduced from this formula).
Indeed, change of orientation of $S^4$ induces change of
orientation of $S^2$ standardly linked with $\psi(S^4)$.
Therefore change of orientation of $S^4$ carries the attaching invariant
$a(\psi,\xi):S^4\times S^2\to S^2$ to
$\sigma_2\circ a(\psi,\xi)\circ(\sigma_4\times\sigma_2)$.

If $n=2$, $l$ is odd, $l+1$ is not a power of 2 and {\it $\pi_{l-1}^S$ has
no $\Z_2$-summands (this is so for $l\in\{5,13,21,45\}$ [To62])}, we can
prove that $w_{l,n-1}\ne0$ (which was used in the alternative proof above)
without use of [Os86] as follows.  Consider the last diagram for $q=2l$.
We have that $\rho_{n-1}$ is an isomorphism and
$\pi_q(V_{l+n,n-1})\cong\pi_q(V^{eq}_{l+n,n-1})\cong\pi_{q+n-2}(S^{l+n-1})
\cong\pi^S_{l-1}$ (identify these groups) .  The homomorphism
$\lambda'=\lambda''$ is the multiplication by 2 [JW54, \S9, Sk02,
Decomposition Lemma 7.1] and the right $\Sigma^{n-1}$ is an epimorphism
with kernel $[\iota_l,\iota_l]$.  Recall that for $l+1\ne2^s$, the element
$[\iota_l,\iota_l]$ is odd of order 2 [Co95].  If to the contrary
$w_{l,1}=0$, then there exists $x\in\pi_{l-1}^S$ such that
$\lambda''x=[\iota_l,\iota_l]$.  Then $2x=\Sigma\lambda''x=0$, so $x$ is
of order 2.  Since $[\iota_l,\iota_l]=\lambda''x$ is odd, it follows that
$x$ is odd.  Hence $\pi_{l-1}^S$ has a $\Z_2$-summand, which is a
contradiction.

$$\minCDarrowwidth{5pt}\CD
F @>> > V^{eq/diff}_{l+n,n}@>>  > V^{eq/diff}_{l+n,n-1}\\
@VV V @VV V @VV V \\
S^l @>> > V_{l+n,n}@>> r > V_{l+n,n-1}\\
@VV V @VV \rho_n V @VV \rho_{n-1} V \\
\Omega_{n-1}S^{l+n-1} @>> > V^{eq}_{l+n,n} @>> r_{eq} > V^{eq}_{l+n,n-1}
\endCD.$$
Let us present an idea of how to prove or disprove Homotopy-Theoretical
Lemma 7.3.c for any $l\not\in\{1,3,7\}$.
Consider the above diagram.
There the upper line consists of homotopy fibers of the inclusions of the
spaces in the second line into those in the third line.
Then the upper line is conjecturally a fibration.
Consider a map of exact sequences of fibrations:
$$\minCDarrowwidth{5pt}\CD \dots @>> >
\pi_{q+1}(V^{eq/diff}_{l+n,n}) @>> >
\pi_{q+1}(V^{eq/diff}_{l+n,n-1}) @>> >
\pi_{q-l}(V_{l+n-1,n-1}) @>> >
\pi_q(V^{eq/diff}_{l+n,n}) @ <\beta << \\
@VV V  @VV V  @VV V  @VV V  @VV V \\
\pi_q(S^l)  @>> > \pi_q(V_{l+n,n})
@>> > \pi_q(V_{l+n,n-1}) @>> >
\pi_{q-1}(S^l) @>> >
\pi_{q-1}(V_{l+n,n}) \\
@VV \Sigma^{n-1} V @VV \rho_n V @VV \rho_{n-1} V
@VV \Sigma^{n-1} V @VV \rho_n V \\
\pi_{q+n-1}(S^{l+n-1})@>> > \pi_q(V^{eq}_{l+n,n})@>> >
\pi_q(V^{eq}_{l+n,n-1}) @>> > \pi_{q+n-2}(S^{l+n-1}) @>> >
\pi_{q-1}(V^{eq}_{l+n,n})\\
@VV H V @VV \beta V @VV . V @VV . V @VV . V
 \endCD.$$
Note that the fourth column
of this diagram coincides with the first one.

\bigskip

It follows from the Torus Theorem 2.8 and from existence of a
map, analogous to $\tau$, for $p$ interchanged with $q$, that
$\im\alpha_{DIFF}\supset\pi_q(V_{m-q,p+1})\vee\pi_p(V_{m-p,q+1})$ for
$m\ge\frac{3q}2+p+2$.

\demo{Alternative proof of 2.7.$\bar\alpha$} The surjectivity of
$\overline\alpha$ follows from the existence of the map $\tau$
(even the surjectivity of $\overline\alpha\circ q$ follows). In
order to prove the injectivity take almost embeddings
$f,f_0:\QEmb^{3l+p}_{DIFF}(S^p\times S^{2l-1})$ for which there
exists an equivariant homotopy $\varphi$ between $\t f$ and
$\t{f_0}$. Take the map $p:(S^p\times S^{2l-1}\times I)\t{ }\to
\Sigma((S^p\times S^{2l-1})\t{ }\times I)$ constructed in [Sk02,
Cylinder Lemma 5.1]. Since $2(3l+p)\ge???3(2l-1+p)+1-(p-1)$, by a
boundary version of [Sk02, Theorem 1.1$\alpha\partial$] there
exists a smooth concordance (regular homotopy?) $F:D^p\times
S^{2l-1}\times I\to\R^{3l+p}\times I$ such that $\t
F\simeq_{eq}\Sigma\varphi\circ p$ on $(D^p\times S^{2l-1}\times
I)\t{ }$ ??? $F|_{0\times S^{2l}\times I}$ will have
self-intersections. Since the map $\Sigma\varphi\circ p$ is
defined on $(S^p\times S^{2l-1}\times I)\t{ }$, it follows that
$\t F$ extends over $(S^p\times S^{2l-1}\times I)\t{ }$.

Take any general position PL extension of $F$ to $S^p\times
S^{2l-1}\times I$. We this extension to a PL almost isotopy modulo
$D^p\times S^{2l-1}\times I$. This is done as in [Sk02, Proof of
the almost injectivity in Theorem 2.2.$q/\beta$, cf. Theorem
2.3.q] using a {\it relative} version of [Sk02, Disjunction
Theorem 3.1]. In this relative version a subcomplex $Z$ of $T$ is
given and the sentence 'Suppose that ... and $q\le m-2$' of
[Sk02, Disjunction Theorem 3.1] is replaced by the following
(because $2(3l+p)\ge3(2l+p+1)+2-2(p-1)$ for $p\ge1$):

{\it Suppose that for each simplices $\sigma^p,\tau^q,\nu^n\in T$ such that
$\sigma\times\tau\subset\Cl(E_1-E_0)$ and $\nu\times\tau\subset E_1$ we have
$p+q+n\le 2m-3$ and $\max\{p,q\}\le m-2$.
If $\sigma^p\subset Z$, then $\tau^q\not\subset Z$, $2p+q\le2m-3$ and
$q\le m-2$. And the obtained homotopy is $\rel Z$}

The proof of this relative version is analogous to [Sk02, Proof of
Disjunction Theorem 3.1].
Only the following changes are necessary.
We assume that $p\le q$ only if neither $\sigma\subset Z$ nor $\tau\subset Z$.
Otherwise we assume that $\sigma\subset Z$ and $\tau\not\subset Z$
(so $p\ge q$ is possible).
But '$\sigma^n\times\tau^q\subset Z$ hence $2n+q\le2m-3$ is required',
so this argument would not go.

We need to consider an isotopy on $D^p\times S^{2l-1}$ as $x\pi_{2l-1}(V)$
and choose it so that the deleted product extends to
$\t{S^p\times S^{2l-1}}$.

The obtained PL almost isotopy is smooth on $D^p\times S^{2l-1}$.
So by smoothing theory we may assume that it is smooth.
\qed\enddemo

Take any almost smooth embedding $f:T^{1,q}\to\R^m$.
If $\nu(f)=0$, analogously to Web Lemma 2.1 we can prove that
there is a web, i.e. a PL embedding $D^{q+1}\subset\R^m$ such that
$D^{q+1}\cap gT^{1,q}=\partial D^{q+1}=g(1\times S^q)$.
No, only for $2m\ge3q+6$!!!
Note that this method cannot be applied to prove that
a PL concordance $F$ between $2\tau^2$ and $2\omega_1$ has a web (which is of course false).
Indeed, for such a concordance $F(x\times S^3\times I)$ and $F(y\times S^3\times I)$
can be linked in $\R^7\times I$ modulo the boundary (even though
$F(x\times S^3\times I)$ and $F(y\times S^3\times t)$ are unlinked because
$F(x\times S^3\times0)$ and $F(y\times S^3\times0)$ in $\R^7\times0$ are).
Note that this method cannot be applied to prove that
any PL isotopy $F$ between smooth embeddings is smoothable.
Indeed, besides the above problem an isotopy from $F_0$ to a mirror-symmetric embedding
is not necessarily smooth.

By [To62, Theorem 10.3 and (10.10)']
$[\iota_{15},\iota_{15}]=2\sigma^2_{15}$ 
In general, we know $V_{l+p+1,p+1}=P^{l+p}_l\cup e^{2l+1}\cup\dots$,
where $P^{l+p}_l =P^{l+p}/P^{l-1}$ and
$P^m$ is the real $m$-dimensinal projective space.
So $\pi_{2l-1}(V_{l+p+1,p+1}) = \pi_{2l-1}(P^{l+p}_l)$.
Let $i: S^{15}\to M^{16}$,  $i': M^{16}\to P^{15+p}_{15}$ and
$i'':P^{15+p}_{15}->V_{15+p,p}$ be the canonical inclusion maps.
Then we obtain
$\mu''[\iota_{15}, \iota_{15}] = i\circ 2\iota_{15}\circ \sigma^2_{15} = 0$.
Hence we get that
$w_{15,p} = i''_*(i'_*(\mu''[\iota_{15}, \iota_{15}]))) = 0$.)

Denote by $\Qua$ the set of smooth almost concordance classes of
smooth almost concordances $F:N\to B^m$ between standard
embeddings (with fixed $F(x\times S^{2l-1}\times I)$ disjoint with
$\Sigma(F)$). By the invariance $\beta:\Qua\to\Z_{(l)}$ is
well-defined. Take classes $[F_1],[F_2]\in\Qua$ represented by
maps $F_1,F_2$ which are constant on some neighborhood of the
boundary in $S^p\times S^{2l-1}\times I$. Define
$[F_1]+[F_2]:=[F_1\cup F_2]$, $0\in\Qua$ to be the constant
homotopy and $-[F_1]=[\overline F]$. Clearly, this operation on
$\Qua$ is well-defined and endows $\Qua$ with a group structure
such that $\beta:\Qua\to\Z$ is a homomorphism. This operation
corresponds to the ordinary sum operation on almost embeddings.

\demo{Fourth proof of the DIFF independence}
Consider $S^{m-1}\times I$ as a subset of $B^m$.
Take any extension $F':N\to B^m$ of $F|_{\partial N}$.
Let $X=S^p\times S^{2l-1}$ and $Y=\t N\cap\partial(N\times N)$.
We assume that $B^n=B^{n-1}\times I$ for some PL ball $B^{n-1}\subset X$.
First we present a proof which work for $p\ge2$ and then
we modify it to work for the general case.
Consider the upper line of the following diagram.
$$\minCDarrowwidth{0pt}\CD
\Qua @>>\alpha> \pi^{m-1}_{eq}(\t N,Y,\t{F'}|_Y)
@>>\cong> \pi^{m-2}_{eq}(\t X\times I,\t X\times\{0,1\},\t{F|_{\partial N}})
@>>\cong> \pi^{m-3}_{eq}(\t X)\\
@VV V   @VV V    @VV V   @VV V \\
\Emb @>>\alpha> \pi^{m-1}_{eq}(\t{N_0},Y_0,\t{F'}|_{Y_0})
@>>\cong> \pi^{m-2}_{eq}(\t X_0\times I,\t
X_0\times\{0,1\},\t{F|_{\partial N}}) @>>\cong> \pi^{m-3}_{eq}(\t
X_0) \endCD$$ Analogously to the injectivity in [HH62, Sk02,
Theorem 1.1.$\beta\partial$] and almost injectivity in [Sk02,
Theorem 2.3.q$/\beta$] it is proved that $\alpha^m(N)$ is almost
injective for
$$2(3l+p+1)\ge3(2l+p)+1-(p-1)\quad\text{and}
\quad 2(3l+p+1)\ge3(2l+p)+2-2(p-1).$$ These inequalities are fulfilled for
$p>1$.  Therefore the upper $\alpha$ in the diagram is injective.  The
first upper isomorphism follows by the relative equivariant Suspension
Theorem analogously to [Sk02, proof of Theorem 5.2.$\alpha$], because
$2(p+2l-1)+1\le2(p+3l-1)-2$.
Since $\t{F|_{\partial N}}$
is equivariantly homotopic to a mapping whose image is in $S^{m-3}$, the
second upper isomorphism analogously follows by the relative equivariant
Suspension Theorem, because $2(p+2l-1)\le2(p+3l-2)-2$.

By [Sk02, Lemmas 6.1 and 7.3.a]
$\pi^{m-3}_{eq}(\t X)\cong\Pi^{m-3}_{p,2l-1}$ is finite for $p\le l-2$, so
$\Qua$ is finite.
If $p=l-1$, then $m=4l$ and $m-2=p+(2l-1+2)$.
Then analogously to [Sk02, Torus Lemma 6.1] it
is proved that there exists a surjection
$\Pi^{4l-3}_{l-1,2l-1}\oplus\Pi^{4l-3}_{2l-1,l-1}\to\pi^{m-3}_{eq}(\t X)$.
The domain of this surjection is finite by [Sk02, Lemma 7.3.a], so
$\pi^{m-3}_{eq}(\t X)$ and then $\Qua$ are finite.

In general case the finiteness of $\Qua$ follows from (a) below.

{\it (a) The set $\Emb$ of smooth isotopy classes of smooth proper
embeddings $F:N_0\to B^m$, extending given $F|_{\partial N}$ and
extendable to a PL almost concordance $F:N\to B^m$, is finite.}

Denote $Y_0=\t{N_0}\cap\partial(N\times N)$ and consider
the entire diagram above, where the vertical homomorphisms are restrictions.
By the boundary versions of the injectivity of [HH62, Remark in \S5,
Ha63, 6.4, cf.  Sk02, Theorems 1.1.$\alpha\partial$ and
1.1.$\beta\partial$], it follows that the bottom $\alpha$ is in injective.
The dimension restrictions
$2(3l+p+1)\ge3(2l+p)+1-(p-1)$ and $2(3l+p+1)\ge3(2l+p)-3(p-1)$
are fulfilled even for $p=1$, on the contrast to the above particular case.
Observe that $$X_0\times I\subset N_0\quad\text{and}\quad N_0-X_0\times
I\subset\partial N, \quad\text{so}\quad\t{N_0}-\t{X_0\times
I}\subset\partial(N\times N).$$ Therefore in the bottom line $N_0$ and
$Y_0$ can be replaced by $X_0\times I$ and $\t{X_0\times
I}\cap\partial(N\times N)$, respectively.  Hence the bottom isomorphisms
follow analogously to the corresponding upper isomorphisms.
Since $\pi^{m-3}_{eq}(\t X)$ is finite, it follows that $\Emb$ is finite.
\qed\enddemo

\proclaim{Lemma 5.6} $\Emb^{3k+p}(S^p\times S^{2k-1})$ is infinite for
each odd $k$ and $1\le p\le k-2$. false, for $p=0$ -- definitely!
\endproclaim

\demo{Proof} Let $q=2k-1$ and $m=3k+p$, then $m-p-q-1=k$.

First consider the case when CAT = DIFF and $p\not\equiv2,3(4)$.
Take the standard embedding $g_k:S^p\times S^q\times D^{k+1}\to S^m$.
Let $C_k=S^m-g_k(S^p\times S^q\times D^{k+1})$.
Clearly, $C_0\simeq S^p\sqcup S^q$.
Therefore
$$C_k\simeq\Sigma^k(S^p\sqcup S^q)\simeq S^k\vee S^{k+p}\vee S^{k+q}$$ for
each $k>0$.
Suppose that $f:S^p\times S^q\times D^{k+1}\to S^m$ is a smooth embedding such
that there exists a homotopy equivalence
$h:S^m-f(S^p\times S^q\times\delet D^{k+1})\to C_k$ and a bundle equivalence
$b$ of the trivial $\R^k$-bundle over $S^p\times S^q$ and $\nu(f)$.
Define {\it the attaching map} $a(f,h,b)$ to be the homotopy class of the
composition
$$S^p\times S^q\times S^k\overset{Sb}\to\to\nu(f)
\subset S^m-f(S^p\times S^q\times\delet D^{k+1})\overset{h}\to\to
C_k.$$ For a map $x:S^p\times S^q\times S^k\to C_k$ let $\beta
x:S^p\times S^q\times S^{k+1}\to C_{k+1}$ be the suspension on
each fiber $*\times *\times S^{k+1}$ (we have
$C_{k+1}\simeq\Sigma C_k$). Clearly, $a(g_0,h,b)$ is the constant
on $S^p\times S^q\times*$ for some choice of $h$ and $b$.
Therefore $a(g_k,h',b')=\beta^ka(g_0,h,b)$ is null-homotopic on
$S^p\times S^q\times*$ for some choice of $h'$ and $b'$. Hence we
can apply Homotopy-theoretical Lemma 5.7 to obtain an infinite
number of maps $x:S^p\times S^q\times S^k\to C_k$. By the Browder
embedding theorem [Bro68, Theorem 1, Git71, Theorem 3.4], for
each such $x$ there exist $f,b$ and $h$ such that $a(f,b,h)\simeq
x$. Hence the number of triples $(f,h,b)$ is infinite. Clearly,
for given $f$ there is only a finite number of $h$. If
$p\not\equiv2,3(4)$, then the groups $\pi_p(\SO_k)$,
$\pi_q(\SO_k)$ and $\pi_{p+q}(\SO_k)$ are finite [FoFu89, \S25].
Therefore for given $f$ there is only a finite number of $b$, and
the lemma follows.

The obstruction to smoothing a PL isotopy between two PL
embeddings $S^p\times S^{2k-1}\to S^{3k+p}$ lie in $C^{k+1}_{2k-1+p}$.
Since this group is finite when $2k+p$ is not divisible by 4 [Ha66],
the PL case of Lemma 5.6 follows from the smooth case for $p\ne2\mod4$.

Let us present a modification of the above proof which works for both
CAT = DIFF and PL, and without the $p\not\equiv2,3(4)$ restriction.
For an embedding $f:S^p\times S^q\to S^m$ denote by $S\nu(f)$ and
$\Int D\nu(f)$ the boundary and the interior of a tubular neighborhood of
$f(S^p\times S^q)$.
By [RoSa68, Theorem 4 and Corollary 5.9, see also Fad65], $S\nu(f)$ has a
unique structure of $S^k$-fibre bundle over $S^p\times S^q$.
For each smooth embedding $f:S^p\times S^q\to S^m$ there exists a section
$s:S^p\times S^q\to S\nu(f)$.
In the PL case consider the set of embeddings $f$ having such a section.
Since $m\ge2p+4$, we can extend an embedding
$S^p\times*\subset S^p\times S^q\overset{s}\to\to S^m-\Int D\nu(f)$
to an embedding $D^{p+1}\to S^m-\Int D\nu(f)$, unique up to isotopy (note
that by the Irwin-Zeeman theorem, here we may assume $2m\ge3p+q+4$
instead of $m\ge2p+4$).  Let $C=S\nu(f|_{S^p\times*})\cup D^{p+1}$.
Since $k>p$, it follows that $\nu(f|_{S^p\times*})$ is trivial.
Therefore $C\simeq\frac{S^p\times S^k}{S^p\times*}\simeq S^k\vee S^{k+p}$.
Since $k\ge2$, it follows that $C$ is simply-connected.
Since $H_i(C)=\Z$ for $i=0,k+p,k+q$ and is zero for other $i$, by the
Alexander duality and the Hurewicz isomorphism, the inclusion
$C\to S^m-\Int D\nu(f)$ is a $(p+q)$-homotopy equivalence.
Let $\beta(f,s)\in[S^p\times S^q,S^k\vee S^{k+p}]$ be the homotopy class of
compositions
$S^p\times S^q\overset{s}\to\to S^m-\Int D\nu(f)\overset{h}\to\to C$,
(where $h$ is a $(p+q)$-homotopy inverse).
Clearly, for an embedding $f$ such that $a(f,b,h)=x$ from the above proof and
some section $s$ we have $x|_{S^p\times S^q\times*}=\beta(f,s)$.
Therefore the number of pairs $(f,s)$ is infinite.
Since $p<k$, $q=2k-1$ and $k$ is odd, it follows that
$\pi_p(S^k)=0$ and the groups $\pi_q(S^k)$ and $\pi_{p+q}(S^k)$ are finite.
Hence for given $f$ there is only a finite number of $s$, and the lemma
follows.
\qed\enddemo

\proclaim{Homotopy-theoretical Lemma 5.7} (false!)
Let $C_k=S^k\vee S^{k+p}\vee S^{k+q}$ and
$a:S^p\times S^q\times S^k\to C_k$ a map such that the restrictions
$a_{S^p\times*\times*}$ and $a_{*\times S^q\times*}$ are null-homotopic.
If $k$ is odd, $q=2k-1$ and $1\le p\le k-2$, then there is an infinite
number of maps $x:S^p\times S^q\times S^k\to C_k$ which are homotopically
distinct (even $x|_{S^p\times S^q\times*}$ are homotopically distinct) and
such that $\beta^Kx\simeq\beta^Ka$ for $K$ large
(recall that for a map $x:S^p\times S^q\times S^k\to C_k$ the map
$\beta x:S^p\times S^q\times S^{k+1}\to C_{k+1}$ is the suspension on each
fiber $*\times *\times S^{k+1}$.
\endproclaim

Let us illustrate the idea of proof of Lemma 5.7 by proving the following
simpler statement: the set of extensions $S^a\times S^b\to X$ of a given map
map $\psi:*\times S^b\to X$ (up to homotopy) is infinite, if $\pi_{a+b-1}(X)$
is finite and $\pi_a(X)$ has an element $y$ of infinite order.
Indeed, consider the maps $\psi_1,\psi_2,\dots:S^q\to X$ representing the
elements $y,2y,\dots$.
The obstruction to extending the map $\psi_i\vee\psi:S^a\vee S^b\to X$ to
$S^a\times S^b$ is $[\psi,iy]=i[\psi,y]$.
Since this obstruction assume values in a finite group $\pi_{a+b-1}(X)$, it
follows that $i[\psi,y]=0$ for an infinite number of indices $i$, and the
statement follows.

\demo{Proof of Lemma 5.7} Analogously to [MaRo86, \S3], the map
$\rho:[S^{p+q},C_k]\to[S^p\times S^q,C_k]$ is a monomorphism.
Since $q=2k-1$, it follows that
$$[S^{p+q},C_k]\cong\pi_{p+q}(S^k\vee S^{k+p})\supset
\pi_{p+q+1}(S^{2k+p})\cong\Z.$$
So there is $y\in[S^{p+q},C_k]$ such that all the elements $\rho(iy)$ are
distinct for $i\in\Z$.

Let $X=S^p\times S^q\times S^k/(S^p\vee S^q)\times*$.
Consider the standard cell decomposition
$X=e^0\cup e^{p+q}\cup e^k\cup e^{k+p}\cup e^{k+q}\cup e^{k+p+q}$ (by $e^l$
we denote {\it closed} cells).
The elements $y,2y,\dots$ define homotopically distinct maps
$x_1',x_2',\dots$ on $e^{p+q}$.  By the Borsuk Homotopy Extension Theorem
we may assume that $a=a'\circ\pr$ for the factor-projection
$\pr:S^p\times S^q\times S^k\to X$ and some map $a':X\to C_k$.
Extend our maps $x_i'$ onto $e^k\cup e^{k+p}\cup e^{k+q}$ as $a'$ and
denote the extensions again by $x_i'$.

If $a$ is a constant map, then $a'$ is also constant, and the
proof is analogous to the above proof of the simpler statement.
In general case, the obstruction to extending the map $x_i'$ over
$e^{p+q+k}$ is a homotopy operation
$$\theta:\pi_{p+q}(\cdot)\times\pi_k(\cdot)\times\pi_{p+k}(\cdot)
\times\pi_{q+k}(\cdot)\to\pi_{p+q+k-1}(\cdot).$$
By [Po85, Lectures 4 and 5], this operation corresponds to an element from a
homotopy group of a wedge of spheres, and therefore is a
superposition of Whitehead products.
Hence for some $s$ and each $i$ we have
$\theta(iy,z,t,u)=i^s\theta(y,z,t,u)$.
We have
$$\pi_{p+q+k-1}(C_k)\cong\pi_{p+3k-2}(S^k\vee S^{k+p}\vee S^{k+q})\cong
\pi_{p+3k-2}(S^k)\oplus\pi_{p+3k-2}(S^{k+p})\oplus\pi_{p+3k-2}(S^{k+q})\oplus$$
$$\oplus[\pi_k(S^k)\otimes\pi_{p+2k-1}(S^{k+p})]\oplus\dots\oplus
[\pi_{2k-1}(S^k)\otimes\pi_{k+p}(S^{k+p})]\oplus\pi_{p+3k-2}(S^{p+3k-2}),$$
where dots denote a finite group.
Since $p\le k-2$ and $k$ is odd, it follows that this group is INfinite.
If $l$ is its order, then $\theta(ily,z,t,u)=0$ for each $i$, so we may assume
that $y$ is such that $\theta(iy,z,t,u)=0$ for each $i$.
Therefore $x_i'$ extends to $e^{p+q+k}$ for each $i$.
Denote an extension again by $x_i'$ and let $x_i=x_i'\circ\pr$.
Then $x_i|_{S^p\times S^q\times*}=\rho(iy)$ are homotopically distinct.
On the other hand, it is easy to see that $\beta^Kx_i\simeq\beta^Ka$ on
$\pr^{-1}(e^{p+q}\cup e^k\cup e^{k+p}\cup e^{k+q})$.
Clearly, we can choose the extension $x_i'|_{e^{p+q+k}}$ so that
$\beta^Kx_i\simeq\beta^Ka$ on the whole $S^p\times S^q\times S^k$.
\qed\enddemo

Consider the following segment of exact cofibration sequence of the pair
$(S^p\times S^q,S^p\vee S^q)$:
$$[S^{p+1}\vee S^{q+1},C_k]\overset{\delta}\to\to[S^{p+q},C_k]\overset{\pi}\to
\to[S^p\times S^q,C_k]\to[S^p\vee S^q,C_k].$$
Since $q=2k-1$, it follows that
$$[S^{p+q},C_k]\cong\pi_{p+q}(S^k\vee S^{k+p})\supset
\pi_{p+q+1}(S^{2k+p})\cong\Z.$$
Since $q=2k-1$ and $p\le k-2$, it follows that
$$[S^{p+1}\vee S^{q+1},C_k]\cong\pi_{q+1}(S^k\vee S^{k+p})\cong
\pi_{q+1}(S^k)\oplus\pi_{q+1}(S^{k+p})$$ is a finite group.
Therefore $\im\delta=\ker\pi$ is finite.

For an embedding $f:S^p\times S^q\times D^{k+1}\to S^m$ denote by $a(f)$
the {\it attaching map}, i.e. the homotopy class of the composition
$$S^p\times S^q\times\partial D^{k+1}
\subset S^m-g(S^p\times S^q\times\delet D^{k+1})\overset{def}\to=C_k.$$

$$\minCDarrowwidth{8pt}\CD
\pi_{p+q}(S^k\vee S^{k+p}) @>> j_0 > [S^p\times S^q,S^k\vee S^{k+p}] @>> >
[S^p\vee S^q,S^k\vee S^{k+p}]\\
@VV i V @AA \beta A  @AA A \\
[\wedge(S^p,S^q,S^k),C_k] @>> j > [S^p\times S^q\times S^k,C_k] @>> >
[S^p\vee S^q\vee S^k,C_k]\\
@VV \Sigma^{K-k} V  @VV \Sigma^{K-k} V  @VV \Sigma^{K-k} V \\
[\wedge(S^p,S^q,S^K),C_K] @>> J > [S^p\times S^q\times S^K,C_K] @>> >
[S^p\vee S^q\vee S^K,C_K]
\endCD.$$

By [Mas90, Theorem 3.1 and Remark 1],
$$\wedge(S^p,S^q,S^k)\simeq S^{p+q}\vee S^{k+p}\vee S^{k+q}\vee S^{k+p+q}.$$
In the previous paragraph it was proved that $\pi_{p+q}(S^k\vee
S^{k+p})\cong\pi_{p+q}(C_k)$. Therefore on the above diagram we
have the composition $i$ of this isomorphism  and the inclusion
of a direct summand. Define the map $\beta$ to be the restriction
onto $S^p\times S^q\times0$. Clearly,
$\beta(a(f))\in\overline\beta(f|_{S^p\times S^q\times0})$ (the choice
of necessary section $s$ is obvious). Clearly, $\Sigma^{K-k}i=0$.
Hence for each element $b\in\pi_{p+q}(S^k\vee S^{k+p})$,
$$\Sigma^{K-k}(a(g)+ji(b))=\Sigma^{K-k}a(g)+\Sigma^{K-k}ji(b)=
\Sigma^{K-k}a(g)+J\Sigma^{K-k}i(b)=\Sigma^{K-k}a(g).$$

Then $\beta(a(f))=\beta(a(g)+ji(b))=\beta a(g)+j_0(b)$.

\end